\documentclass[11pt]{article}

\usepackage[latin1]{inputenc}
\usepackage[T1]{fontenc}
\usepackage{amsmath}
\usepackage{amsfonts}
\usepackage{amssymb}
\usepackage[margin=1in]{geometry}
\usepackage[url=false,backend=biber,style=trad-abbrv,doi=false,isbn=false]{biblatex}
\usepackage{graphicx}
\usepackage{xparse}
\usepackage[colorlinks=true,citecolor=magenta]{hyperref}
\usepackage[capitalize,nameinlink]{cleveref}
\usepackage{siunitx}
\usepackage[labelfont=bf]{caption}
\usepackage{setspace}
\DeclareFieldFormat{volume}{volume \textbf{#1}}
\DeclareFieldFormat[article]{volume}{\textbf{#1}}

\usepackage{amsthm}
\theoremstyle{plain}
\newtheorem{example}{Example}[section]

\newtheorem{remark}{Remark}[section]
\numberwithin{equation}{section}

\makeindex

\newcommand{\normal}{\mathcal N}
\newcommand{\real}{\mathbf R}
\DeclareDocumentCommand\abs{s m o}{\IfBooleanTF{#1}{|#2|}{\left|#2\right|}\IfNoValueF{#3}{_{#3}}}
\DeclareDocumentCommand\eucnorm{s m o}{\IfBooleanTF{#1}{|#2|}{\left|#2\right|}\IfNoValueF{#3}{_{#3}}}
\DeclareDocumentCommand\ip{s m m o}{\IfBooleanTF{#1}{( #2,#3 )}{\left(#2,#3 \right)}\IfNoValueF{#4}{_{#4}}}
\DeclareDocumentCommand\eip{s m m o}{\IfBooleanTF{#1}{\langle #2,#3 \rangle}{\left\langle #2,#3 \right\rangle}\IfNoValueF{#4}{_{#4}}}
\newcommand{\var}{{\rm Var}}
\newcommand{\expect}{\mathbf E}

\renewcommand{\d}{\mathrm d}
\newcommand{\e}{\mathrm e}

\renewcommand{\div}{{\rm div}}

\renewcommand{\rho}{\varrho}

\newcommand{\grad}{\nabla}

\newcommand{\A}{\mathcal{A}}

\renewcommand{\t}{\mathsf T}

\newcommand{\argmin}{\operatornamewithlimits{argmin}}
\newcommand{\xb}{m_{f}}
\crefname{figure}{Figure}{Figures}

\begin{document}

\title{Consensus-based Optimization and Ensemble Kalman Inversion for Global Optimization Problems with Constraints}

\author{%
    J. A. Carrillo\thanks{%
        Mathematical Institute,
        University of Oxford,
        \texttt{carrillo@maths.ox.ac.uk}
    } \and
    C. Totzeck\thanks{%
        University of Wuppertal,
        School of Mathematics and Natural Sciences,
        \texttt{totzeck@uni-wuppertal.de}
    } \and
    U. Vaes\thanks{%
        MATHERIALS team,
        Inria Paris,
        \texttt{urbain.vaes@inria.fr}
    }
}

\maketitle

\begin{abstract}
    We introduce a practical method for incorporating equality and inequality constraints in global optimization methods based on stochastic interacting particle systems,
    specifically consensus-based optimization (CBO) and ensemble Kalman inversion (EKI).
    Unlike other approaches in the literature,
    the method we propose does not constrain the dynamics to the feasible region of the state space at all times;
    the particles evolve in the full space,
    but are attracted towards the feasible set by means of a penalization term added to the objective function
    and, in the case of CBO, an additional relaxation drift.
    We study the properties of the method through the associated mean-field Fokker--Planck equation
    and demonstrate its performance in numerical experiments on several test problems.
\end{abstract}

\vspace*{12pt}


\section{Introduction}
Finding efficiently the global minimizer of a function in large dimensions is crucial in today's data science activity,
with applications in various fields of research ranging from
modern areas such as machine learning, neural networks, inverse problems, data assimilation and model parameter estimation
to more classical areas such as economics, computational biology, mathematical finance and statistical physics.
In the most general form, the global optimization task is simply
\begin{equation}
    \label{eq:optimization_problem}
    \argmin\limits_{x\in B} f(x)
\end{equation}
for a given objective function $f$ and state space $B\subset \real^d$.
This simple-to-describe problem is highly nontrivial for general nonconvex functions~$f$
due to the existence of usually many,
and in some cases even infinitely many, local minima.
Most of the approaches have followed metaheuristic arguments in order to find robust algorithms;
see for instance \cite{aarts1988simulated,back1997handbook,reeves2003genetic,blum2003metaheuristics,dorigo2005ant,kennedy2010particle} and the references therein.
Among global optimization methods,
several algorithms were inspired by physical or biological models using the idea of interacting particle/agent systems:
ant colony optimization \cite{ant}, artificial bee colony optimization \cite{ABC}, firefly optimization \cite{Firefly},
wind driven optimization \cite{WDO}, and particle swarm optimization \cite{PSO,PoliKennedyBlackwell}.
All of these methods have in common the use of interactions of particles with nonlocal terms,
the exchange of information among agents
and the use of noise terms to explore the landscape of the function $f$;
see~\cite{bianchi2009survey} for a review.
Another family of global optimization algorithms is known as simulated annealing methods~\cite{holley1988simulated,holley1989asymptotics}.
The simplest examples of simulated annealing methods use the cost function~$f$ as the potential in an overdamped Langevin diffusion,
which is solved with a time-dependent temperature~\cite{hwang1990large,marquez1997convergence,pelletier1998weak} -- a cooling schedule --
so that the unique probability measure in the kernel of the Fokker--Planck operator associated with the dynamics at time $t$ converges to a Dirac mass at the global minimizer in the limit as $t \to \infty$.
Despite their simplicity, these methods usually have a slow convergence rate.

Concerning machine learning applications and large scale optimization,
there have been a flurry of works in SGD methods, where stochasticity is again introduced as a way of improving exploration of the landscape.
An essential ingredient of SGD methods \cite{RM51,bottou1998online,bubeck2015convex} is the mini-batch strategy,
also introduced recently in \cite{JinLiLiu,JinLiLiu2} for interacting particle systems,
which reduces the computational cost of the iterations.
The SGD family of methods relies on the computation of gradients of functions in high-dimensional settings,
which may be computationally expensive or even impossible, if the loss function is nondifferentiable.
In addition, even when gradients can be calculated at a reasonable computational cost,
they may fail to provide useful information on the overall behavior of the loss function
if this function exhibits many small-scale oscillations.
Gradient-free methods, which rely only on evaluations of the loss function,
are therefore an attractive alternative to SGD methods in these settings.
In this paper, we consider two particular classes of methods belonging to this category
which received a lot of attention lately:
the Consensus-Based Optimization~(CBO) methods~\cite{pinnau2017consensus,carrillo2018analytical,carrillo2019consensus,carrillo2021consensus},
reviewed recently in~\cite{totzeck2021trends},
and methods based on the Ensemble Kalman Filter (EnKF)~\cite{Chen2012,MR3036193,MR3041539,SS17,2019arXiv190308866G,garbuno2020affine,Carrillo_2021},
which are mainly employed in the context of inverse problems
but have also proved useful for machine learning tasks~\cite{MR3998631}.
It is shown in~\cite{2021arXiv210403384D}
that gradient-free ensemble Kalman methods perform better than their gradient-based counterparts in noisy likelihood landscapes.
Within the class of ensemble Kalman methods,
we focus specifically on the Ensemble Kalman Inversion (EKI) gradient-free method,
in the form proposed in~\cite{SS17}.

While the EKI method as discussed in~\cref{sec:eki_with_constraints} relies on a nonlocal approximation of gradients,
and can be viewed as an exact preconditioned gradient descent in the case of a quadratic objective function~\cite{2019arXiv190308866G},
the central mechanism of CBO is a relaxation drift towards the weighted average
\begin{gather}
    \label{eq:weightedmean}
    \xb = \frac{\sum_{j=1}^J x^{(j)} \e^{-\alpha f(x^{(j)})} }{\sum_{j=1}^J \e^{-\alpha f(x^{(j)})}}
    =: \sum_{j=1}^J w(x^{(j)}) x^{(j)}.
\end{gather}
Here, $x^{(j)}$ are the positions of a number $J$ of explorers evolving in the landscape of the objective function $f$.
As a metaphor for the CBO method,
one may think of butterflies in the mountains able to communicate their height and position instantaneously by wireless communication
and to integrate this information in the form of the weighted average $\xb$,
towards which they move with the hope of eventually reaching the global minimum of the landscape given by the function $f$.

The motivation for the CBO reweighting $\e^{-\alpha f(x)}$,
which is the Gibbs distribution corresponding to~$f(x)$,
 comes from statistical mechanics:
the probability distribution function $\frac{1}{Z_{\alpha}}\e^{-\alpha f(x)}$,
where $Z_{\alpha}$ is the normalization constant,
is the unique invariant measure of overdamped Langevin dynamics
at temperature $\alpha^{-1}$ in an external potential given by the cost function~$f(x)$~\cite[Chapter 6]{pavliotis2011applied}.
For sufficiently large $\alpha$,
the weighted average~\eqref{eq:weightedmean} provides a good approximation of the particle position $x^{(j_*)}$ minimizing~$f(x^{(j)})$,
assuming this is unique.
Indeed, in this limit all the weights in~\eqref{eq:weightedmean} converge to 0,
except for the one associated with particle $j_*$, which converges to 1.

Both the CBO and EKI methods are based on a system of interacting stochastic differential equations,
and it is possible to show that their behavior is well described by nonlocal deterministic Fokker--Planck equations when the number of particles is sufficiently large.
Studying these continuous equations brings considerable insight,
and turns out to be much simpler than working with the particle systems.
The mean-field Fokker--Planck equation associated with CBO includes a weighted average which,
despite no longer being a finite sum as in~\eqref{eq:weightedmean},
can still be understood in the limit as $\alpha \to \infty$.
This is achieved through Laplace's method~\cite{miller,dembo2009large},
a classical asymptotic method for integrals.

In its original form~\cite{pinnau2017consensus},
the CBO method combines the relaxation drift towards the weighted average $\xb$ with isotropic multiplicative noise
proportional in amplitude to the Euclidean distance to the weighted average~$\xb$.
The specific form of the multiplicative factor of the noise,
which may be viewed as a temperature scheduling,
is improved in~\cite{carrillo2021consensus} in order to decrease the computational cost of the method for large dimensional problems;
see \cref{sec:cbo_for_constrained_global_optimization_problems} for details.
In~\cite{carrillo2021consensus},
the authors also demonstrate that approximating $\xb$ based on a randomly selected subset of the particles
-- a mini-batch~--~leads to large computational savings.
Using mini-batches also introduces extra stochasticity,
which is observed to be useful for promoting exploration of the landscape given by $f$.

In this work, we propose new extensions of the CBO and EKI methods to constrained optimization and constrained inverse problems, respectively.
In the case of CBO, we consider problems of the general form
\begin{equation}
    \label{eq:our_problem}
    \argmin_{x \in B} f(x),
    \quad  \text{with} \quad B = \bigl\{ x : \mathcal E(x)=0, \, \mathcal I(x) \geq 0 \bigr\},
\end{equation}
where $\mathcal E\colon \real^d \to \real^{N_e}$ and $\mathcal I\colon \real^d \to \real^{N_i}$ are continuously differentiable functions
and the inequality $\mathcal I(x) \geq 0$ is understood componentwise.
In the case of EKI,
we consider subsets of problems of the form~\eqref{eq:our_problem} for which the objective function
is of the specific form $f(x) = \frac{1}{2}\bigl\lvert \Gamma^{-1} \bigl(G(x) - y\bigr) \bigr\rvert^2$,
where $G: \real^d \to \real^K$ is a map,
$\Gamma \in \real^{K \times K}$ is a positive definite matrix,
and $y \in \real^K$.
In applications, objective functions of this form are used as a simple measure of the misfit between some data $y$ and a model $G(x)$.

The literature on global optimization with constraints is abundant and an extensive review is beyond the scope of this paper,
so we summarize hereafter only recent contributions that are specifically related to the CBO and~EKI methods.
Several recent works propose extensions of CBO strategies to equality constrained problems,
based on imposing that the dynamics is restricted to the feasible manifold at all times~\cite{fornasier2020consensus,FornasierEtAl2,fornasier2021consensusbased}.
Most of these works consider the specific case of the Euclidean sphere as the feasible region,
with applications to robust subspace detection and efficient eigenvalue computations in machine learning as prime examples.
The evolution of particles is restricted to the feasible manifold by appropriate projection onto the tangent space.
Although this method performs relatively well in very simple settings like the sphere,
its implementation is difficult and error-prone for more general constraint manifolds.
In addition, since the weighted average is computed in the ambient space,
there are cases in which the method does not produce good results,
for example when the feasible region is a closed hypersurface enclosing a nonconvex domain.
Let us finally emphasize that manifolds other than the sphere are important in optimization problems in machine learning. SGD methods have been studied in Riemannian optimization applications,
see for instance \cite{zhang2016riemannian,weber2021projectionfree} and the references therein.
Here, the state space is usually a set of matrices with certain constraints,
like positive definite matrices in Wasserstein barycenters or Riemannian centroids.

Several methods have also been proposed in the literature for incorporating constraints in EKI.
In~\cite{MR3998633}, the authors propose a generalization of EKI such that the iterates produced by the method are guaranteed to lie in the feasible region at each iteration.
The method proposed leverages the fact that the update step in the usual ensemble Kalman method can be formulated as an optimization problem,
in which linear constraints can be integrated.
In~\cite{2019arXiv190800696C}, a variant of EKI incorporating a projection step is developed for the specific case of box constraints,
and the continuous time limit of the method is studied.
See also~\cite{MR4121318},
where the continuous-time and mean field limits of the method proposed in~\cite{MR3998633} are studied.

In this work, we follow an orthogonal approach:
we propose extensions of the CBO and EKI methods in which the particles evolve in the full space~$\real^d$,
instead of being constrained to the feasible manifold at all times.
A requirement for this approach to work,
which we take as a standing assumption in the rest of this paper,
is that the objective function~$f(x)$ can be evaluated anywhere in~$\real^d$,
which may or may not hold in applications.
The unifying idea behind our extensions of CBO and EKI is the addition of a penalization to the objective function.
In the case of only one equality constraint $\mathcal E(x) = 0$, for example,
we seek a solution to
\begin{equation}
    \label{eq:our_problem_one_equality}
    \argmin_{x \in \real^d} \left( f(x) + \frac{1}{\nu} \abs{\mathcal E(x)}^2 \right).
\end{equation}
When the penalization parameter $\nu$ tends to 0,
we expect the solution of this problem to be close to the zero-level set of $\mathcal E$.
We emphasize that the use of penalty functions for constrained optimization problems is a standard idea
and can be employed, in principle, together with any method for unconstrained global optimization.
In this paper,
we demonstrate the effectiveness of this simple idea when used in conjunction with~CBO and~EKI,
and we propose specific approaches for its implementation.

In the specific case of CBO,
we also propose to use, in addition to a penalization of the form~\eqref{eq:our_problem},
an extra drift term based on the constraint
that imposes that particles are asymptotically confined to the manifold.
This idea is reminiscent of swarming problems where the Vicsek model with noise on the sphere can be retrieved from the Cucker-Smale model on the whole space in an appropriate limit~\cite{bostan2013asymptotic}.

The rest of this paper is organized as follows.
In \cref{sec:cbo_for_constrained_global_optimization_problems,sec:eki_with_constraints},
we explain our strategy for including constraints in the CBO and the EKI methods, respectively.
In \cref{sec:cbo_for_constrained_global_optimization_problems},
we also compare the solution to the constrained CBO method and associated the mean-field equations qualitatively for a toy problem in two dimensions,
in order to illustrate the method and gain a deeper understanding of its behavior in the many-particle limit.
In \cref{sec:numerical_results,sec:numerics_for_eki},
we investigate the potential of the proposed methods by means of numerical experiments in some reference examples.
To this end, we use typical benchmarks~\cite{benchmarks} in optimization.
\Cref{sec:conclusion} is reserved for conclusions and perspectives for future work.


\section{CBO for constrained global optimization problems}%
\label{sec:cbo_for_constrained_global_optimization_problems}
In order to illustrate our strategy for incorporating constraints into the CBO methods and its variants,
we consider the CBO scheme with component-wise Brownian motion proposed in~\cite{carrillo2019consensus},
which reads as
\begin{equation}
    \label{eq:cbo_without_constraint}
    \d x^{(j)}_t = - (x^{(j)}_t- \xb) \, \d t + \sqrt{2}\sigma (x^{(j)}_t - \xb) \circ \d W^{(j)}_j,
    \quad j=1, \dots,J,
\end{equation}
where $\{W^{(j)}\}_{1 \leq j \leq J}$ are independent Wiener processes in $\real^d$ and
$v_1 \circ v_2$, for vectors $v_1$ and $v_2$ in $\real^d$, is the Hadamart (component-wise) product,
that is to say $v_1 \circ v_2 = {\rm diag}(v_1) v_2$.
The particles are initialized independently according to some given probability density $\rho_0,$
and so the law of the initial ensemble is $\rho_0^{\otimes J}$.
The main ingredient of the method,
indeed the only mechanism by which the particles interact with the objective function,
is the weighted mean $\xb$ given by
\begin{equation}
    \label{eq:weighted_mean_finite}
    \xb(\mu^J_t) = \frac{\int x \e^{-\alpha f(x)} \mu^J_t(\d x)}{\int \e^{-\alpha f(x)} \mu^J_t(\d x)},
    \qquad \mu^J_t = \frac{1}{J} \sum_{j=1}^{J} \delta_{x^{(j)}_t},
\end{equation}
where $\alpha > 0$ is a parameter.
Equivalently,
\begin{align*}
    \xb(\mu^J_t) = \frac
    { \sum_{j=1}^J x^{(j)}_t \exp\bigl(-\alpha f(x^{(j)}_t)\bigr)}
    { \sum_{j=1}^J \exp\bigl(-\alpha f(x^{(j)}_t)\bigr)}.
\end{align*}
As mentioned in the introduction,
for fixed $t$ the weighted mean~$\xb(\mu^J_t)$ in~\eqref{eq:weighted_mean_finite} converges to the global minimizer of $f$ constrained to the support of the empirical measure~$\mu^J_t$,
provided this minimizer is unique,
in the limit as $\alpha$ tends to infinity.

Assume the problem we aim to solve is~\eqref{eq:optimization_problem}.
In order to illustrate the incorporation of equality and inequality constraints in a unified manner,
we introduce $\mathcal A\colon \real^d \to \real^{N_e + N_i}$ given by
\[
    \bigl(\mathcal A(x)\bigr)_i =
    \begin{cases}
        \bigl(\mathcal E(x)\bigr)_i & \text{ if $i \leq N_e$,} \\
        \min \bigl\{\bigl(\mathcal I(x)\bigr)_{i-N_e}, 0\bigr\}  & \text{ if $i > N_e$.} \\
    \end{cases}
\]
We can then write $B = \{x \in \real^d: \mathcal A(x) = 0\}$;
that is, we can assume without loss of generality that all the constraints are of equality type.
Notice that $\mathcal A$ may not be $C^1$ at the manifold even though $\mathcal E$ and $\mathcal I$ are,
but this is not an issue given that only $\grad \abs{\mathcal A}^2$ appears in the method we propose,
and the function~$x \mapsto \abs{\mathcal A(x)}^2$ is continuously differentiable.
To solve the optimization problem \eqref{eq:optimization_problem} with constraint $\mathcal A(x) = 0$,
we modify the scheme~\eqref{eq:cbo_without_constraint} in the following manner:
\begin{enumerate}
\item
    The first modification is a penalization of the objective function in the weighted average.
    More precisely, we substitute $\xb$ in~\eqref{eq:cbo_without_constraint} by
    \begin{equation}
        \label{eq:weight_w_constraint}
        m_g(\mu^J_t) = \frac{\int x \e^{-\alpha g(x)} \mu^J_t(\d x)}{\int \e^{-\alpha g(x)} \mu^J_t(\d x)},
        \qquad g(x) := f(x) + \frac{1}{\nu} \lvert \mathcal A(x) \rvert^2,
    \end{equation}
    where $\nu$ is a small parameter.
    By the Laplace principle,
    for fixed $t$ this new weighted average tends to the minimizer of $g$ (within the support of the empirical measure and assuming this minimizer is unique) in the limit as $\alpha \to \infty$,
    which is expected to lie close to the feasible set $\{x: \A(x) = 0\}$ when the parameter~$\nu$ is sufficiently small.

    \

\item
    The second modification is the introduction of a drift term
    that drives particles towards the constraint manifold.
    For $0< \varepsilon \ll 1$ we propose the dynamics
    \begin{equation}
        \label{eq:cbo_with_constraints}
        \begin{aligned}[b]
            \d x^{(j)}_t &= - \frac{1}{\varepsilon} (\nabla \lvert \mathcal A \rvert^2)\bigl(x^{(j)}_t\bigr) \, \d t - \left(x^{\left(j\right)}_t- m_g(\mu^J_t)  \right) \, \d t \\
                         &\qquad + \sqrt{2} \sigma \bigl(x^{(j)}_t - m_g(\mu^J_t)  \bigr) \circ \d W^{(j)}_t.
        \end{aligned}
    \end{equation}
\end{enumerate}
Both mechanisms are useful for driving the particles to the constraint manifold.
The first modification is straightforward to implement and does not generate stiffness of the stochastic differential equations driving the particles.
Consequently, these can be discretized in time with a step of the same order of magnitude as that used for unconstrained problems.
The second modification is helpful for ensuring that the particles move towards the weighted average in a manner independent of the other particles.
The additional drift pushes the particles to the manifold in an asymptotically orthogonal manner
(indeed $\grad \abs{\mathcal A}^2$ is orthogonal to the contour lines of~$\abs{\mathcal A}^2$),
giving us more function evaluations along the manifold.
Despite being detrimental for stability at the discrete level,
this additional drift is observed to be useful in our numerical experiments.

\begin{example}
For global optimization problems subject to just one inequality constraint,
    i.e. problems of the form
    \[
        \argmin_{x \in B} f(x) \quad \text{with} \quad B = \{ x : \mathcal I(x) \ge 0 \}
    \]
    with $\mathcal I\colon \real^d \to \real$,
    the above scheme reads
    \begin{align*}
        \d x^{(j)}_t
        &= - \frac{1}{\varepsilon} \nabla (\mathcal I^2)(x^{(j)}_t)  \chi_{(-\infty, 0)}\bigl(\mathcal I(x^{(j)}_t)\bigr) \, \d t
        - (x^{(j)}_t- m_g(\mu^J_t) ) \, \d t \\
        &\qquad+ \sqrt{2} \sigma (x^{(j)}_t - m_g(\mu^J_t) ) \circ \d W^{(j)}_t,
    \end{align*}
    where $\chi_{(-\infty, 0)}$ denotes the indicator function of the set $(-\infty, 0)$.
    The weighted average is given by
    \[
        m_g(\mu^J_t) = \frac
        {\displaystyle \sum_{j=1}^J x^{(j)}_t \exp\Bigl(-\alpha f(x^{(j)}_t)- \alpha \nu^{-1} \mathcal I(x^{(j)}_t)^2 \chi_{(-\infty, 0)}\bigl(\mathcal I(x^{(j)}_t)\bigr)\Bigr)}
        {\displaystyle \sum_{j=1}^J \exp\Bigl(-\alpha f(x^{(j)}_t)- \alpha \nu^{-1} \mathcal I(x^{(j)}_t)^2 \chi_{(-\infty, 0)}\bigl(\mathcal I(x^{(j)}_t)\bigr)\Bigr)}.
    \]
\end{example}

\subsection{Mean-field CBO}\label{sec:mfcbo}
It is possible to prove a propagation of chaos result for~\eqref{eq:cbo_without_constraint}; see~\cite{huang2021mean}.
In the many-particle limit $J \gg 1$, the law of $(x^{(1)}_t, \dotsc, x^{(J)}_t)$ approximately tensorises,
i.e. it holds approximately that
\begin{align}
    \label{eq:tensorization_of_law}
    \bigl(x^{(1)}_t, \dotsc, x^{(J)}_t\bigr) \sim \mu_t^{\otimes J},
\end{align}
where $\mu_t$ is a probability measure.
In addition, this measure satisfies the partial differential equation
\begin{equation}
    \label{eq:mean_field_cbo_without_constraints}
    \partial_t \mu = \grad \cdot \Bigl(\bigl(x- \xb(\mu)\bigr)\mu\Bigr) + \sigma^2 \sum_{i=1}^d \partial_{ii}\Bigl(\bigl(x-\xb(\mu)\bigr)_i^2 \mu\Bigr),
    \qquad \mu_0 = \rho_0,
\end{equation}
in the distributional sense.
In order to formally derive this equation from~\eqref{eq:tensorization_of_law},
one can take the limit $J \to \infty$ in~\eqref{eq:cbo_without_constraint} and use the law of large numbers to obtain
\begin{equation}
    \d x^{(1)}_t = - \bigl(x^{(1)}_t- \xb(\mu_t)\bigr) \, \d t + \sqrt{2}\sigma \bigl(x^{(1)}_t - \xb(\mu_t)\bigr) \circ \d W^{(1)}_t.
\end{equation}
The Fokker--Planck equation corresponding to this equation is then~\eqref{eq:mean_field_cbo_without_constraints}.
A similar formal argument for~\eqref{eq:cbo_with_constraints} leads to the following nonlocal Fokker--Planck equation
\begin{align}
\begin{split}
\label{eq:mean_field_cbo_with_constraints}
    \partial_t \mu
    &= \frac{1}{\varepsilon} \grad \cdot (\nabla \lvert \mathcal A\rvert^2(x) \mu)  + \div_x \Bigl(\bigl(x- m_g(\mu)\bigr)\mu\Bigr) \\
    & \qquad + \sigma^2 \sum_{i=1}^d \partial_{ii}\Bigl(\bigl(x-m_g(\mu)\bigr)_i^2 \mu\Bigr),
    \qquad \mu_0 = \rho_0.
\end{split}
\end{align}

\begin{remark}
    If $g$ in \eqref{eq:weight_w_constraint} satisfies the assumptions in \cite{carrillo2018analytical},
    then all the theoretical results proved there for the mean-field Fokker--Planck equation associated with unconstrained CBO
    apply mutatis mutandis to~\eqref{eq:mean_field_cbo_without_constraints} when $\varepsilon = \infty$,
    i.e.~with only the first modification -- the penalization.
\end{remark}

\subsection{Numerical experiments: mean-field versus particles and the penalizations}

We illustrate the CBO approach with and without constraints by comparing the solution to the mean-field PDE~\eqref{eq:mean_field_cbo_without_constraints}
with the solution to the particle system~\eqref{eq:cbo_without_constraint} in simple two dimensional problems.
\begin{figure}[ht!]
    \centering
    \includegraphics[scale=0.4]{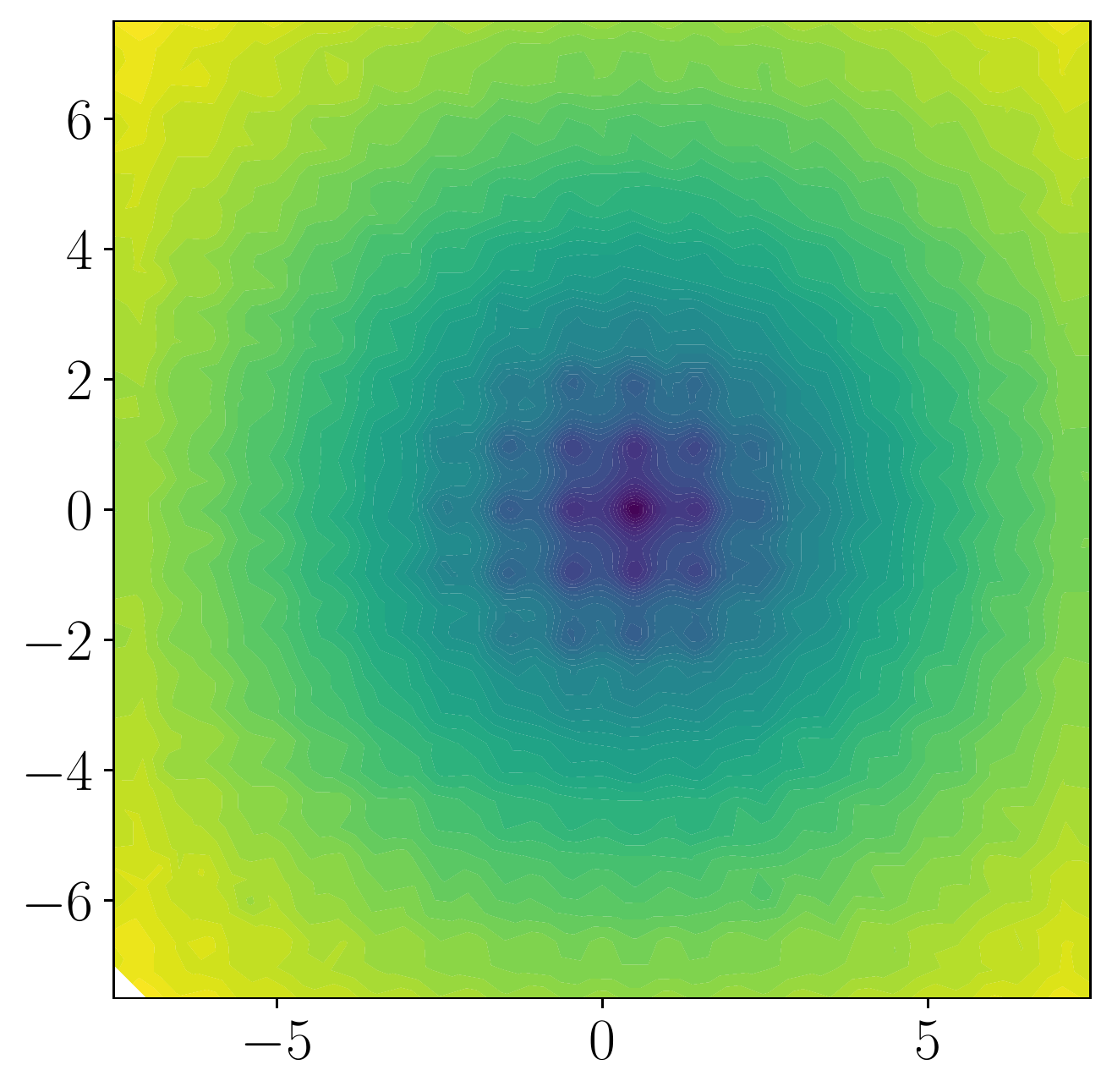}
    \includegraphics[scale=0.4]{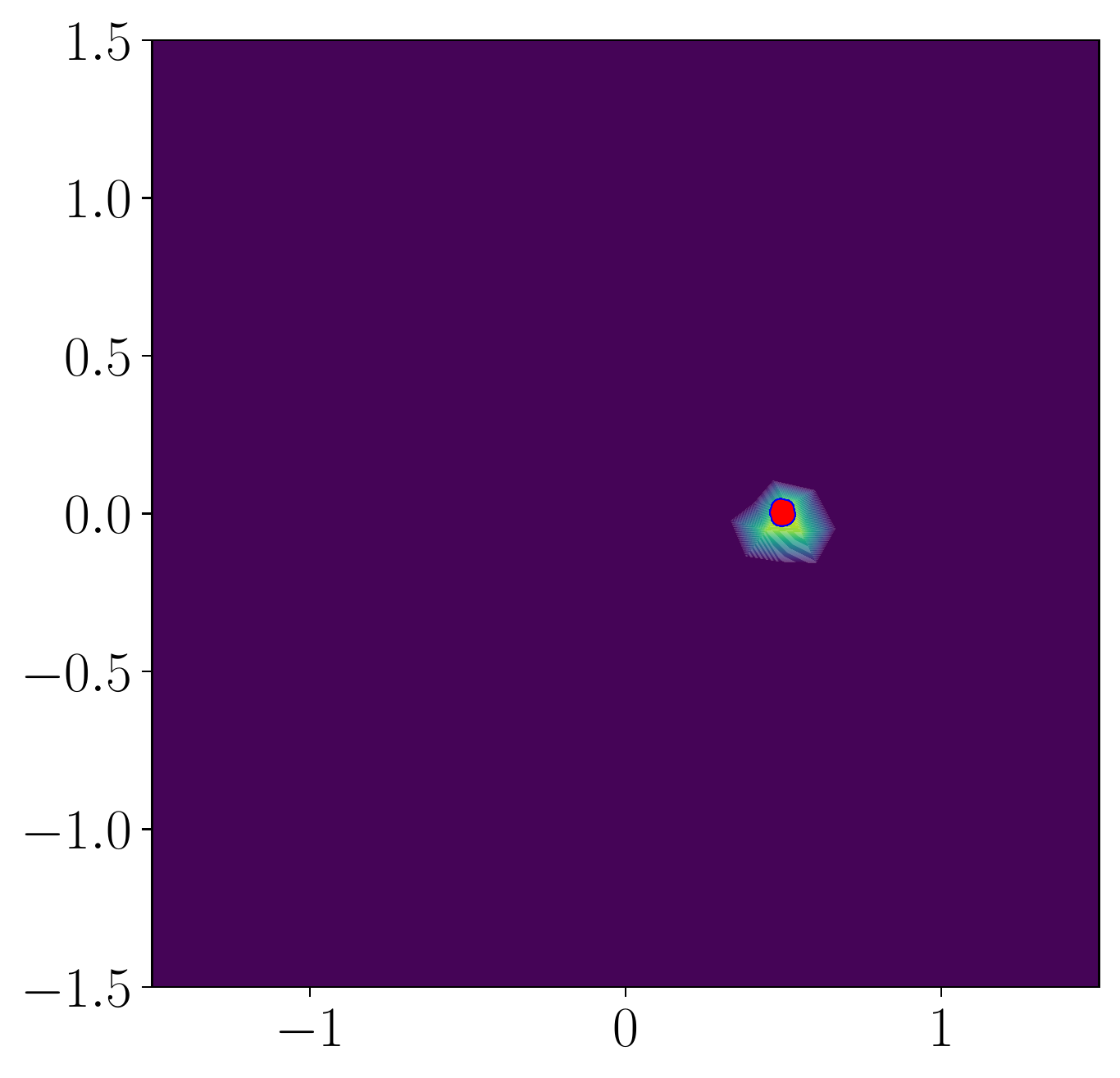}
    \caption{$J=50$, $M=100$, $\Delta t_{\rm SDE}=0.0005$, $\Delta t_{\rm PDE}=0.005$, $h_{\min}=0.125$, $h_{\max}=0.5$, $\alpha = 30$, $\nu = 1$, $\sigma = 0.7$, $T = 50.$
    Here $h_{\min}$ and $h_{\max}$ denote the minimum and maximum characteristic lengths of the mesh employed for the PDE simulation.
    Left: contour plot of the Ackley function shifted by $(0.5,0)$.
    Right: contour plot of the PDE solution which is concentrated at one triangle.
    The red points show the weighted averages of $100$ samples of particle simulations after convergence.}
    \label{fig:pde-sde-comparison}
\end{figure}

We begin with a comparison in the unconstrained setting.
To this end, we use the Ackley benchmark~\cite{benchmarks} shifted by $(0.5,0)$;
see~\cref{fig:pde-sde-comparison} (left) for a contour plot.
The SDE simulation uses $50$ particles,
which are drawn independently from the normal distribution $\mathcal N(0, 3 I_2)$,
where $I_2$ is the 2 by 2 identity matrix,
at the beginning of the simulation.
The PDE solution is initialized according to the same distribution.
In view of the shift,
the global minimizer of the objective function is not at the center of mass of the initial distribution.
The mesh for the PDE solver has a maximal meshfield size of~$h_{\max} = 0.5$ and is refined around the origin with a minimal meshfield size of~$h_{\min} = 0.125.$
In particular, the region around the global minimum is contained in the refined area.
Additional details on the numerical implementation and the parameters employed in the simulations are given in \cref{sec:numerical_results}.

\begin{figure}[ht!]
    \centering
    \includegraphics[scale=0.4]{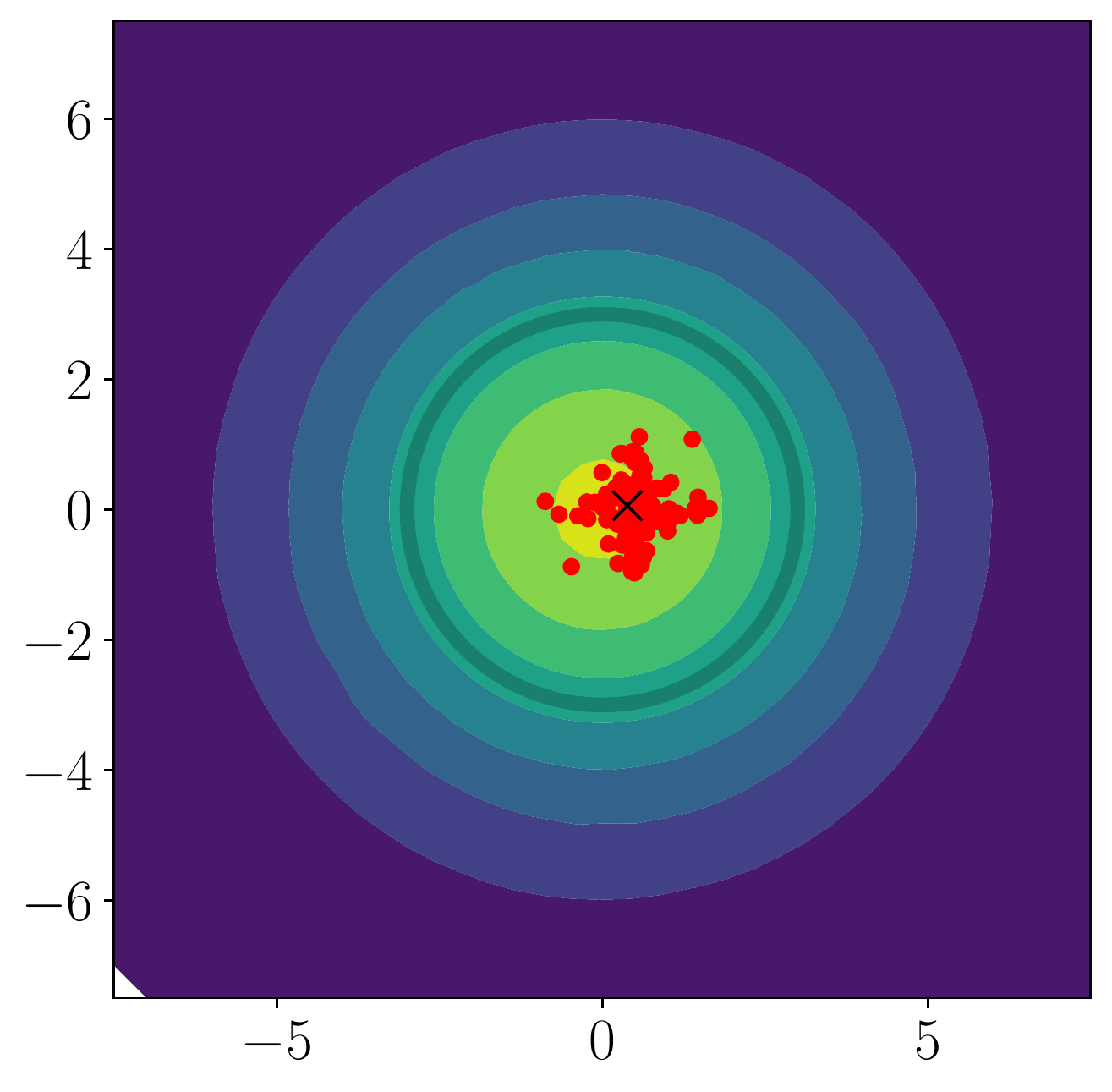}
    \includegraphics[scale=0.4]{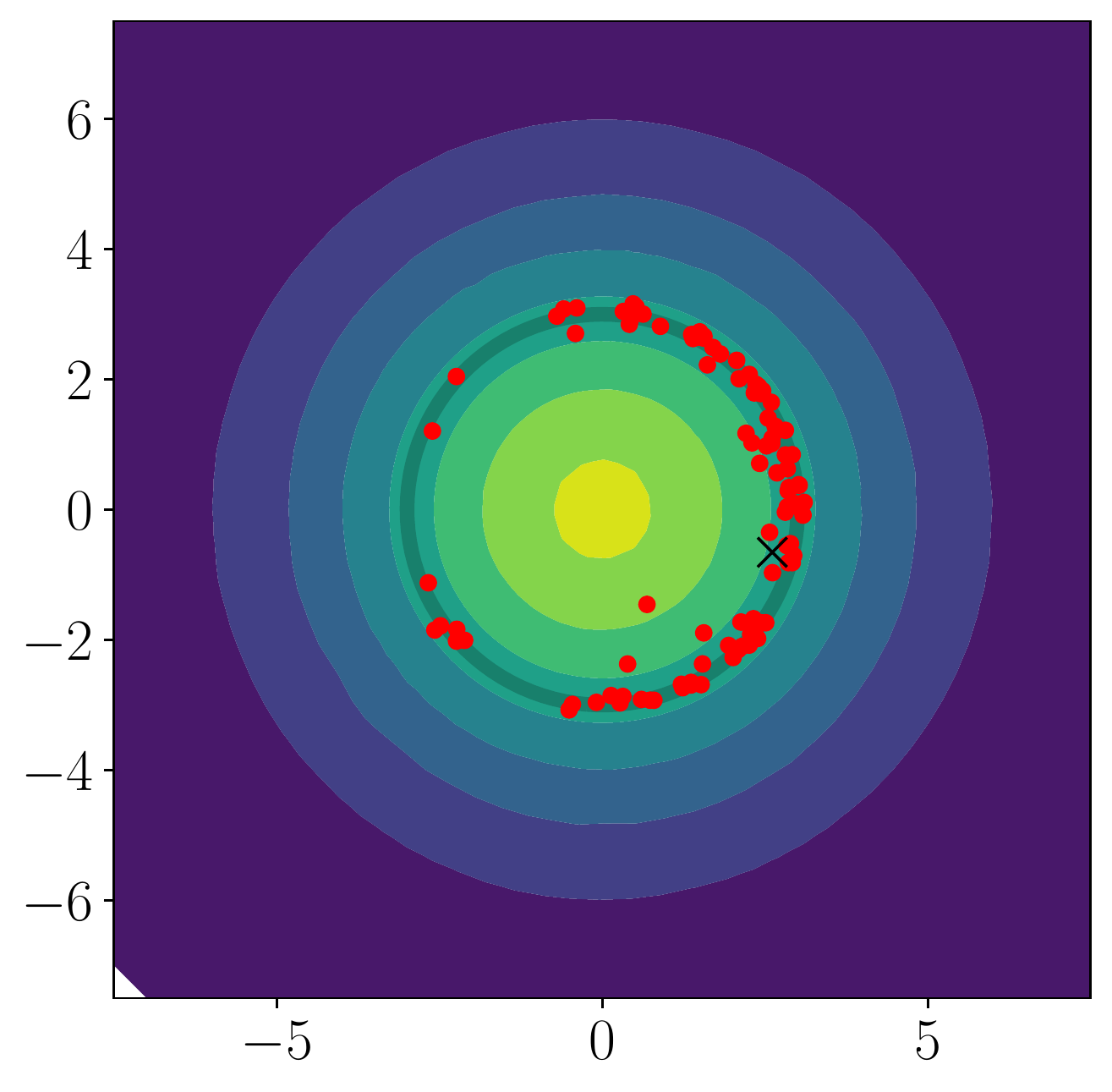}
    \caption{$J=50, M=100, T = 0, \alpha = 30, \nu = 0/1.$ Left: weighted averages with $\nu = 0$ of $100$ samples at the beginning of the simulation. Right: weighted averages with $\nu = 1$ of $100$ samples at the beginning of the simulation. The contour plots show the initial distribution of the mean-field solution. The influence of the constraint in the weight is clearly visible.}
    \label{fig:influenceWeightConstraint}
\end{figure}

\Cref{fig:pde-sde-comparison} (right) illustrates the weighted means after convergence of~$M=100$ independent simulations,
as well as the contour plot of the PDE solution after convergence to one triangle.
The weighted means lie at the minimizer of the shifted Ackley function and in particular in the support of the PDE solution.
This demonstrates that the SDE and PDE models are in good agreement,
even with as few as 50 particles.
Note that the left panel in \cref{fig:pde-sde-comparison} covers the domain $[-7.5,7.5]^2$,
while the right panel covers the zoomed-in region~$[-1.5,1.5]^2.$

Next, we study the influence of the penalization~\eqref{eq:weight_w_constraint}.
We consider again the Ackley function shifted by $(0.5,0)$,
with now a constraint given by the circle $B = \{ x \in \real^2 \colon x^2 = 9\},$
which is depicted in gray in~\cref{fig:influenceWeightConstraint}.
In this figure, the filled contours illustrate the initial density of the PDE solution,
while the red dots are the weighted averages corresponding to $M=100$ independent particle simulations at the initial time.
The black crosses indicate the positions of the weighted averages of the PDE solution at the initial time,
without~(left, calculated using~\eqref{eq:weightedmean}) or with (right, calculated using~\eqref{eq:weight_w_constraint}) penalization.
We observe that the weighted average with penalization of the initial PDE solution is already very close to the true global minimizer of the constrained problem.
This numerical experiment motivates the strategy of enforcing the constraints through the penalization on the objective functions.


\section{EKI with constraints}
\label{sec:eki_with_constraints}
The EKI was initially proposed as a method for solving inverse problems of the following form:
find an unknown parameter $x \in \real^d$ from data $y \in \real^K$ given that
\begin{align}
\label{eq:inverse_problem}
y = G(x) + \eta, \qquad \eta \sim \mathcal N(0, \Gamma),
\end{align}
where $G: \real^d \rightarrow \real^K$ is a nonlinear map known as the \emph{forward operator}.
A point estimator for the solution to~\eqref{eq:inverse_problem} is usually defined as the minimizer of a least-squares functional of the form
\begin{align}
\label{eq:regularized_least_squares}
\Phi_R(x; y) = \frac{1}{2}\eucnorm{y - G(x)}_{\Gamma}^2 + \frac{1}{2} \eucnorm{x - a}_{\Sigma}^2.
\end{align}
Here $\eucnorm{\cdot}_{A}$ is a short-hand notation for $\eucnorm{A^{-1/2}\cdot}$ and
the second term is a regularization parametrized by a vector $a \in \real^d$ and a matrix $\Sigma \in \real^{d \times d}$.
In the Bayesian approach to inverse problems,
these parameters can be viewed as the mean and variance of a Gaussian \emph{prior distribution} $\normal(a, \Sigma)$ on the unknown $x$.
In this case,
the minimizer of $\Phi_R(\,\cdot\,; y)$ is the pointwise maximizer of the Bayesian posterior distribution associated with the problem~\cite{S10}.
The EKI is an optimization scheme for the functional $\Phi_R(\,\cdot\,; y)$
that can be derived from the ensemble Kalman filter.
Its continuous-time version is based on the following interacting particle system:
\begin{equation}
    \label{eq:eki:evolution_of_thetas}%
    \begin{aligned}[b]
    \d x^{(j)}_t
    &= - \frac{1}{J} \sum_{k=1}^{J} \eip*{G(x^{(k)}_t) - \bar G_t}{G(x^{(j)}_t) - y}[\Gamma] \bigl(x^{(k)}_t - \bar x_t\bigr) \, \d t \\
    &\qquad - \mathcal C(\mu^{J}_t) \Sigma^{-1} (x^{(j)}_t - a) \, \d t, \qquad j  = 1, \dotsc, J,
    \end{aligned}
\end{equation}
where $\mu^J_t = \frac{1}{J} \sum_{j=1}^{J} \delta_{x^{(j)}_t}$ denotes as before the empirical measure associated with the ensemble,
and the quantities $\bar x_t$ and $\bar G_t$ are defined as
\begin{align*}
    \bar x_t &= \int_{\real^d} x \, \mu^J_t(\d x) = \frac{1}{J} \sum_{j =1}^{J} x^{(j)}_t, \\
    \bar G_t &= \int_{\real^d} G(x) \, \mu^J_t(\d x) = \frac{1}{J} \sum_{j=1}^{J} G(x^{(j)}_t).
\end{align*}
The matrix $\mathcal C(\mu^J_t)$ is the covariance under the empirical distribution:
\begin{align*}
    \mathcal C(\mu^J_t) = \int_{\real^d} (x \otimes x) \mu^J_t(\d x) = \frac{1}{J} \sum_{j=1}^{J} (x^{(j)}_t - \bar x_t) \otimes (x^{(j)}_t - \bar x_t).
\end{align*}
The first argument in the inner product on the right-hand side of \eqref{eq:eki:evolution_of_thetas} is related to consensus,
whereas the second argument measures the mismatch with the observed data.
In the limit~$t \to \infty$,
the particles are expected to concentrate at the global minimizer of $\Phi_R$ given in~\eqref{eq:regularized_least_squares}.

In the case of a linear forward map $G$,
the EKI algorithm takes the form of a gradient descent preconditioned by the ensemble covariance:
\begin{equation}
    \label{eq:gradient_eki}
    \d x^{(j)}_t = -\mathcal C(\mu^J_t) \nabla \Phi_R (x^{(j)}_t) \, \d t, \qquad j  = 1, \dotsc, J.
\end{equation}
In the case of a general forward map,
the EKI~\eqref{eq:eki:evolution_of_thetas} can be viewed as a derivative-free approximation of~\eqref{eq:gradient_eki}.
This viewpoint, which is adopted in~\cite{2019arXiv190308866G,garbuno2020affine},
is useful below:
we first show how to include constraints for~\eqref{eq:gradient_eki},
and then perform a derivative-free approximation of the resulting system.

In the Bayesian approach to inverse problems,
it would be natural to incorporate equality and inequality constraints in the prior distribution on the unknown parameter $x$.
The EKI method, however, applies only to the case of a Gaussian prior distribution,
i.e.~regularization by a quadratic function in~\eqref{eq:regularized_least_squares},
and so extensions of this algorithm are required in order to handle nonlinear constraints.

In this section,
we demonstrate how a simple penalty-based method similar to that in \cref{sec:cbo_for_constrained_global_optimization_problems} for CBO
can be employed for including constraints in EKI.
For the sake of simplicity, we assume there is only one equality constraint $\mathcal A(x) = 0$.
This constraint can be incorporated in the gradient method~\eqref{eq:gradient_eki} by
adding a penalization term to the regularized least-squares functional:
\begin{align*}
    \d x^{(j)}_t
    &= - \mathcal C(\mu^J_t) \nabla \left( \Phi_R  + \frac{1}{\nu} \abs{\mathcal A}^2 \right)(x^{(j)}_t) \, \d t, \\
    &= - \mathcal C(\mu^J_t) \nabla \Phi_R (x^{(j)}_t) \, \d t
    - \frac{2}{\nu} \mathcal C(\mu^J_t)  \mathcal A (x^{(j)}_t) \nabla\mathcal A(x^{(j)}_t) \, \d t, \,\,\, j  = 1, \dotsc, J.
\end{align*}
A derivative-free version of this methodology is obtained by approximating gradients in the same manner as in standard EKI:
\begin{align}
    \label{eq:eki:evolution_of_thetas_constraints}%
    \begin{split}
        \d x^{(j)}_t
        &= - \frac{1}{J} \sum_{k=1}^{J} \eip*{G(x^{(k)}_t) - \bar G_t}{G(x^{(j)}) - y}[\Gamma] (x^{(k)} - \bar x_t) \d t \\
        &\qquad - \mathcal C(\mu^J_t) \Sigma^{-1} (x^{(j)} - a) \, \d t
        - \frac{2}{\nu} \mathcal C(\mu^J_t) \mathcal A(x^{(j)}_t) \nabla \mathcal  A(x^{(j)}_t) \, \d t,
    \end{split}
\end{align}
for $j = 1, \dotsc, J$.
In contrast with the other methods proposed in the literature,
the particles forming the ensemble produced by this method are not confined to the feasible region.
Therefore, as already emphasized in the introduction,
a prerequisite of the method is that the forward map may be evaluated for any $x \in \real^d$,
and not only in the feasible region.


The formulation~\eqref{eq:eki:evolution_of_thetas_constraints} also assumes that
the gradient of the function~$\mathcal A$ can be calculated analytically.
In the general case, one may employ a gradient-free approximation of $\nabla \lvert \mathcal A \rvert^2$,
\begin{equation}
    \label{eq:eki:evolution_of_thetas_constraints_derivative_free}%
    \begin{aligned}
        \d x^{(j)}_t
        =\,& - \frac{1}{J} \sum_{k=1}^{J} \eip*{G(x^{(k)}_t) - \bar G_t}{G(x^{(j)}_t) - y}[\Gamma] x^{(k)}_t \, \d t
        - \mathcal C(\mu^J_t) \Sigma^{-1} (x^{(j)}_t - a) \, \d t, \\
        &- \frac{2}{\nu J} \sum_{k=1}^{J} \eip*{\mathcal A(x^{(k)}_t) - \bar {\mathcal A}_t}{\mathcal  A(x^{(j)}_t)} (x^{(k)}_t - \bar x_t) \, \d t,
        \qquad j = 1, \dotsc, J,
    \end{aligned}
\end{equation}
where $\bar {\mathcal A}_t = \frac{1}{J} \sum_{j=1}^{J} \mathcal A(x^{(j)}_t)$.
The gradient-based and gradient-free constraint terms coincide in the case where $\mathcal A$ is linear.

\begin{remark}
The scheme \eqref{eq:eki:evolution_of_thetas_constraints_derivative_free} can be rewritten as
\begin{align}
    \label{eq:unified_eki}
    \begin{aligned}[b]
    \d x^{(j)}_t
    &= - \frac{1}{J} \sum_{k=1}^{J} \eip*{\mathcal G(x^{(k)}_t) - \bar {\mathcal G}_t}{\mathcal G(x^{(j)}_t) - \widetilde y}[\widetilde \Gamma] (x^{(k)}_t - \bar x_t) \, \d t, \qquad j  = 1, \dotsc, J,
    \end{aligned}
\end{align}
where $\mathcal G(x) = \bigl(G(x), x, \mathcal A(x)\bigr)^\t$,
$\widetilde y = \bigl(y, a, 0_d\bigr)^\t$ and
\[
    \widetilde \Gamma =
    \begin{pmatrix}
        \Gamma & 0 & 0 \\
        0 & \Sigma & 0 \\
        0 & 0 & \nu
    \end{pmatrix}.
\]
Numerically, implementing~\eqref{eq:eki:evolution_of_thetas_constraints_derivative_free} presents the same level of difficulty as implementing unconstrained EKI.
The constraint is interpreted here as an additional observation,
with very small associated noise.
\end{remark}

\section{Numerical results for CBO}
\label{sec:numerical_results}
In this section,
we discuss numerical results for CBO with constraints.
We first present the numerical scheme employed throughout this section in \cref{sub:numerical_schemes},
and then we investigate the qualitative behavior of the method.

In \cref{sub:cbo_equality} and \cref{sub:cbo_inequality},
we present numerical results for the case of equality and inequality constraints, respectively.
In \cref{sub:cbo_many_particle}, we investigate the convergence to the mean-field equation~\eqref{eq:mean_field_cbo_with_constraints} in the many-particle limit.
In \cref{sub:cbo_comparison}, we compare the proposed method with another approach in the literature~\cite{fornasier2020consensus,FornasierEtAl2}.
Finally, in \cref{sub:numerical_experiments_for_the_particle_systems},
we present a brief parameter study,
with the aim of highlighting the role of the parameters entering in the method.

In all the numerical experiments presented in this section,
the objective function is defined from the Ackley function in two dimension:
\begin{equation*}
    \label{eq:ackley}
    f_A(x) = -20 \exp \left( - \frac{1}{5} \sqrt{\frac{1}{2} \sum_{i=1}^{2} |x_i|^2} \right)
    -\exp \left( \frac{1}{d} \sum_{i=1}^{2} \cos\bigl(2 \pi x_i\bigr) \right) + \e \,+\, 20.
\end{equation*}

\subsection{Numerical schemes}%
\label{sub:numerical_schemes}

As mentioned, the introduction of a relaxation drift  towards the feasible manifold with $\varepsilon \ll 1$ leads to stiffness of the CBO particle system~\eqref{eq:cbo_with_constraints}.
This leads to a stringent limitation on the time step when using an explicit method such as Euler--Maruyama~\cite{Higham}.
In the case of a single quadratic equality or inequality constraint,
this issue can be mitigated by using
a semi-implicit scheme.
Let us consider, for example, an equality constraint of the form
\begin{equation}
    \label{eq:constraint_numerics}
    \mathcal A(x) = x^\t A x - c = 0
\end{equation}
for a positive definite matrix $A \in \real^{d \times d}$ and a positive real number $c$.
This general form of the constraint encompasses elliptic and, in particular, the circular equality constraints considered throughout this section. Moreover, it is simple to extend the method to similar inequality constraints.
In order to use CBO with the constraint~\eqref{eq:constraint_numerics},
we propose the discretization
\begin{subequations}
    \label{eq:semi-implicit}
\begin{align}
    & y^{(j)}_n = x^{(j)}_n - \bigl(x^{(j)}_n - m_g(\mu^J_n)\bigr) \Delta t + \sqrt{2 \Delta t} \, \sigma \, \bigl( x^{(j)}_n - m_g(\mu^J_n) \bigr) \circ \xi^{(j)}_n, \\
    &x^{(j)}_{n+1} =  \left(I_d + \frac{4 \Delta t}{\varepsilon} \mathcal A(x^{(j)}_n) A\right)^{-1} y^{(j)}_n,
\end{align}
\end{subequations}
where $\xi^{(j)}_n \sim \mathcal N(0,1)$,
for $n=1, 2, \dotsc$ and $j = 1, \dotsc, J$, are independent.
This discretization arises when treating the term $\grad \mathcal A$ in~\eqref{eq:cbo_with_constraints} implicitly
and all the other terms on the right-hand side explicitly.
The PDE associated with CBO is solved with the help of the Python free software library \texttt{Fenics}~\cite{Fenics}
on a mesh refined locally around the feasible manifold using the open source software \texttt{gmsh}~\cite{gmsh}.
We use the discontinuous Galerkin method with elements of degree zero.
The drift term is computed with upwind flux and for the diffusion we use the discretization proposed in \cite{Kulkarni}.

In order to make the SDE and PDE results comparable,
we initialize the particle SDE~\eqref{eq:cbo_with_constraints}-\eqref{eq:weight_w_constraint} with $x_0^{(j)} \sim \mathcal N(0,3 I_2)$,
and the corresponding PDE \eqref{eq:mean_field_cbo_with_constraints} is initialized with initial condition $\rho_0(x) = \frac{1}{\sqrt{2\pi}}e^{-\frac{1}{2}\left( \frac{x}{3} \right)^2}$,
i.e. $\rho_0$ is the probability density function of $\mathcal N(0, 3I_2)$.
Unless otherwise stated,
we always assume these initial conditions at the particle and PDE levels.
The following parameters are fixed for all simulations:
$\sigma = 0.7, \alpha = 30, \Delta t_\text{PDE}= \num{5e-3}, \Delta t_\text{\rm SDE}= \num{5e-4}.$
The number of particles $J$,
and the relaxation parameter $\varepsilon$,
are specified on a case-by-case basis.

\subsection{CBO with equality constraint}
\label{sub:cbo_equality}
We now present the results of particle simulations of the CBO method with constraints~\eqref{eq:cbo_with_constraints},
which incorporates both the relaxation drift towards the constraint and the penalization~\eqref{eq:weight_w_constraint},
and of PDE simulations of the associated mean field equation~\eqref{eq:mean_field_cbo_with_constraints}.
We begin by considering the following optimization problem with circular equality constraint,
\[
    \argmin_{x \in B} f_{A}(x - x_{*}),
    \qquad x_{*} = \begin{pmatrix} 3 \\ 0 \end{pmatrix},
    \qquad B = \{ x \in \real^2 \colon x^2 = 9\}.
\]
Here $f_A$ is the Ackley function, which is minimized at $x = 0$,
so that the solution of this problem is given by $x_* = (3, 0)$.
(Note that this corresponds to the minimizer of the unconstrained global optimization problem.)

\begin{figure}[ht!]
    \centering
    \includegraphics[width=.35\linewidth]{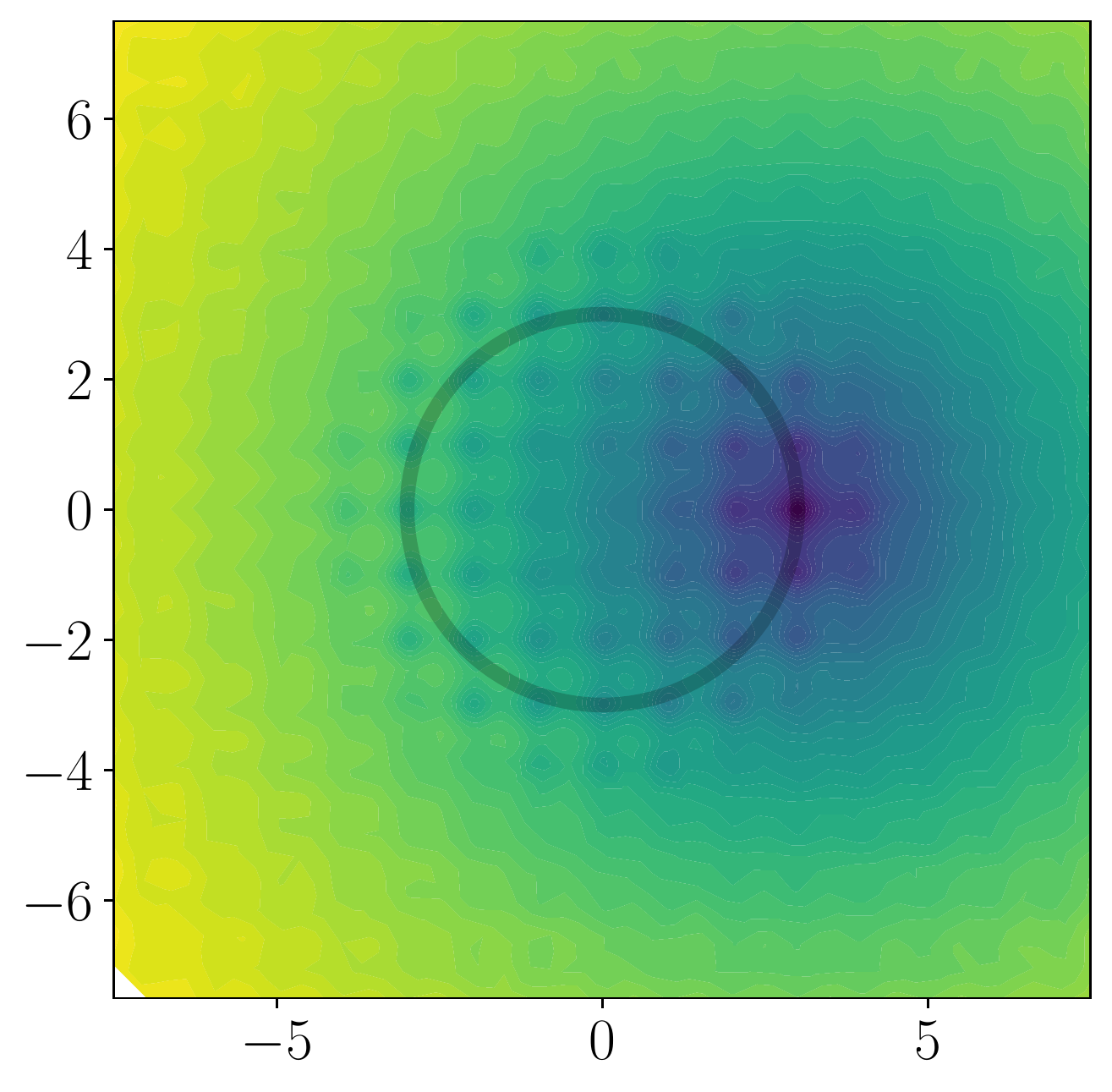}
    \includegraphics[width=.35\linewidth]{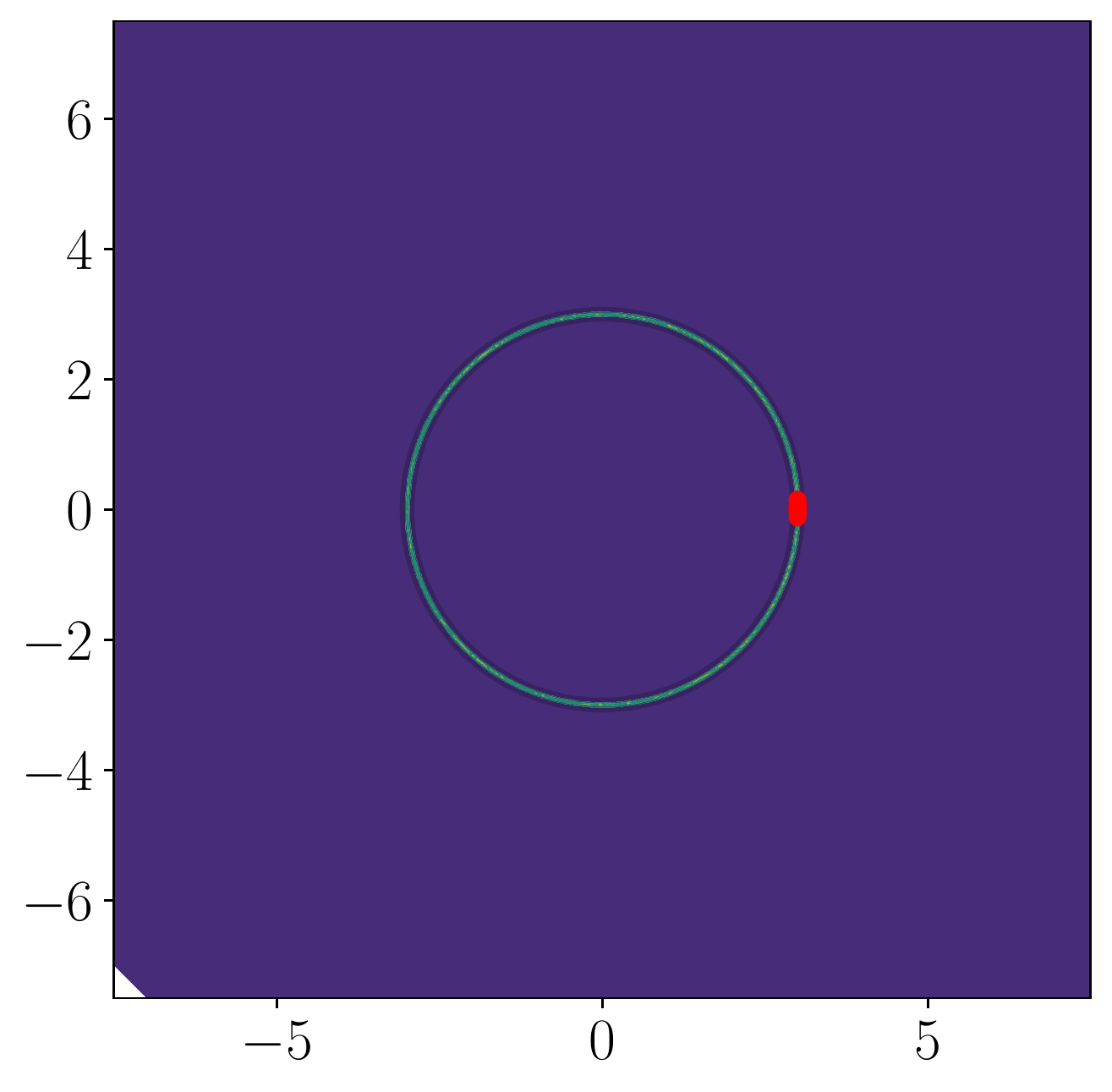}
    \includegraphics[width=.35\linewidth]{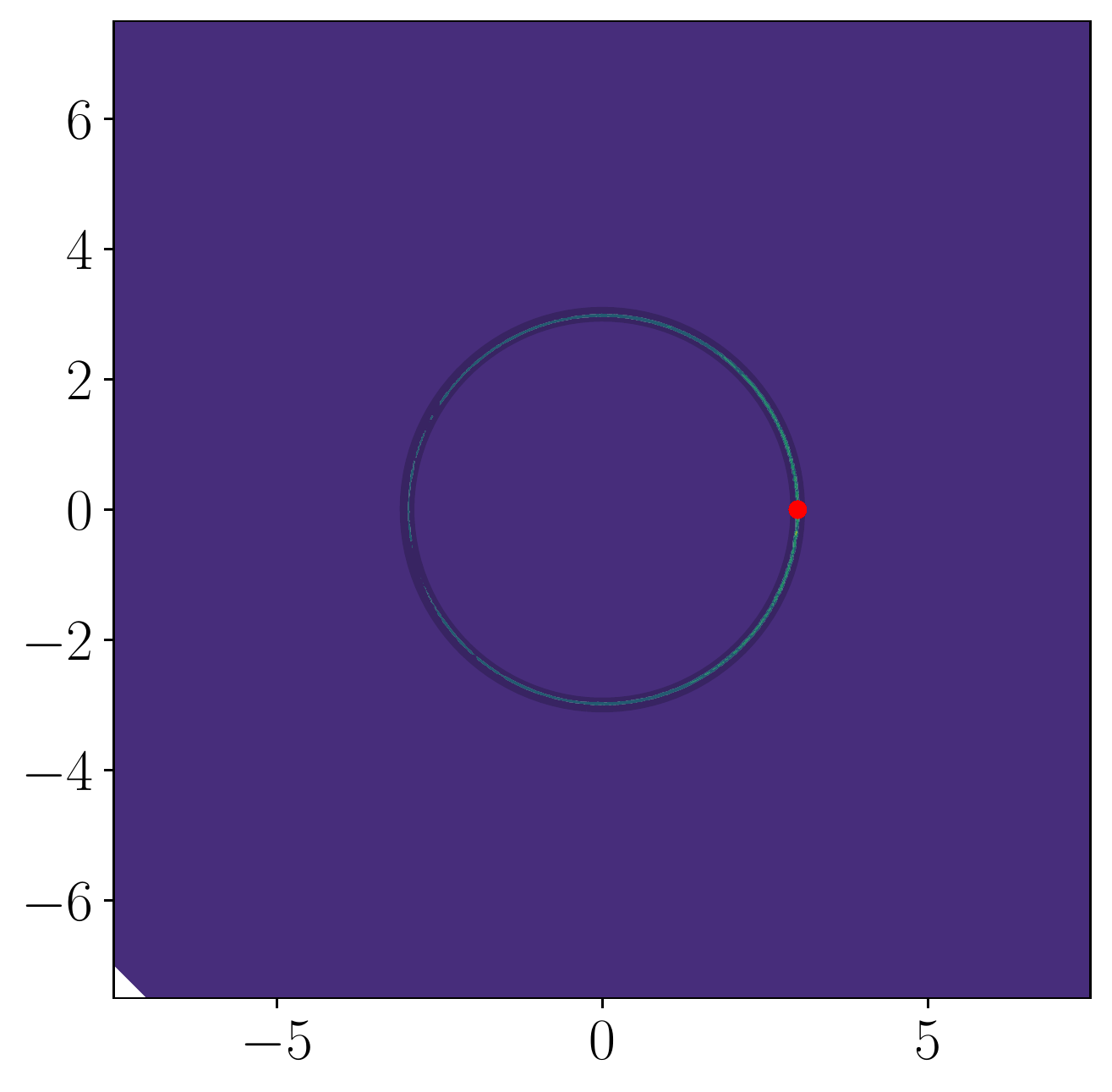}
    \includegraphics[width=.35\linewidth]{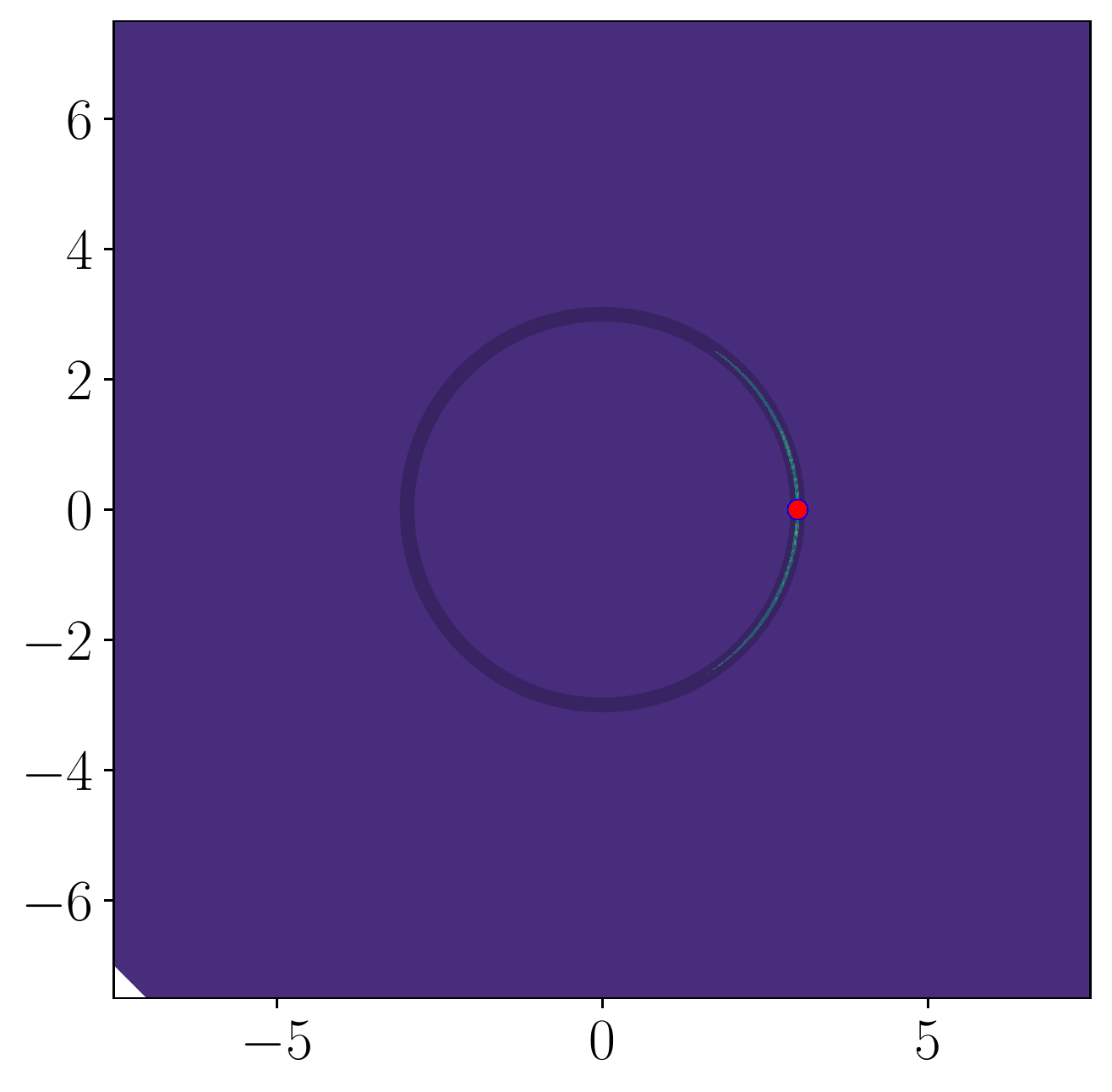}
    \caption{%
        \textbf{Top left panel}: Ackley function shifted by $(3,0),$ and feasible set
        \(
            B = \{ x \in \real^2 \colon x^2= 9\}.
        \)
        \textbf{Other panels}: PDE solution at $t=0.01$, $t=0.15$ and $t=0.5$ (filled contour),
        and weighted averages of $M=100$ independent simulations with $J=50$ particles each (red dots).
        We observe that the weighted averages lie in the support of the mean-field solution and provide a good approximation of the global minimizer.
        The parameters employed for the simulation are the following:
        $J=50, M=100, \Delta t_{\rm SDE}=0.0005, \Delta t_{\rm PDE}=0.005, h_{\min}=0.125, h_{\max}=0.5, \alpha = 30, \nu = 1, \sigma = 0.7, \varepsilon=0.1$.
    }
    \label{fig:pdeSdeConstraint}
\end{figure}

The top left panel in~\cref{fig:pdeSdeConstraint} illustrates the shifted Ackley function (filled contour)
and the constraint set $B$ (gray line).
The top right, bottom left and bottom right panels depict the PDE solution (filled contour),
as well as the weighted averages corresponding to 100 independent simulations of the particle system (one red dot per independent simulation),
at times $0.01$, $0.15$ and $0.5$.
We observe that, already for $t = 0.01$,
the weighted averages provide a good approximation of the global minimizer.

\subsection{CBO with inequality constraint}
\label{sub:cbo_inequality}
\begin{figure}[ht]
    \centering
    \includegraphics[width=.35\linewidth]{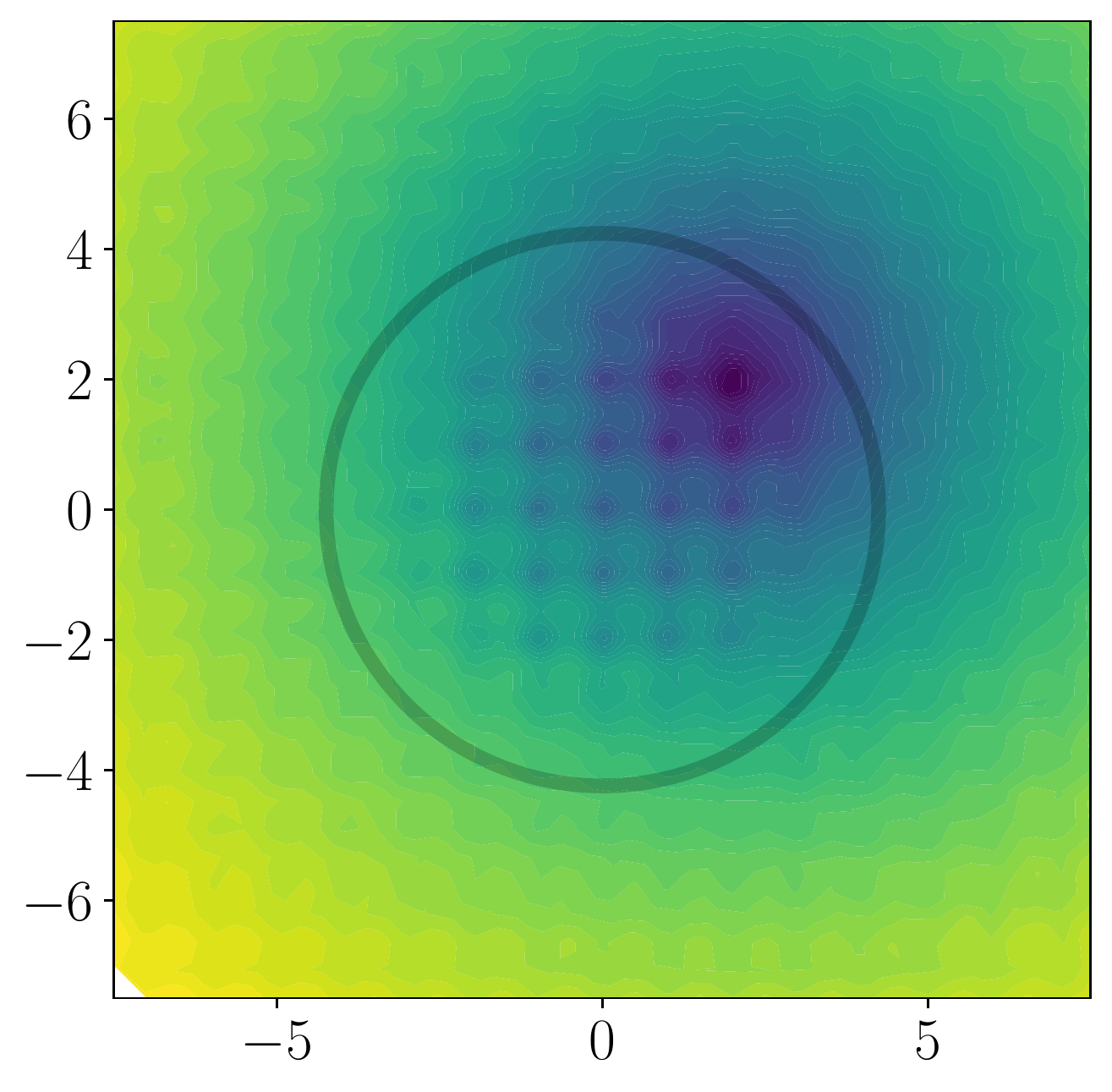}
    \includegraphics[width=.35\linewidth]{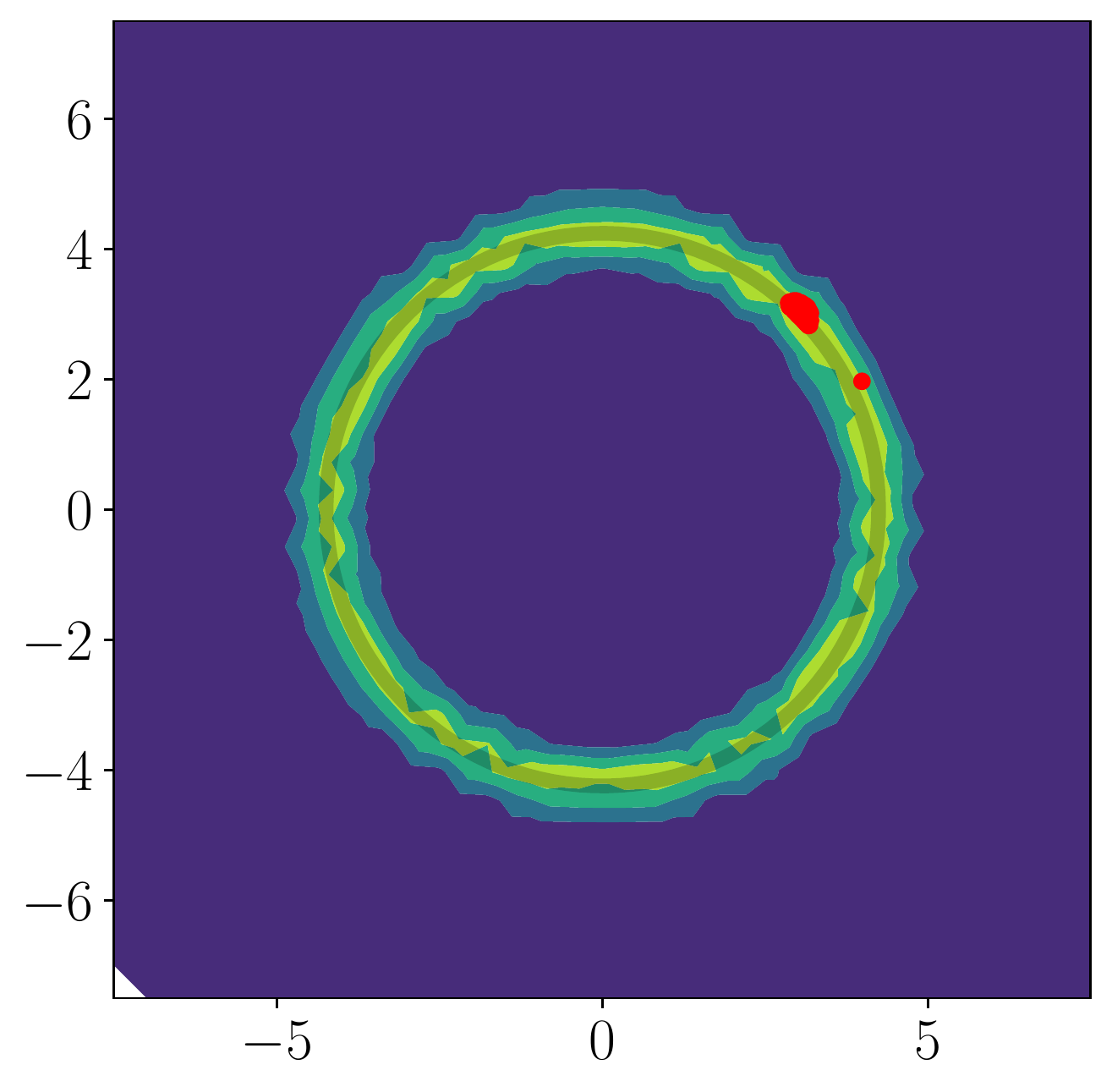}
    \includegraphics[width=.35\linewidth]{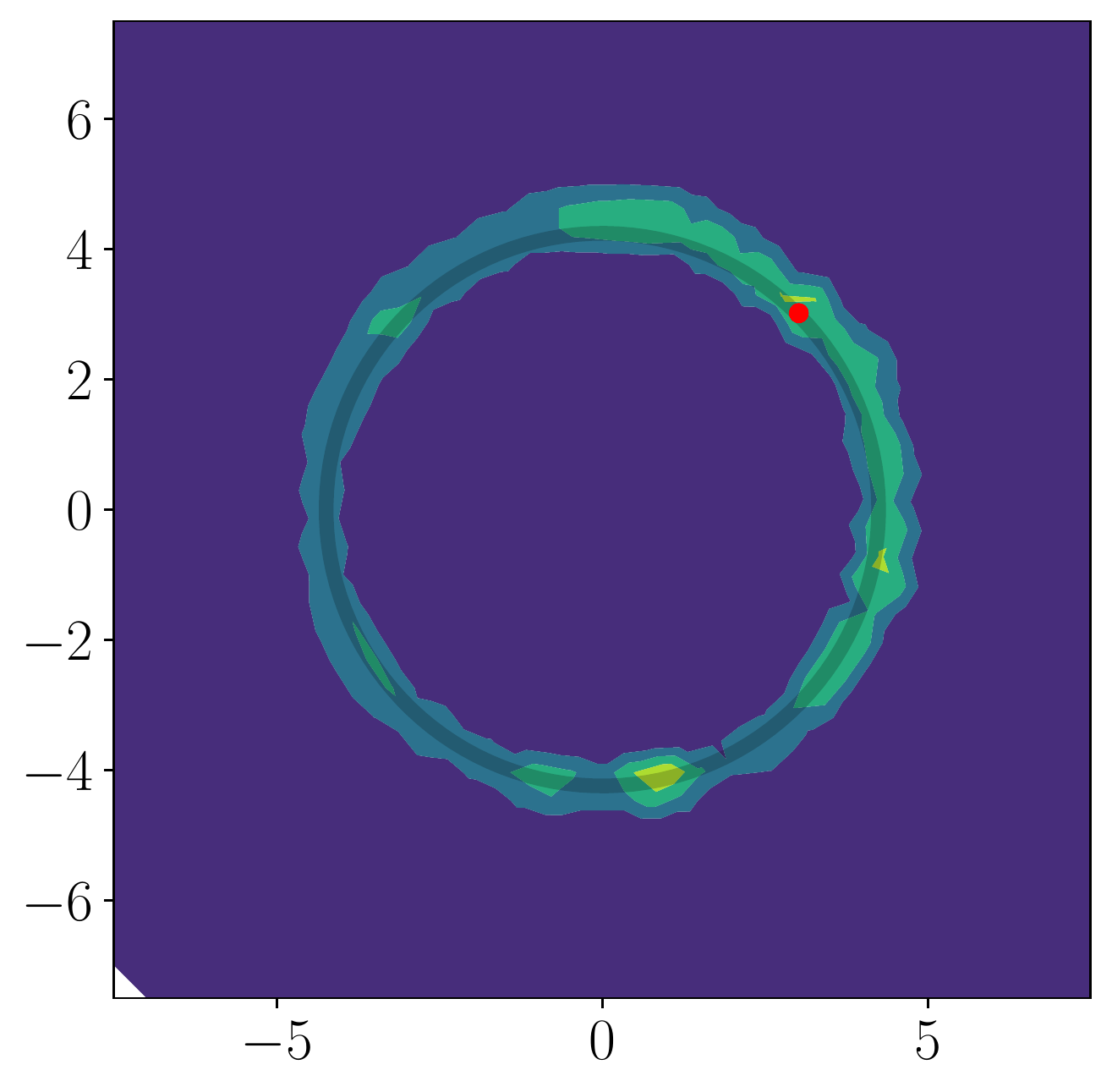}
    \includegraphics[width=.35\linewidth]{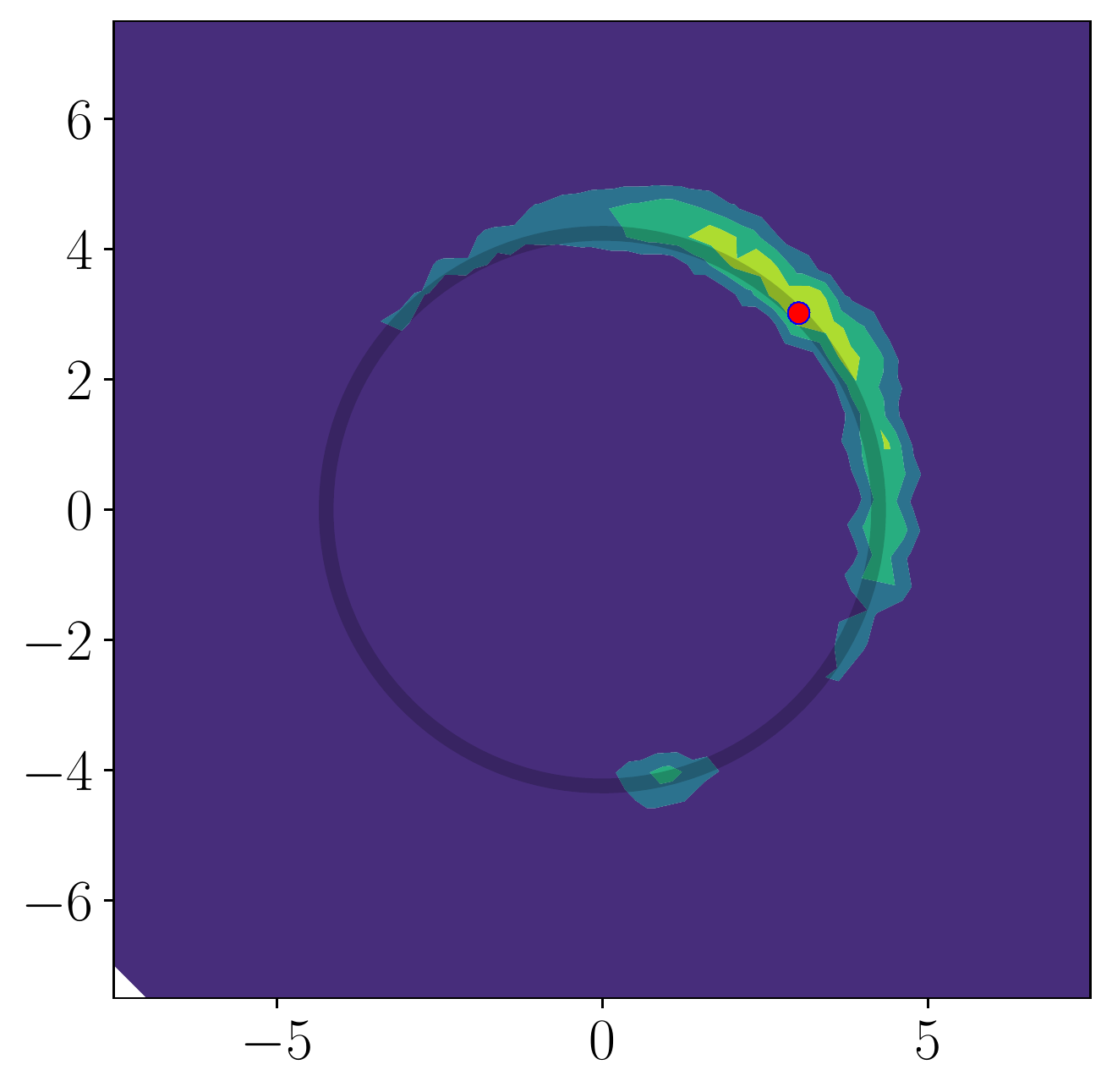}
    \caption{Parameters: $J=50, M=100, \Delta t_{\rm SDE}=0.0005, \Delta t_{\rm PDE}=0.005, h_{\min}=0.125, h_{\max}=0.5, \alpha = 30, \nu = 1, \sigma = 0.7, \varepsilon=0.1$ Top-left: Ackley function shifted by $(2,2)$ and the constraint $x^2 \ge 18$ (gray). Top-right to bottom-right: mean-field density and weighted means of the particle solution at $t=0.01, 0.15$ and $0.5.$ The weighted averages lie in the support of the mean-field solution and very close to the global minimum of the function. }
    \label{fig:pdeSdeInequalityOutside}
\end{figure}

Now, we investigate the following example with inequality constraint:
\[
    \argmin_{x \in B} f_{A}(x - x_{*}),
    \qquad x_{*} = \begin{pmatrix} 2 \\ 2 \end{pmatrix},
    \qquad B = \{ x \in \real^2 \colon x^2 \ge 18\},
\]
where $f_A$ is the standard Ackley function in two dimensions.
The associated numerical results are illustrated in~\cref{fig:pdeSdeInequalityOutside}.
Again, the red particles indicate the positions of the weighted means of $M=100$ independent simulations with $J=50$ particles each.
The top left panel of \cref{fig:pdeSdeInequalityOutside} shows the shifted Ackley function (filled contour),
and the boundary of the feasible set (gray line).
The top right, bottom left and bottom right panels depict the evolution of the PDE solution as well as the weighted means of the SDE simulations.
Compared to the previous example, the solution to the mean-field PDE is less concentrated,
which is expected since the constraint is of inequality type.
At $t=0.01$ (top right plot),
not all weighted averages are close to the minimizer.
At $t = 0.15$ and $t = 0.5$ (bottom plots),
the weighted averages are very close to the minimizer,
and the density concentrates around it.

We emphasize that the accuracy of the PDE results is strongly restricted by the mesh size;
since this is fixed during a simulation,
the PDE solver is unable to capture the evolution of the ensemble after this has concentrated to a very small area.
Comparisons between the SDE and PDE simulations are therefore meaningful only in the beginning of the simulation.
For this reason we restrict the following numerical studies to particle simulations.

\subsection{\texorpdfstring
    {Many particle limit $J \rightarrow \infty$ for CBO}
    {Many particle limit for CBO}}
\label{sub:cbo_many_particle}
\begin{figure}[ht]
    \centering
    \includegraphics[scale=0.5]{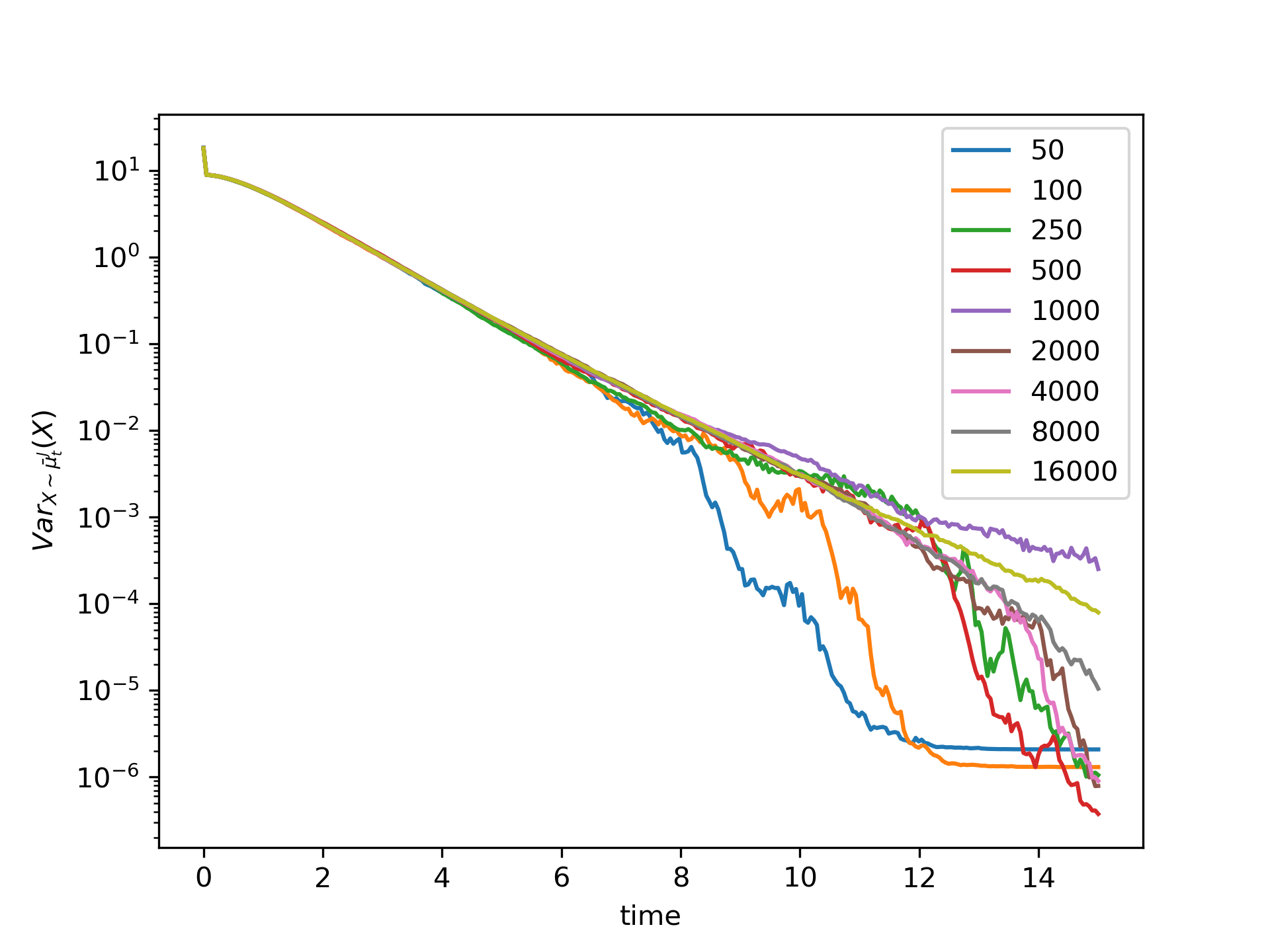}
    \caption{%
        Ensemble variance for different ensemble sizes. For the simulation we used the following parameters: $M=100, \Delta t_{SDE} = 0.0005, \alpha= 30, \nu=1, \sigma=0.7, \epsilon=0.1. $
    }
    \label{fig:variance}
\end{figure}

In the following we investigate the behavior of the proposed CBO scheme as $J$ tends to infinity.
For simplicity, we consider toy examples where the true minimizer~$x_*$ is known.
For the particle scheme we perform $M=100$ independent simulations for every setting considered,
each with a number~$J$ of particles.
For $j=1,\dots,J$ and $m=1,\dots,M,$ the vector $x^{(m,j)}_t \subset \real^d$ denotes the position of the $j$-th particle in the $m$-th independent simulation at time $t.$
First, we study in \cref{fig:variance} the evolution of the variance (more precisely, the trace of the covariance) of particle positions under an averaged empirical measure,
\[
    \var_{X \sim \bar\mu^J_t}(X)
    = \int \left\lvert x - \int y \,\bar \mu^{J}_t(\d y) \right\rvert^2 \, \bar\mu^{J}_t(\d x) ,
    \qquad \bar\mu^{J}_t = \frac{1}{M} \sum_{m=1}^M \mu^{J,m}_t ,
\]
where $\mu^{J,m}_t$ denotes the empirical measure for the $m$-th independent simulation at time $t$, for different ensemble sizes.
During the first few time steps, we observe a fast decrease of all variances,
as all the particles quickly move toward the feasible manifold.
Thereafter, the ensembles continue concentrating but more slowly.
As the figure shows,
starting from $t \approx 8$ ensembles with a smaller number of particles begin to collapse faster than ensembles with more particles.
This plot does not reveal any information on the proximity of the convergence point to the global minimizer, however.

\begin{figure}[ht]
    \centering
    \includegraphics[scale=0.32]{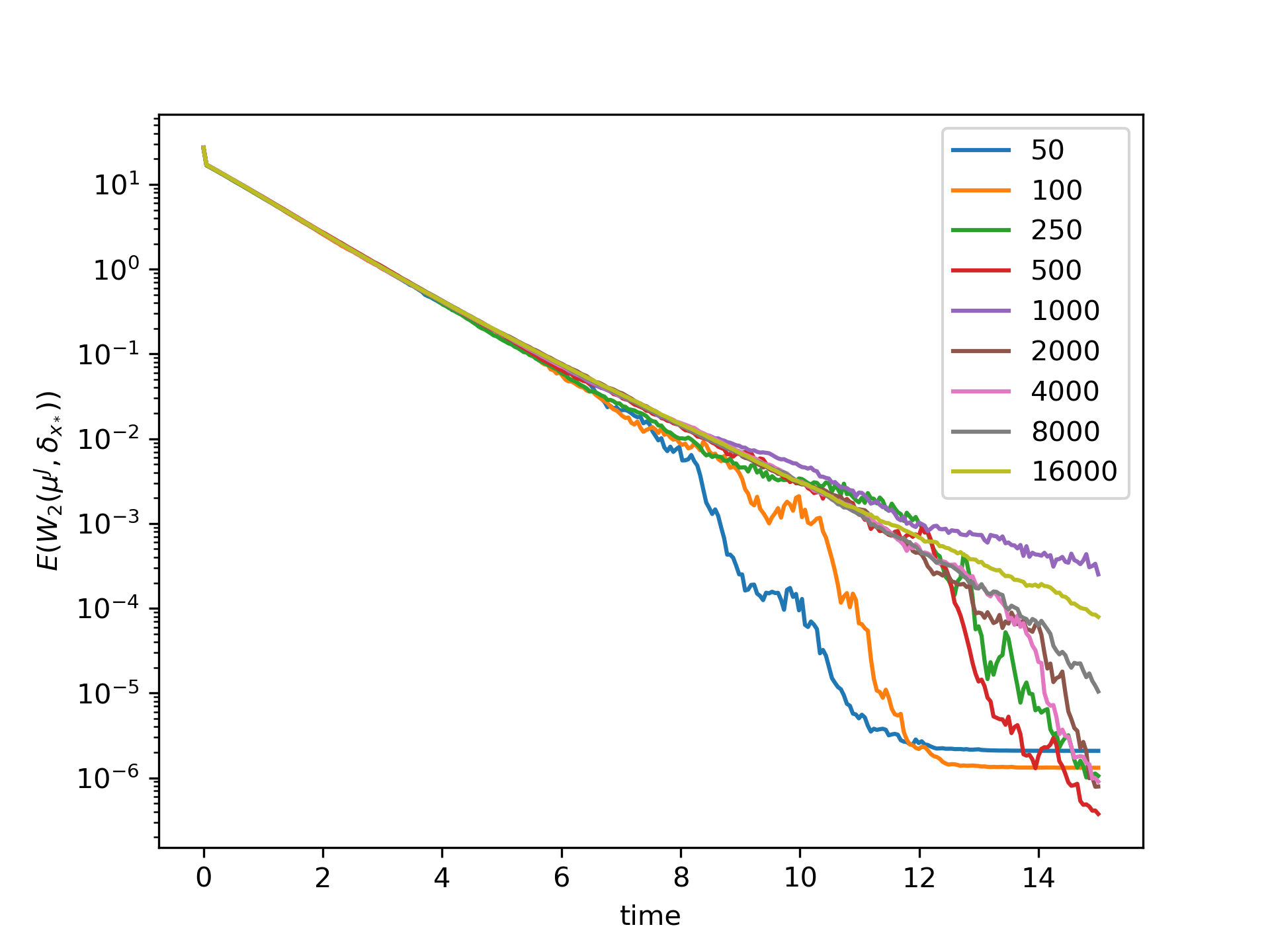}
    \includegraphics[scale=0.32]{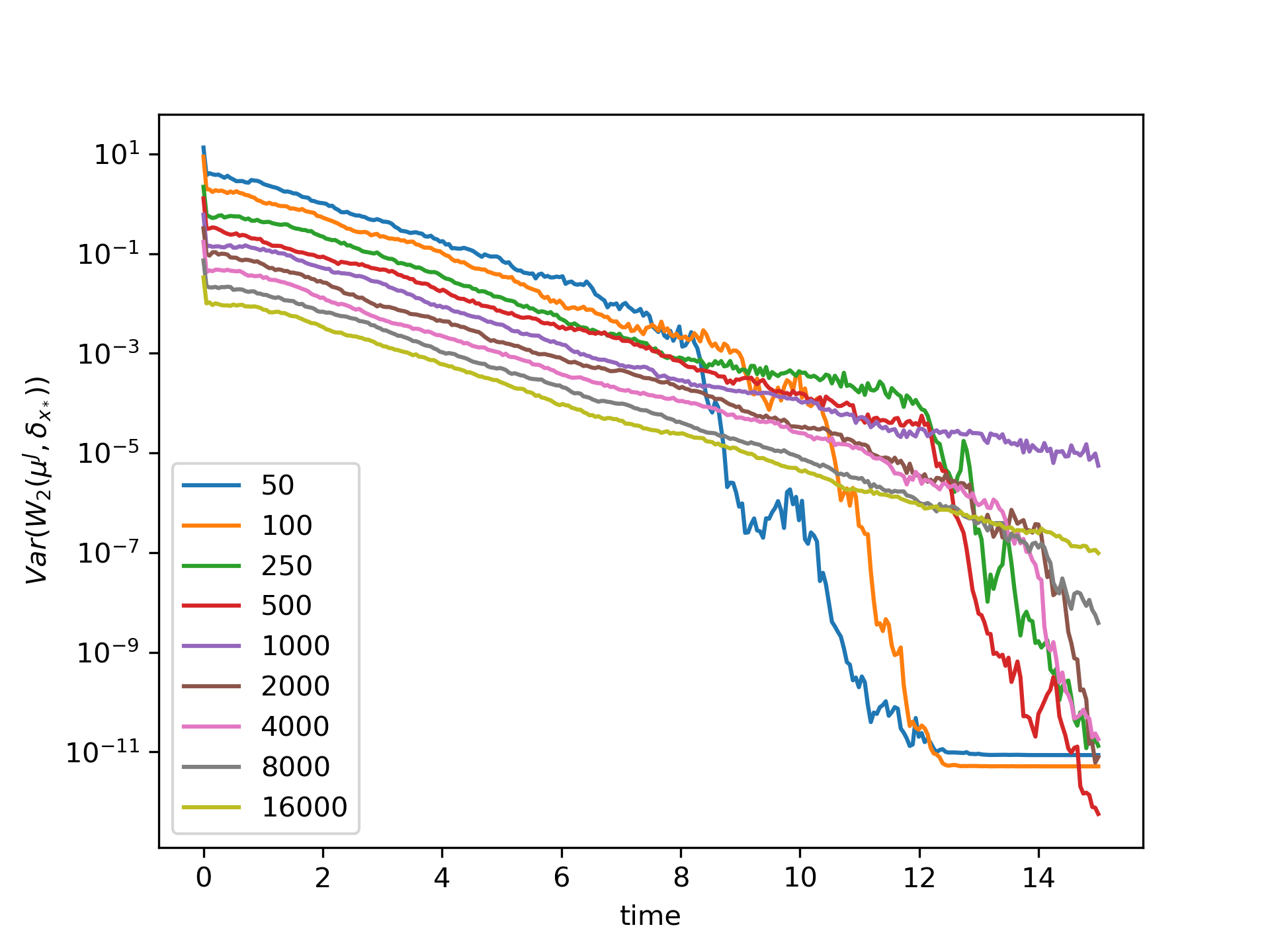}
    \caption{%
        Evolution in time of the expectation (\textbf{left}) and variance (\textbf{right}) of the Wasserstein distances.
        The parameters employed here are the same as in~\cref{fig:variance}.
    }
    \label{fig:WassVariance}
\end{figure}

The accuracy of the convergence point is discussed in~\cref{fig:WassVariance},
which depicts the time evolution of the expected value of the Wasserstein distance~$W_2(\mu_t^J,\delta_{x_*})$,
approximated in practice based on 100 independent simulations.
Again, we see a strong decrease at the beginning of the simulation as the particles are driven towards the manifold.
Then the expected values of the Wasserstein distances decrease exponentially until the ensembles collapse.
As we already observed in \cref{fig:variance},
smaller ensembles collapse earlier than larger ones.
The approximation of the true minimizer is very good already for small ensemble sizes.

\begin{figure}
    \centering
    \includegraphics[scale=0.5]{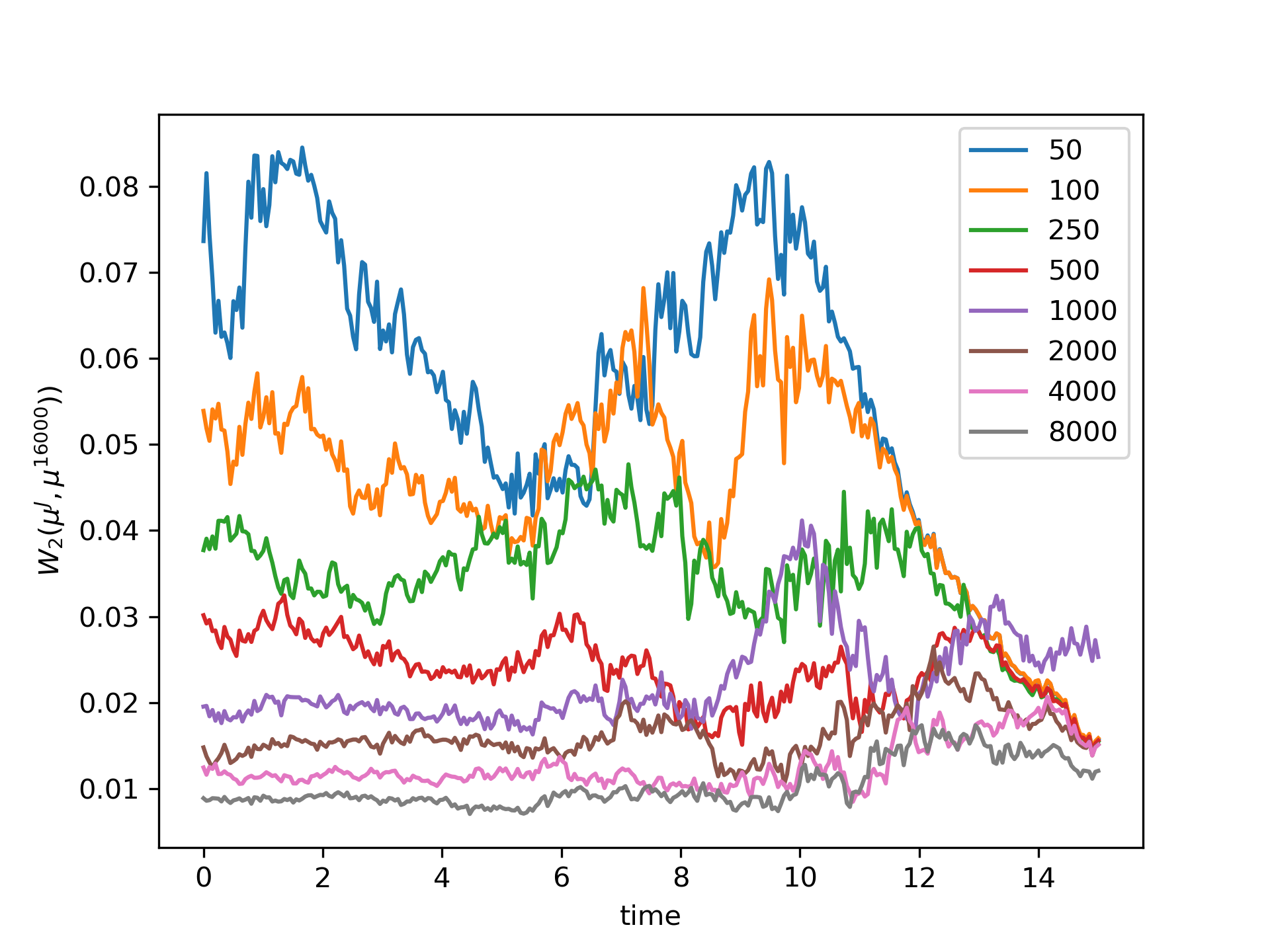}
    \caption{Wasserstein distance to reference solution with $J=16000.$ The parameters employed here are the same as in~\cref{fig:variance}.}
    \label{fig:Wass16}
\end{figure}

To illustrate the convergence to the deterministic mean-field dynamics in the limit as~$J\to\infty$,
we plot the Wasserstein distances of $W_2(\mu^J, \mu^{16000})$ for different values of $J$ in~\cref{fig:Wass16},
with this time only one simulation per value of $J$.
At $t=0$ we see the expected decrease of the Wasserstein distance along the $y$-axis. The distances vary over time but decrease at the end of the simulation as all ensembles converge towards the global minimizer.

\begin{figure}
    \centering
    \includegraphics[scale=0.5]{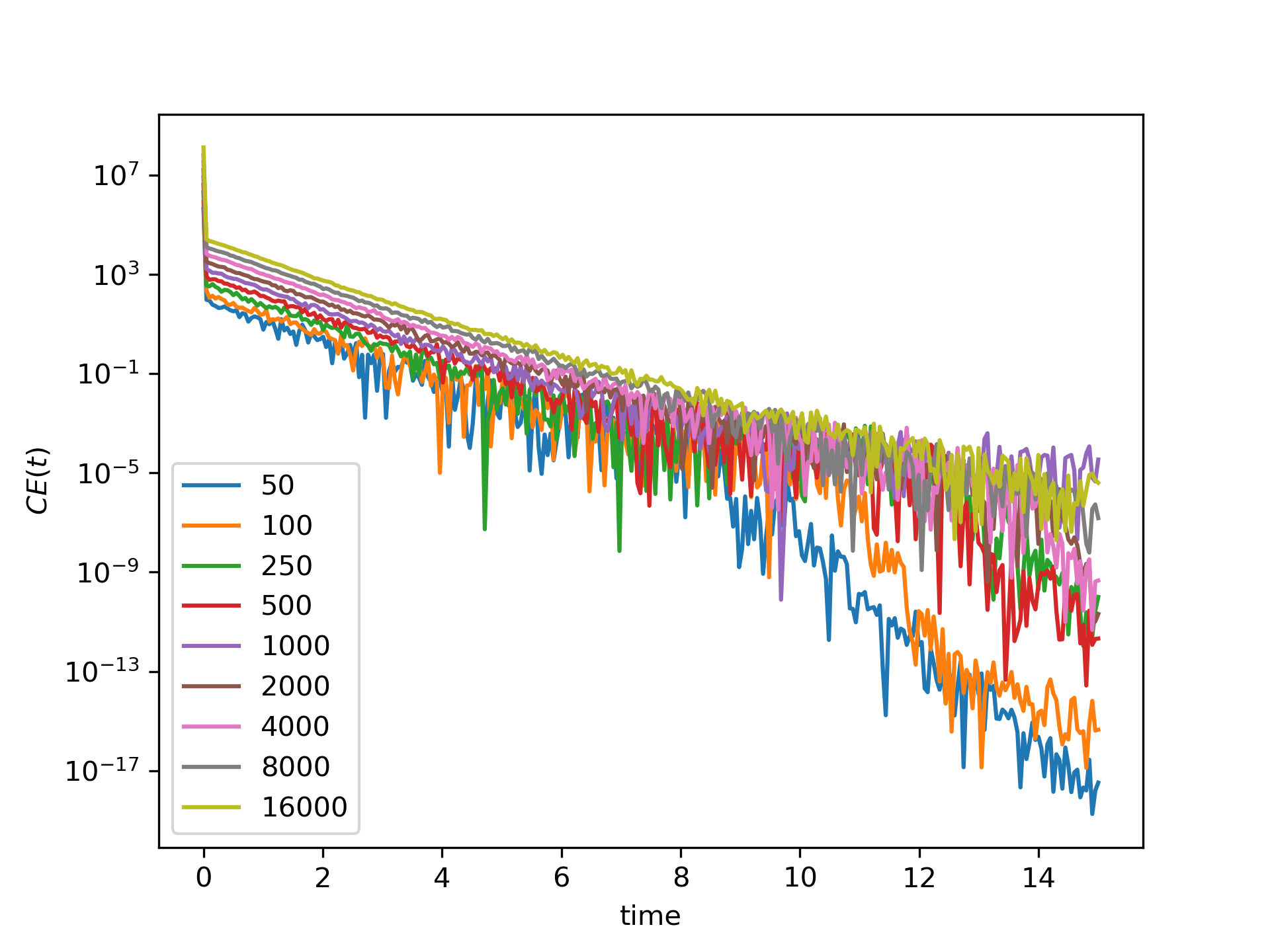}
    \caption{Distance to constraint measured with the quantity defined in \eqref{eq:CE}. The parameters employed here are the same as in~\cref{fig:variance}.}
    \label{fig:dist2con}
\end{figure}

In~\cref{fig:dist2con} the average distance of the ensembles to the constraint is illustrated with the help of the quantity
\begin{equation}\label{eq:CE}
    {\rm CE}\bigl(t)= \frac{1}{M}\sum_{m=1}^M \left(\frac1J \sum_{i=1}^J\A\bigl(x^{(m,j)}_t\bigr)^2\right)
        \approx \expect \left( \int_{\real^d} \lvert \mathcal A(x) \rvert^2 \, \mu^J_t(\d x) \right).
\end{equation}
This may be viewed as a ``constraint energy''; the smaller ${\rm CE}(t)$ is,
the closer the ensemble is to the constraint.
At the beginning of the simulation,
the quantity ${\rm CE}(t)$ decreases very fast thanks to the relaxation drift towards the manifold.
As time increases, we observe small fluctuations of the constraint energy,
but it is likely that these oscillations would disappear in the limit~$M \to \infty$ of many independent simulations.
The oscillations are due to the diffusion term acting on the particles,
which is stronger for particles far away for the weighted mean.
At $t=15$, all ensembles are close to the constraint.

\subsection{Comparison with a projection-based method}
\label{sub:cbo_comparison}
In the introduction,
we mentioned that, in the case of CBO,
the approach discussed here incorporates the constraint via two distinct mechanisms:
a penalization added to the objective function, and a relaxation drift towards the manifold.
This is in contrast with other methods that ensure that the constraint is satisfied at all times.
We therefore make a comparison of the method proposed here and a method using projection to the constraint in every time step, as proposed in~\cite{fornasier2020consensus}.
\begin{figure}
    \includegraphics[width=0.32\linewidth]{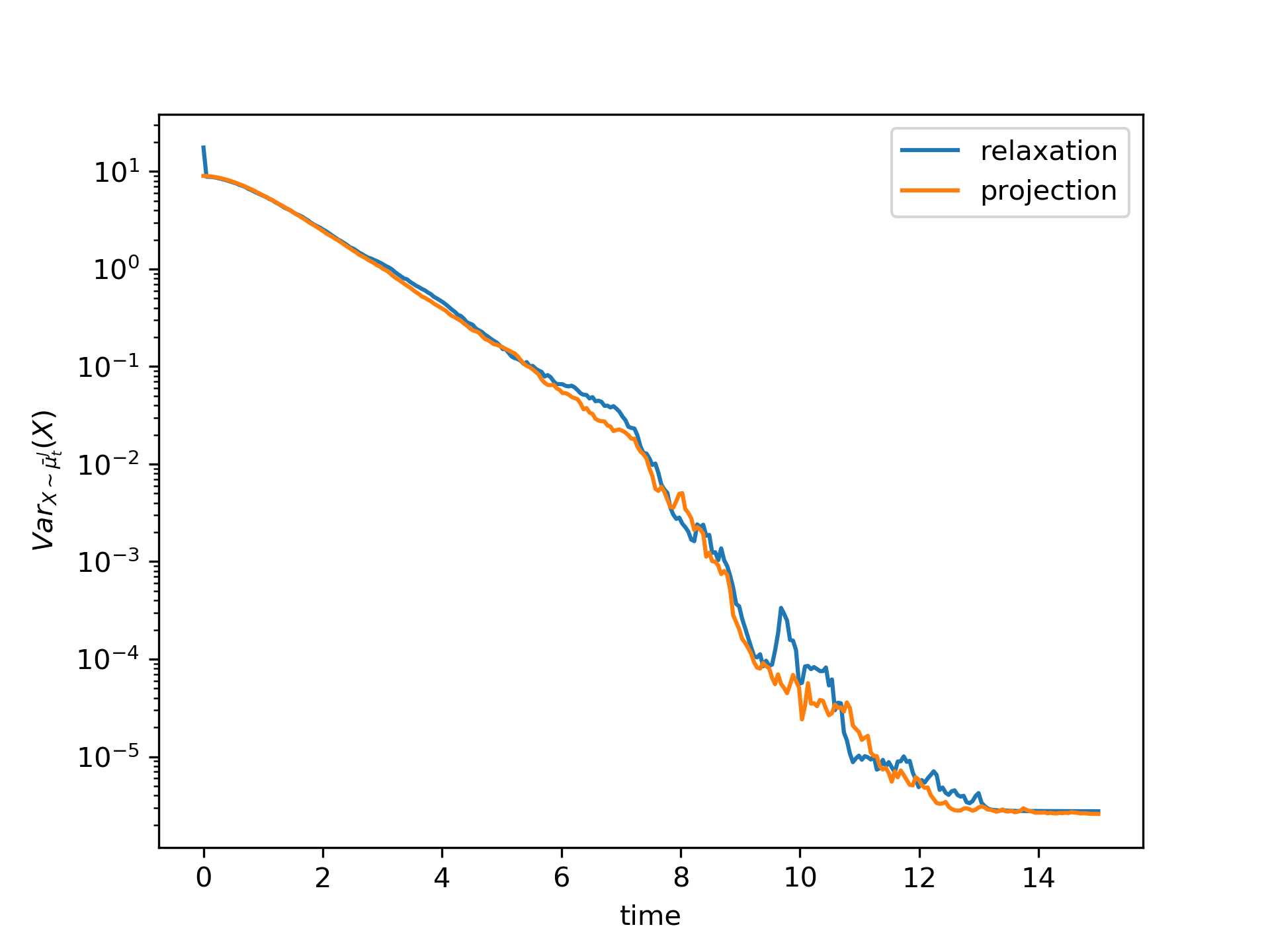}
    \includegraphics[width=0.32\linewidth]{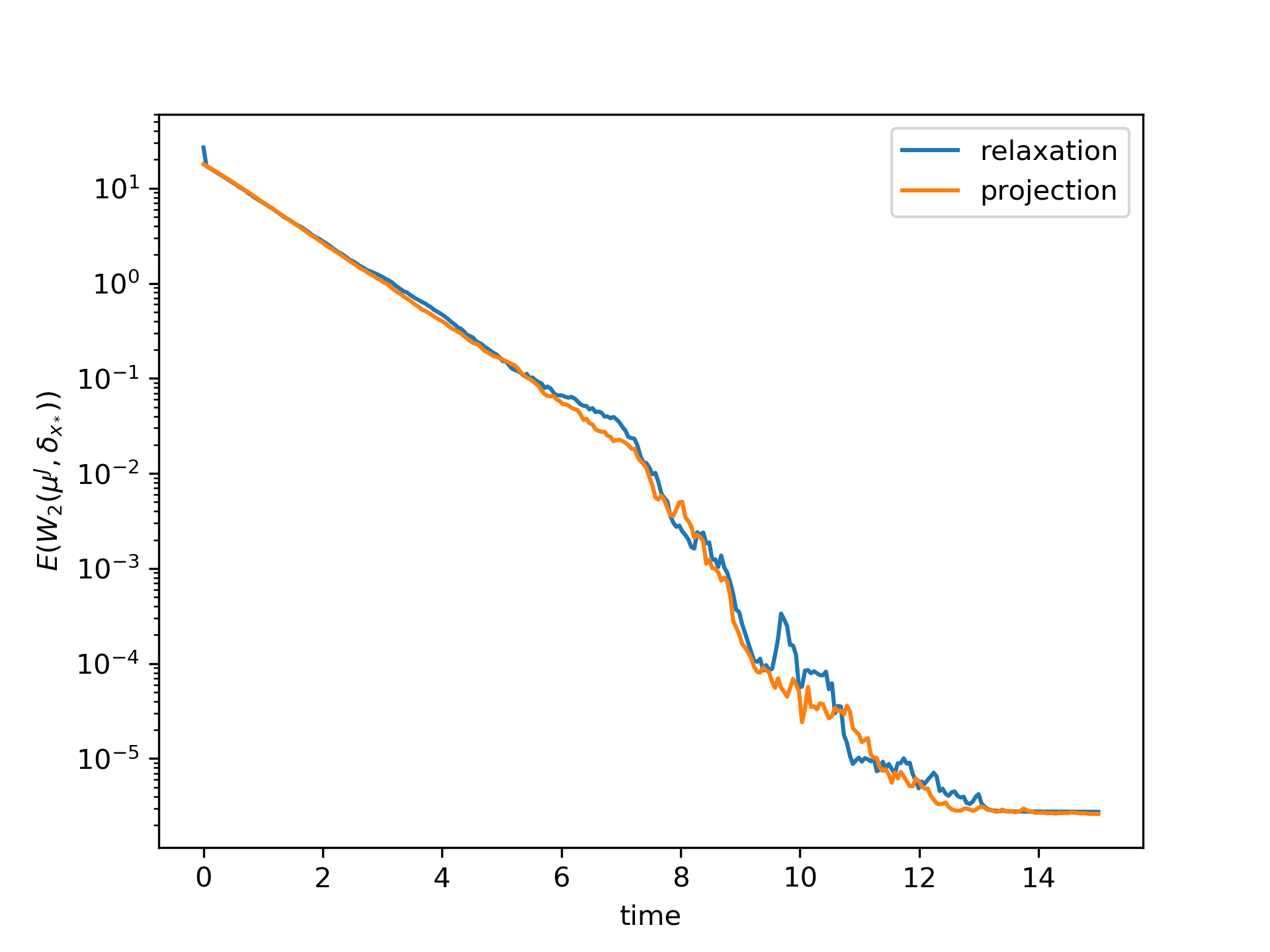}
    \includegraphics[width=0.32\linewidth]{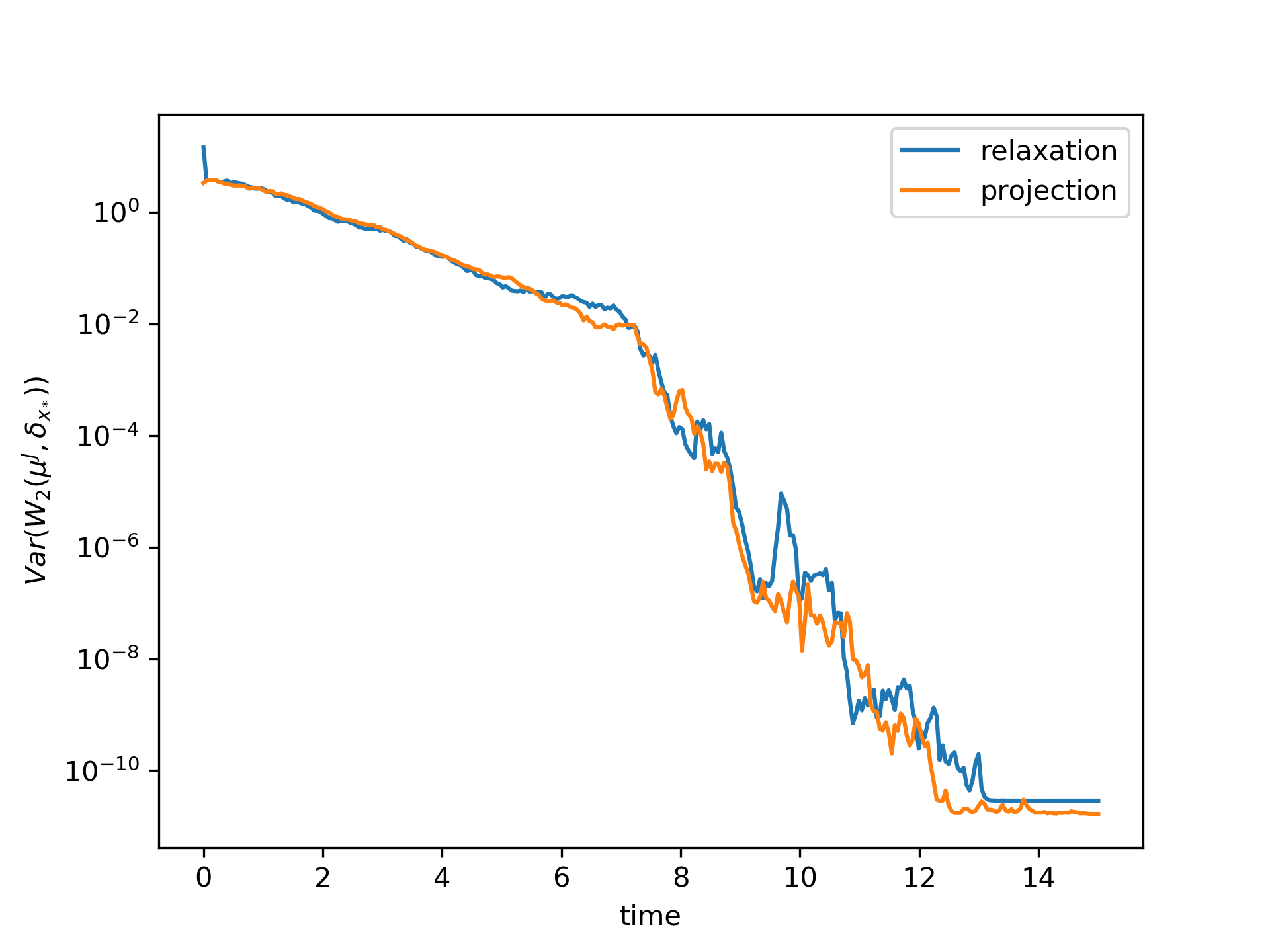}
    \caption{%
        Comparison of our method~\eqref{eq:cbo_with_constraints} with the projection-based method.
        The parameters employed for the simulation are the following: $J=50, M=100,\alpha = 30, \nu = 1, \Delta t = 0.0005, \sigma = 0.7, \varepsilon =0.1$.
    }
    \label{fig:projection}
\end{figure}
The results of the comparison are presented in~\cref{fig:projection}.
In the left panel, which depicts the evolution of the expected value of the variance,
we observe that the two graphs differ at the beginning,
which is not surprising given that the ensembles are not initially confined to the manifold with our method,
and then the two approaches perform comparably.
One advantage of our approach is that the initial ensemble can be drawn from any probability distribution over $\real^d$,
whereas for the projection method of~\cite{fornasier2020consensus} the initial ensemble needs to be supported on the constraint.
In addition, calculating projections to the feasible manifold at each time step can be computationally costly.
The panel in the  middle illustrates the evolution of the expected value of the Wasserstein distance to the minimizer,~$W_2(\mu^{50}, \delta_{x_*})$,
approximated for each method from 100 independent simulations and the evolution of the variance of the Wasserstein distance to the minimizer, again approximated for each method from 100 independent simulations, is shown on the right.
As expected, the ensembles corresponding to both methods eventually collapse and provide reasonable approximations of the global minimizer.

\begin{figure}
    \centering
    \includegraphics[width=0.75\linewidth]{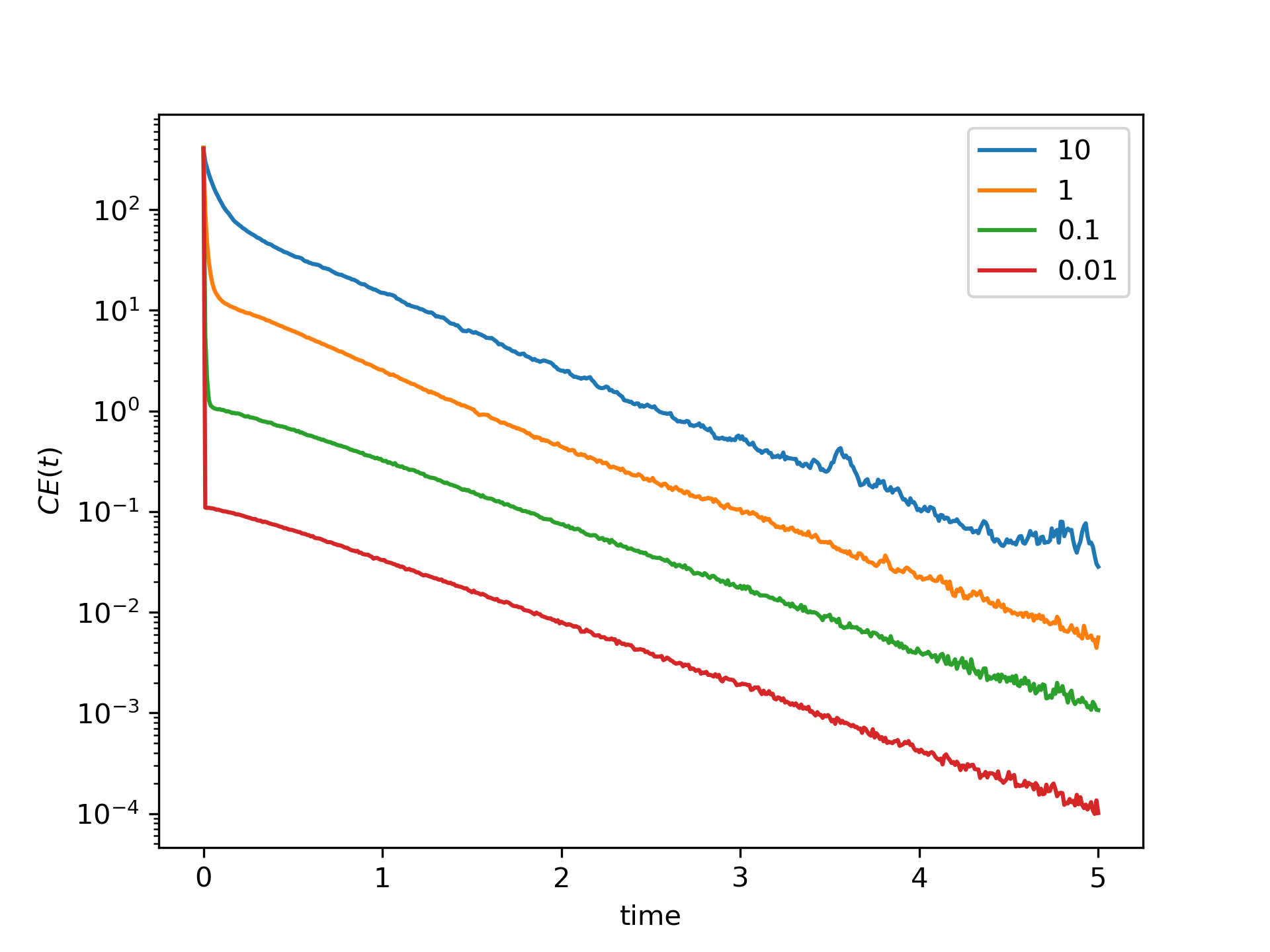}
    \caption{Distance to constraint for different relaxation parameters $\epsilon.$ Other parameters of the simulation are: $J=1000, M=100,\alpha = 30, \nu = 1, \Delta t = 0.0005, \sigma = 0.7.$ }
    \label{fig:relaxation}
\end{figure}
We notice that, in all the above plots concerning CBO,
the time evolution of the particle ensemble can be divided into roughly three phases.
In the first phase, at the beginning of the simulation,
the particles are quickly driven to the manifold.
In the second phase, the particles move along the feasible manifold towards the solution of the constrained optimization problem,
and few fluctuations are observed in the plots.
The size of the fluctuations, in logarithmic scale, in the final phase are consistent with the fact that the amplitude of the noise decreases as we get closer to the collapse of the ensemble near the global minimizer at the final stages of the simulation.


To conclude this section, we examine the influence of the parameter $\varepsilon$, which enters in the relaxation drift,
on the speed of relaxation towards the constraint.
This is illustrated in~\cref{fig:relaxation},
where the evolution of ``constraint energy'',
measured using~\eqref{eq:CE},
is shown for different values of $\varepsilon$.
For all the values of $\varepsilon$ considered,
the ``constraint energy'' is approximated based on $M=100$ independent simulations with $J=1000$ particles each.

\subsection{Numerical experiments for the CBO particle system}%
\label{sub:numerical_experiments_for_the_particle_systems}

In this section, we present additional numerical experiments concerning the particle system.
We investigate, in particular, the influence of the parameters $\nu$, $\varepsilon$ and $J$ on the accuracy of the convergence point of CBO,
first for the following optimization problem:
\begin{equation}
    \label{eq:optimization_problem_ackley}
    \argmin_{x \in B} f_A(x - x_*),
    \quad x_{*} = \begin{pmatrix} 2 \\ 2 \end{pmatrix},
    \quad B = \bigl\{ x \in \real^2 \colon x^2 + y^2 = 18 \bigr\}.
\end{equation}
All the numerical results presented below are generated from CBO with the parameters $\sigma = 0.7$ and $\Delta t = 0.01$.
In each simulation, the particles forming the initial ensembles are always drawn from $\mathcal N(0, 100I_2)$.
The other parameters are specified on a case-by-case basis.

\begin{figure}[ht]
    \centering
    \includegraphics[width=.24\linewidth]{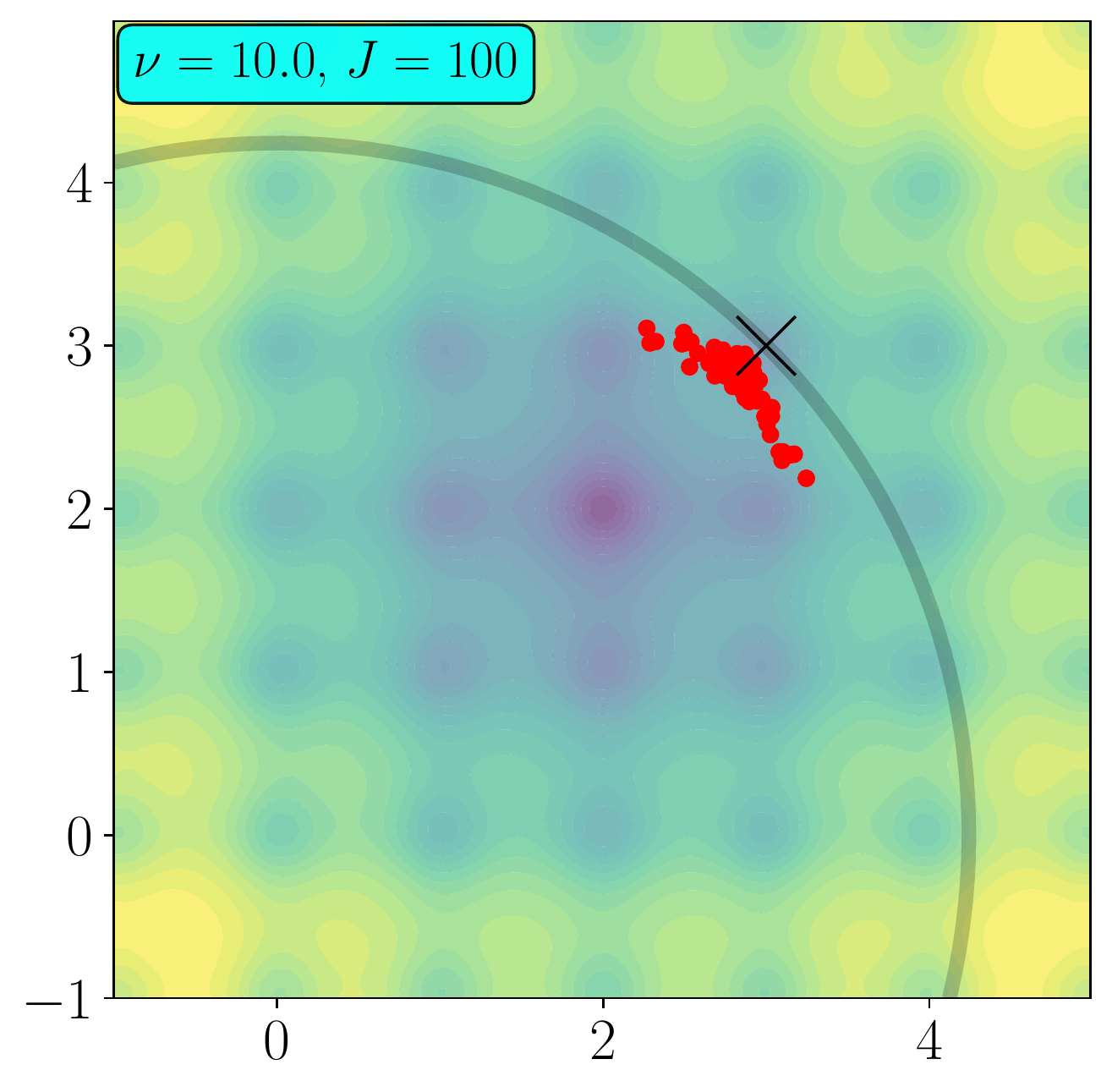}
    \includegraphics[width=.24\linewidth]{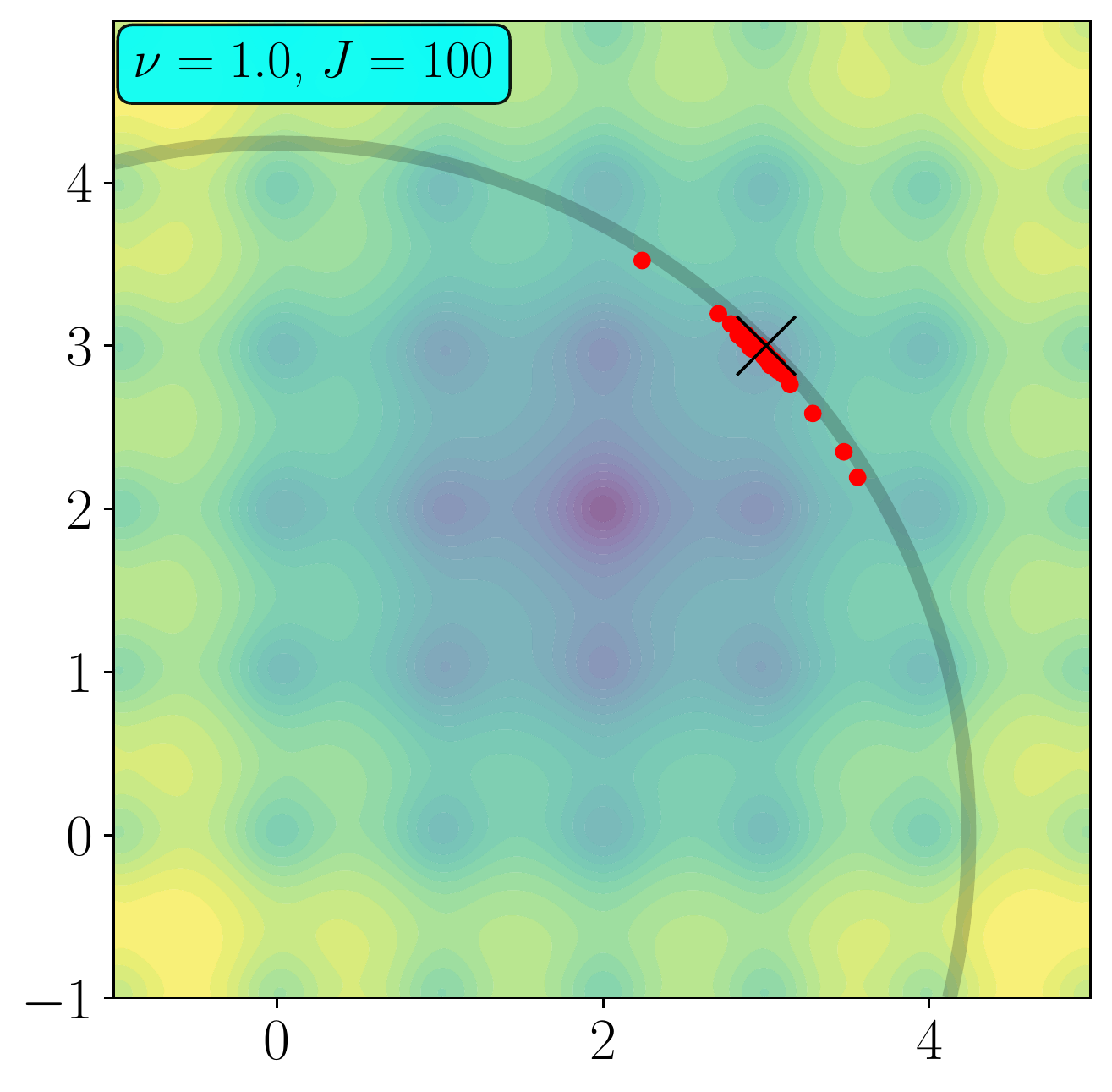}
    \includegraphics[width=.24\linewidth]{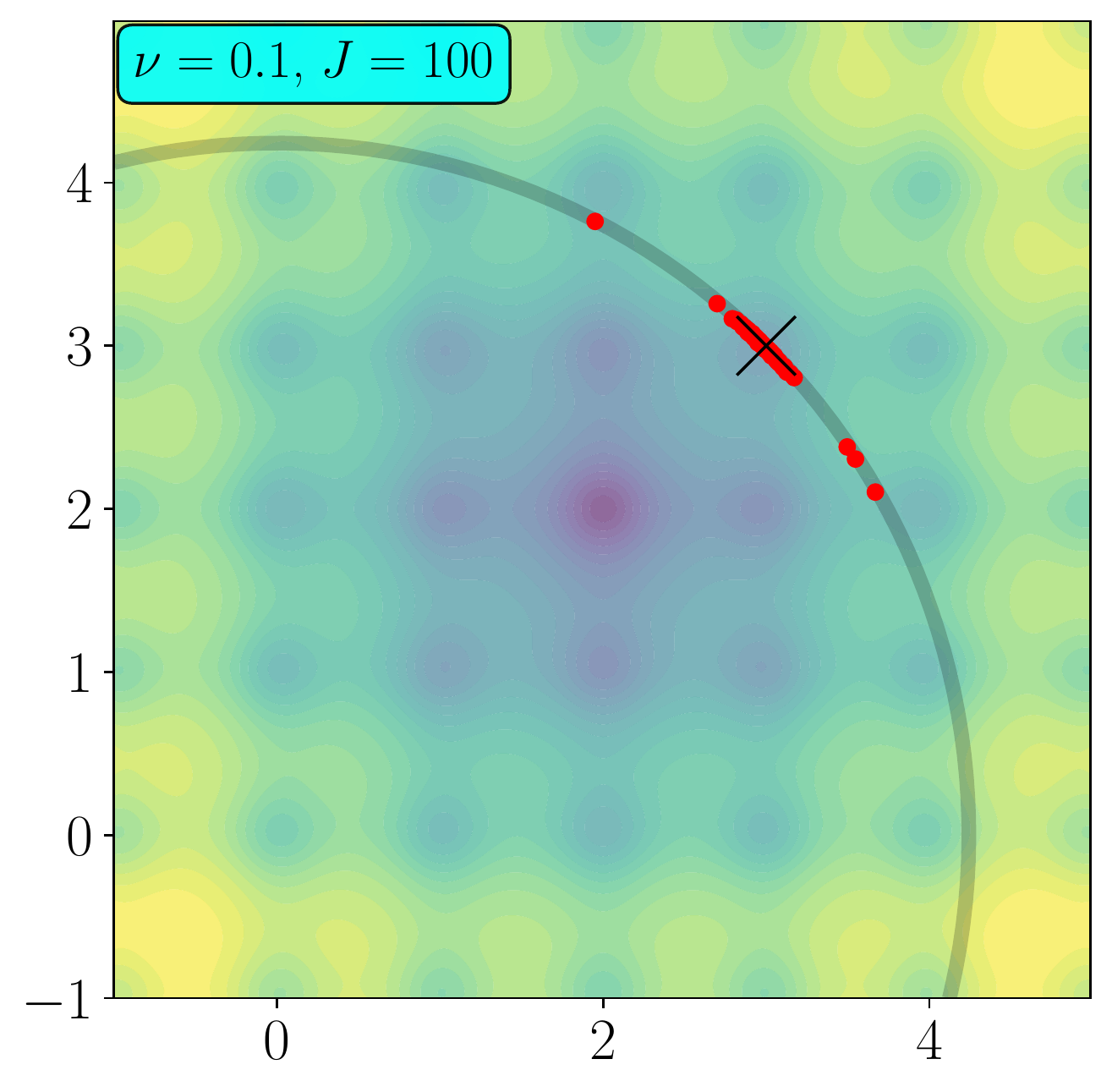}
    \includegraphics[width=.24\linewidth]{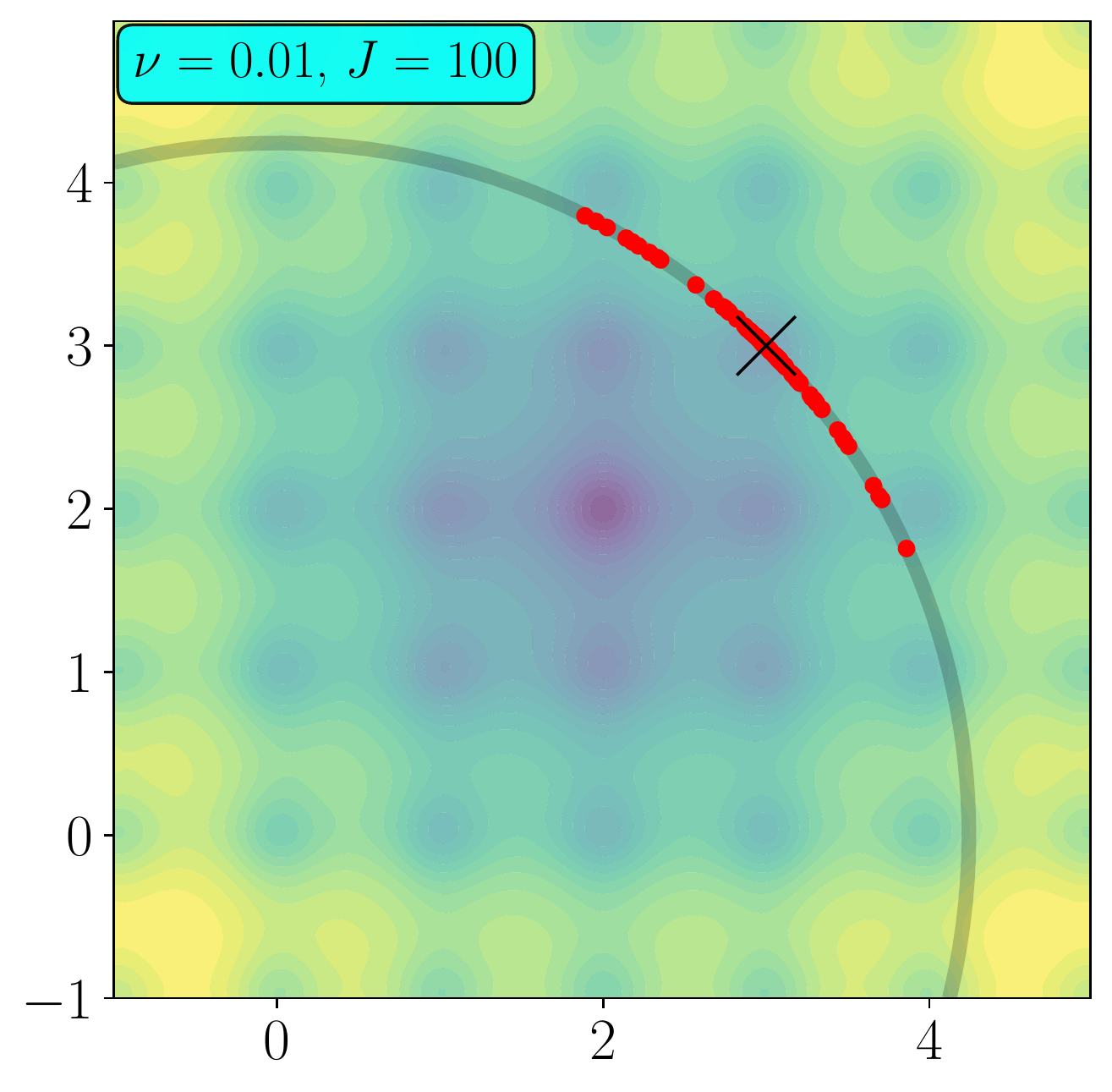}

    \includegraphics[width=.24\linewidth]{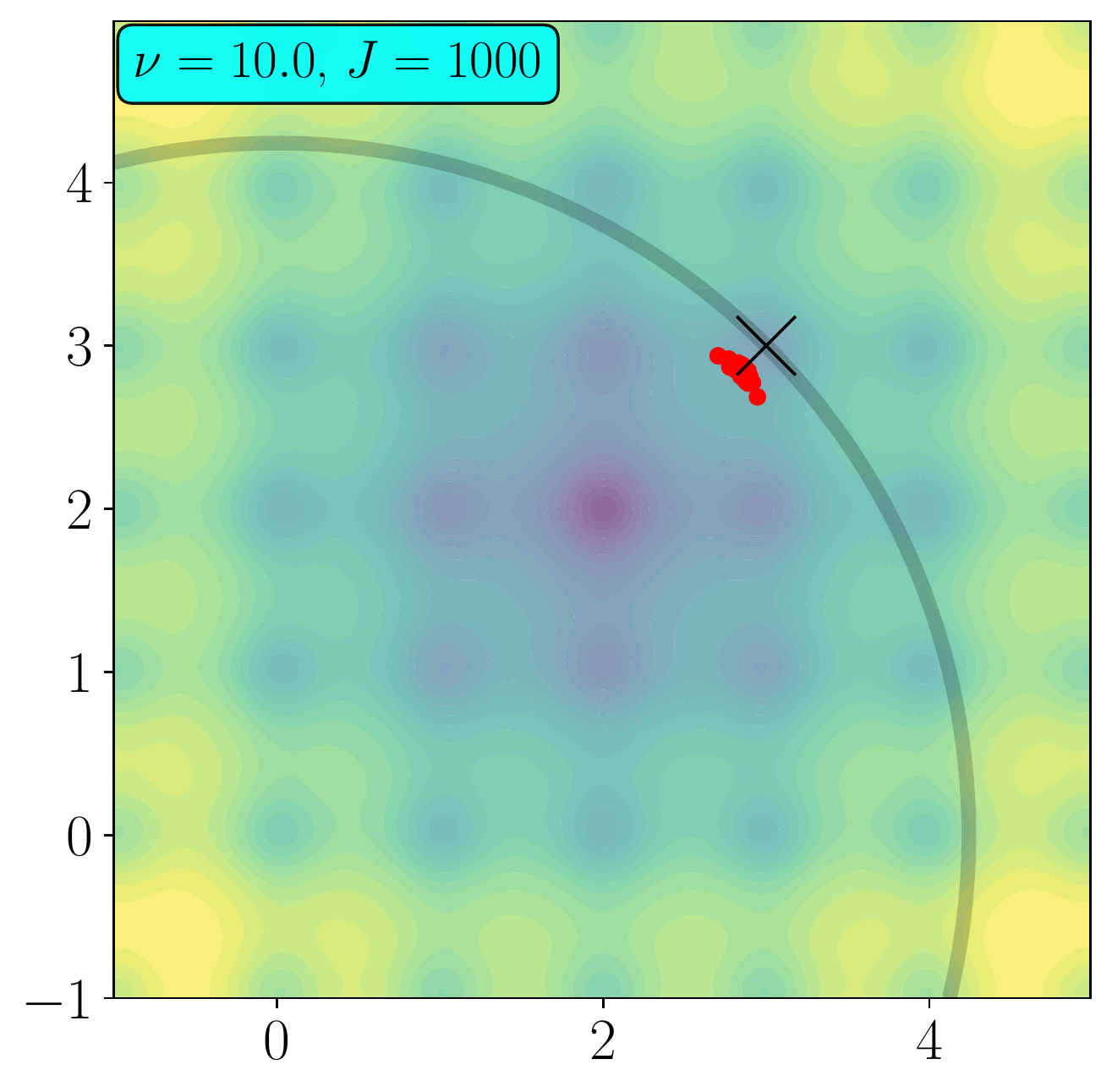}
    \includegraphics[width=.24\linewidth]{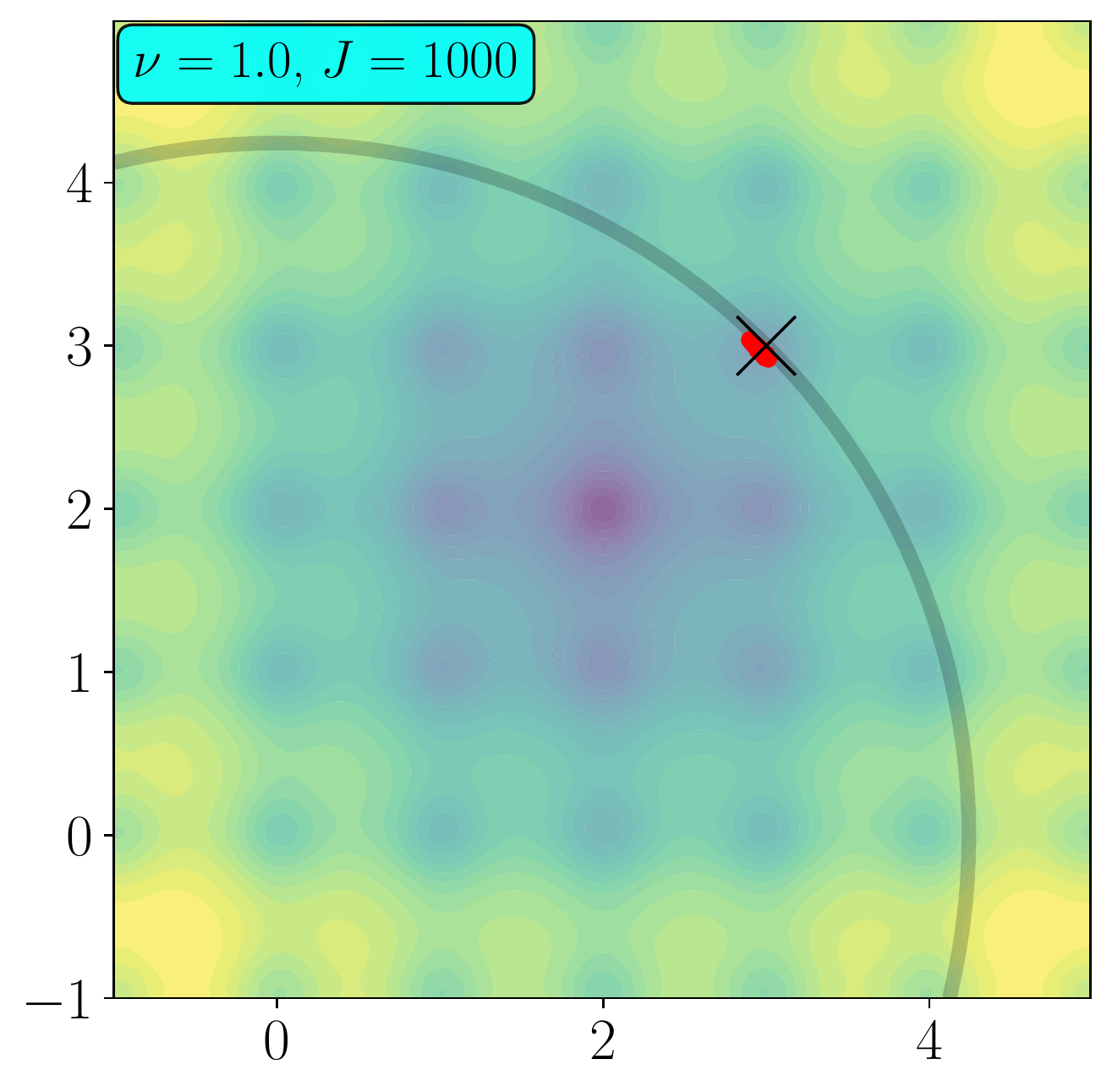}
    \includegraphics[width=.24\linewidth]{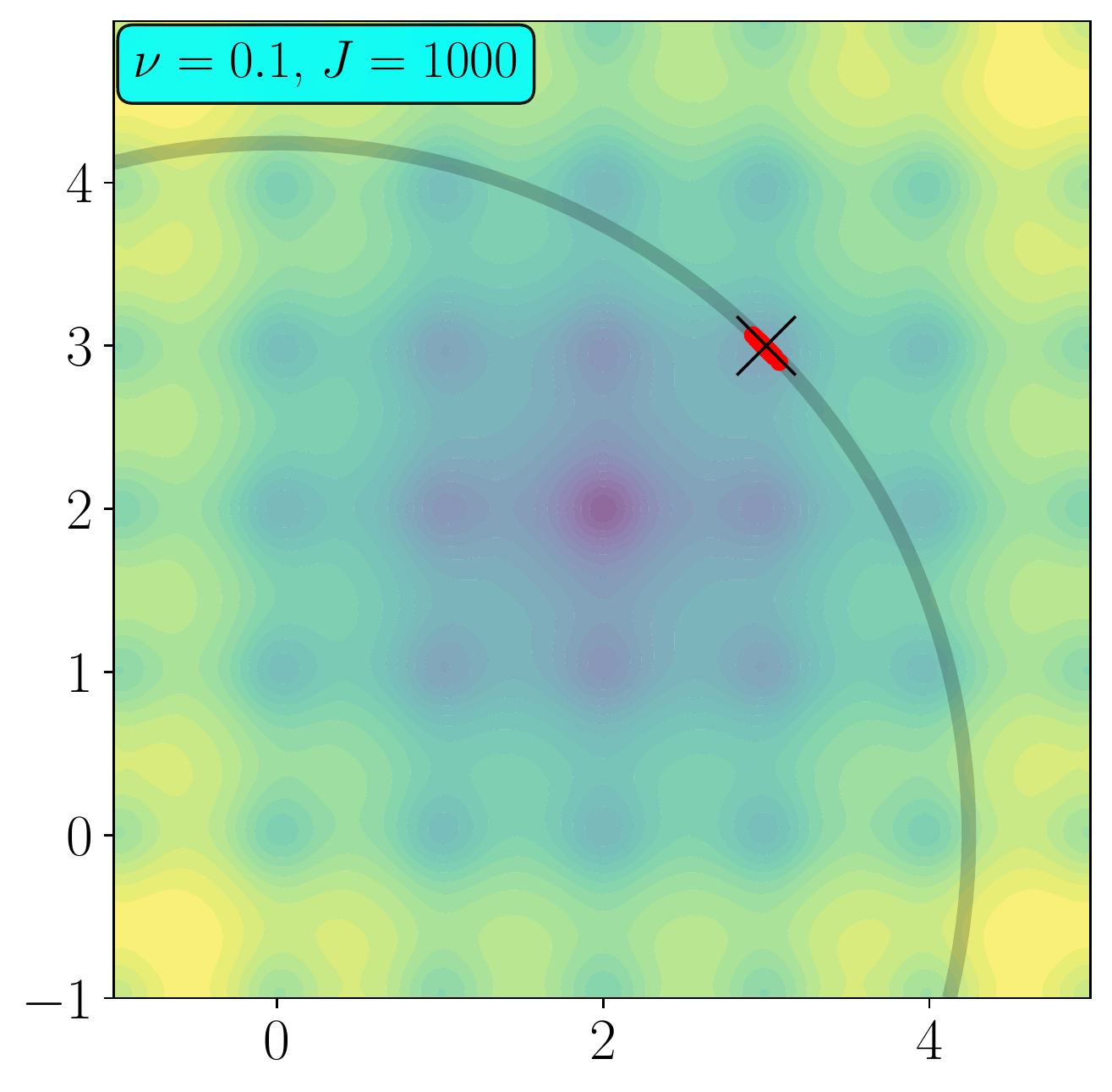}
    \includegraphics[width=.24\linewidth]{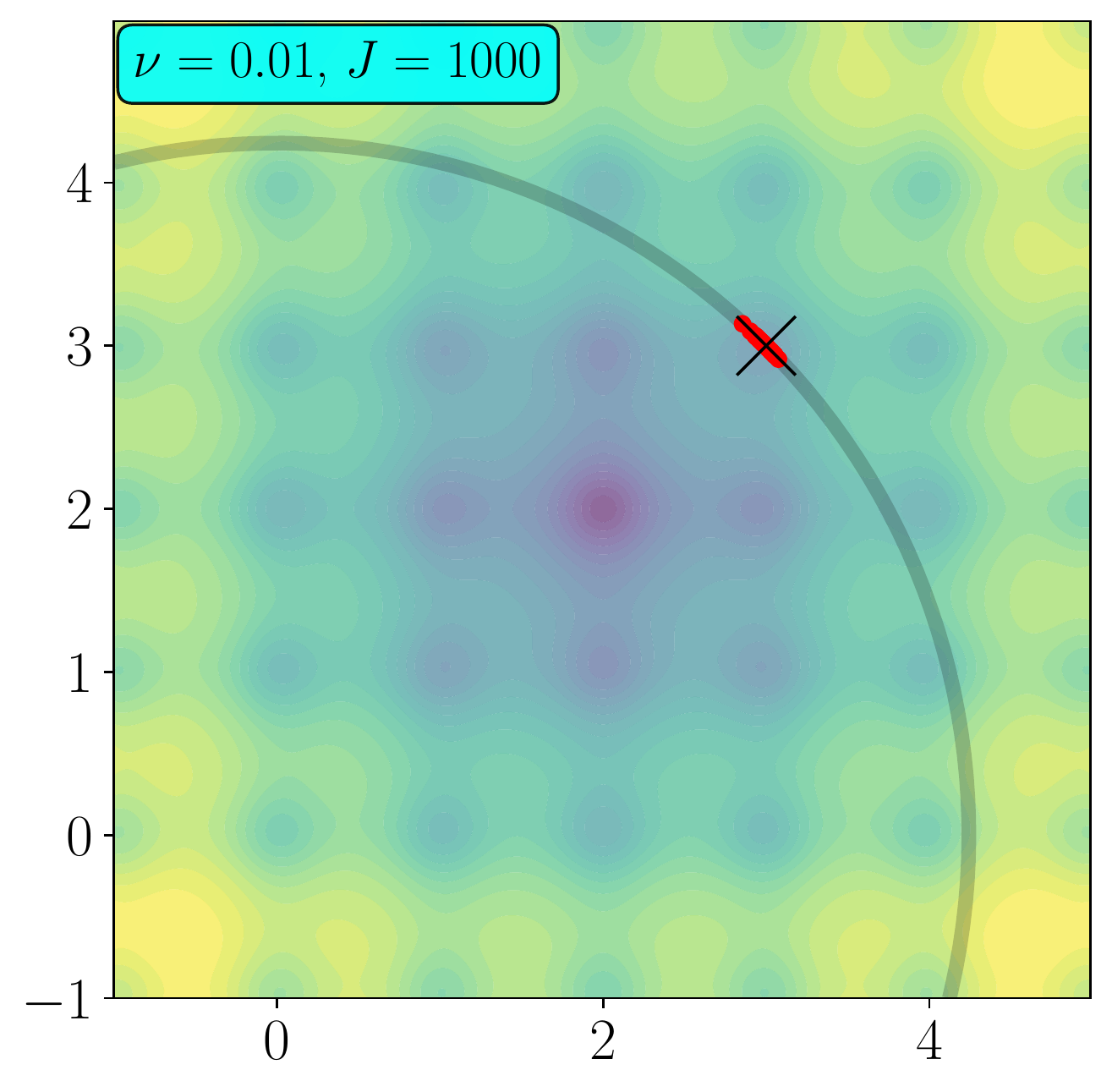}
    \caption{%
        Points of convergence of $M=100$ simulations of CBO for the constrained optimization problem~\eqref{eq:optimization_problem_ackley},
        for different values of $\nu$ and $J$ and for $\varepsilon = 0$.
        The black cross indicates the position of the global minimizer under the constraint.
    }
    \label{figures:cbo_particles}
\end{figure}

In~\cref{figures:cbo_particles},
we illustrate for different values of $\nu \in \{10, 1, .1, .01\}$ and $J \in \{100, 1000\}$,
the convergence points of 100 independent simulations (per set of parameters).
Here we take $\varepsilon = \infty$;
that is, we employ the CBO particle system with penalization but without relaxation drift.
From the figures corresponding to $J = 100$,
we observe that,
although smaller values of~$\nu$ enable to enforce the constraint more effectively,
taking $\nu$ too small is detrimental for the accuracy of the convergence point.
This effect is observed also for $J = 1000$,
but to a much lesser extent,
which may indicate a better behavior of the system in the many-particle limit.
This numerical experiment highlights the considerable influence of the number of particles on the behavior of the method.
In practice, a trade-off needs to be achieved between precision and computational cost.

\begin{figure}[ht]
    \centering
    \includegraphics[width=.24\linewidth]{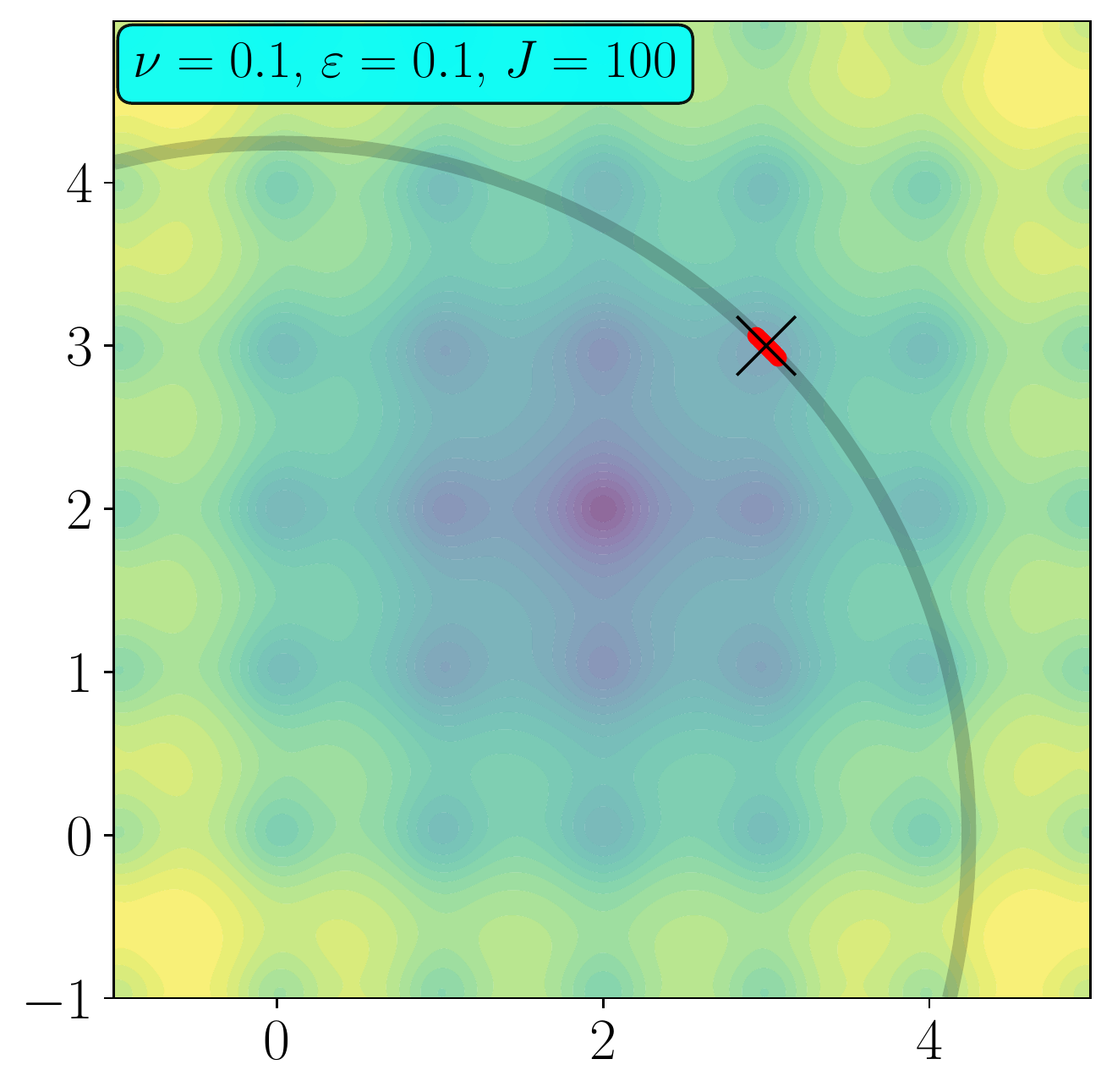}
    \includegraphics[width=.24\linewidth]{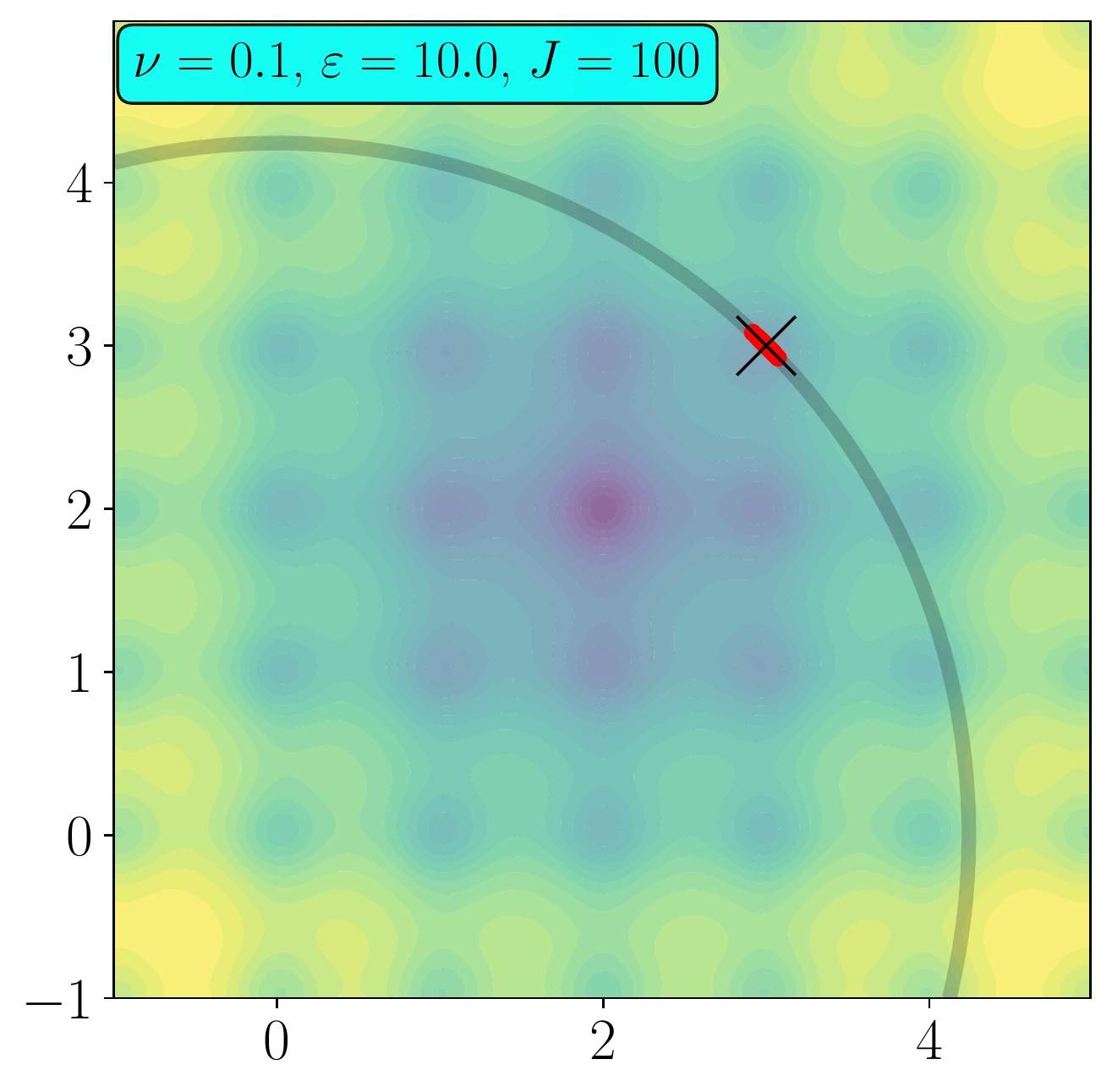}
    \includegraphics[width=.24\linewidth]{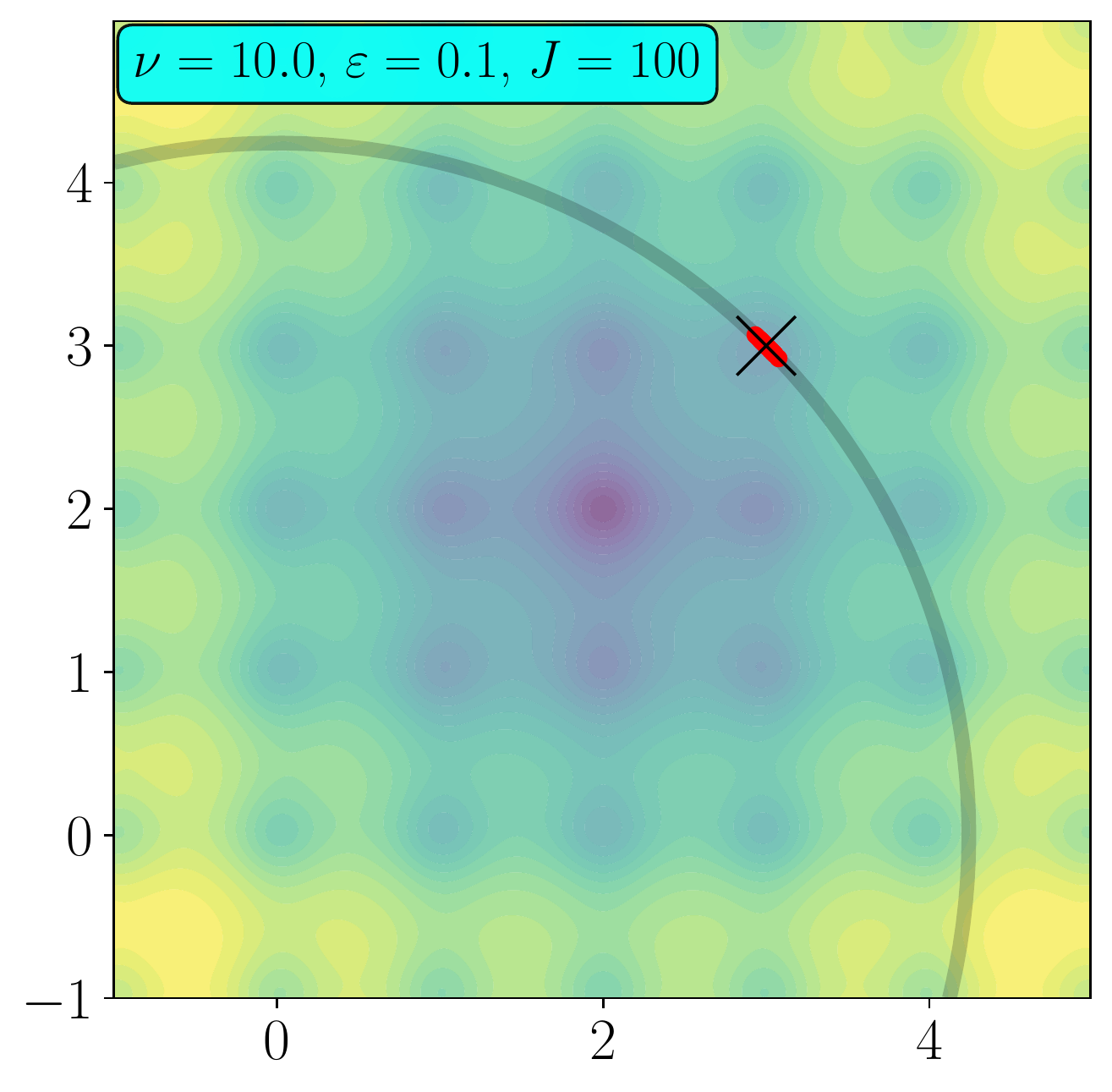}
    \includegraphics[width=.24\linewidth]{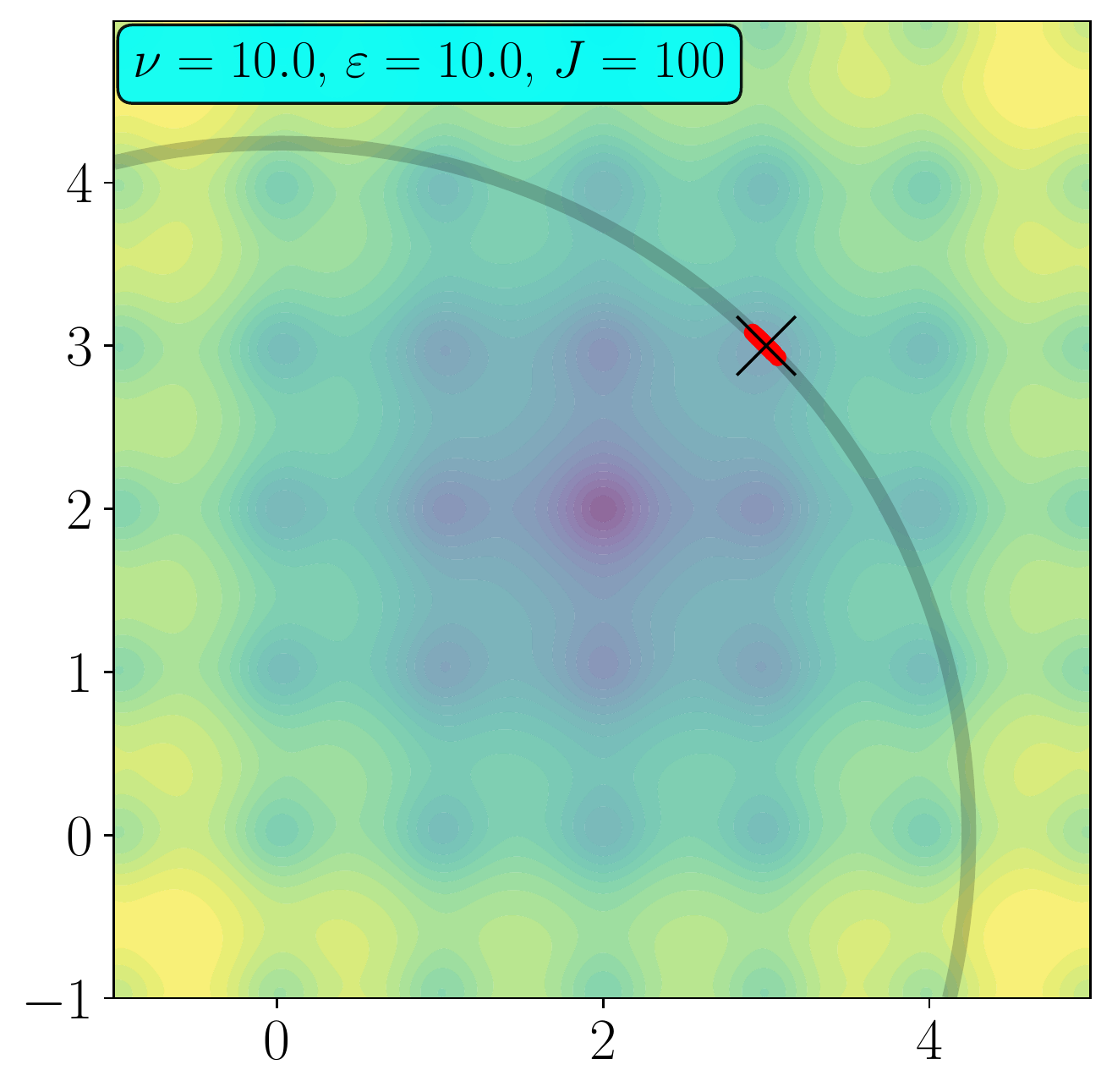}
    \caption{%
        Points of convergence of $M=100$ simulations of CBO for the constrained optimization problem~\eqref{eq:optimization_problem_ackley},
        for different values of $\nu$ and $\varepsilon$.
    }
    \label{figures:cbo_particles_2}
\end{figure}

In~\cref{figures:cbo_particles_2},
we illustrate the convergence points of 100 independent simulations with this time $(\nu, \varepsilon) \in \{.1, 10\}^2$,
and for fixed $J = 100$.

For the simulations corresponding to $\varepsilon = .1$,
we take a smaller time step $\Delta t = 10^{-3}$ in order to avoid stability issues,
and in all the simulations, we employ the semi-implicit scheme~\eqref{eq:semi-implicit}.
Interestingly, the presence of the relaxation drift,
even with a relatively small amplitude $1/\varepsilon$,
considerably improves the performance of the method compared to the case without penalization.

\begin{figure}[ht]
    \centering
    \includegraphics[width=.32\linewidth]{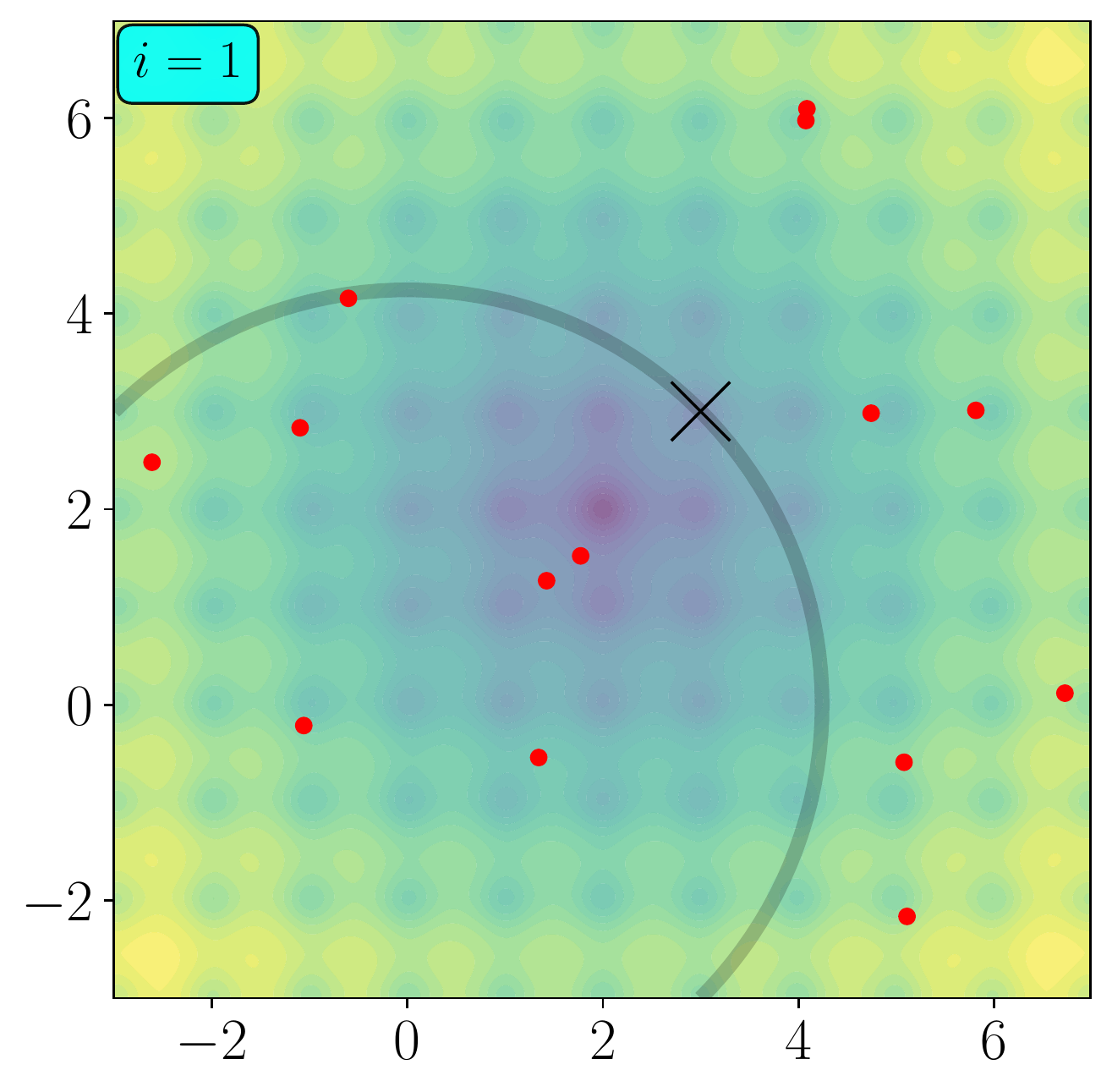}
    \includegraphics[width=.32\linewidth]{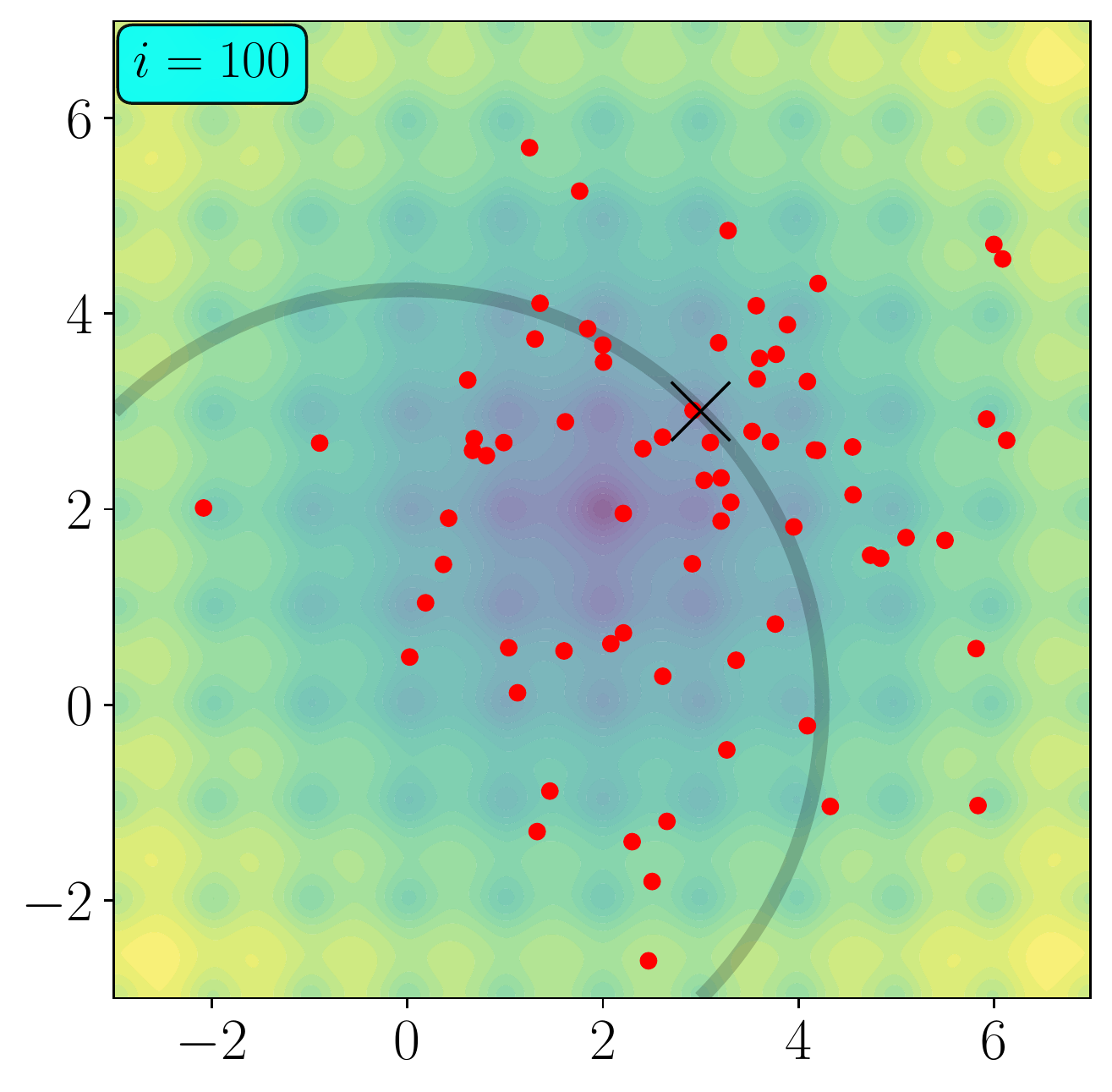}
    \includegraphics[width=.32\linewidth]{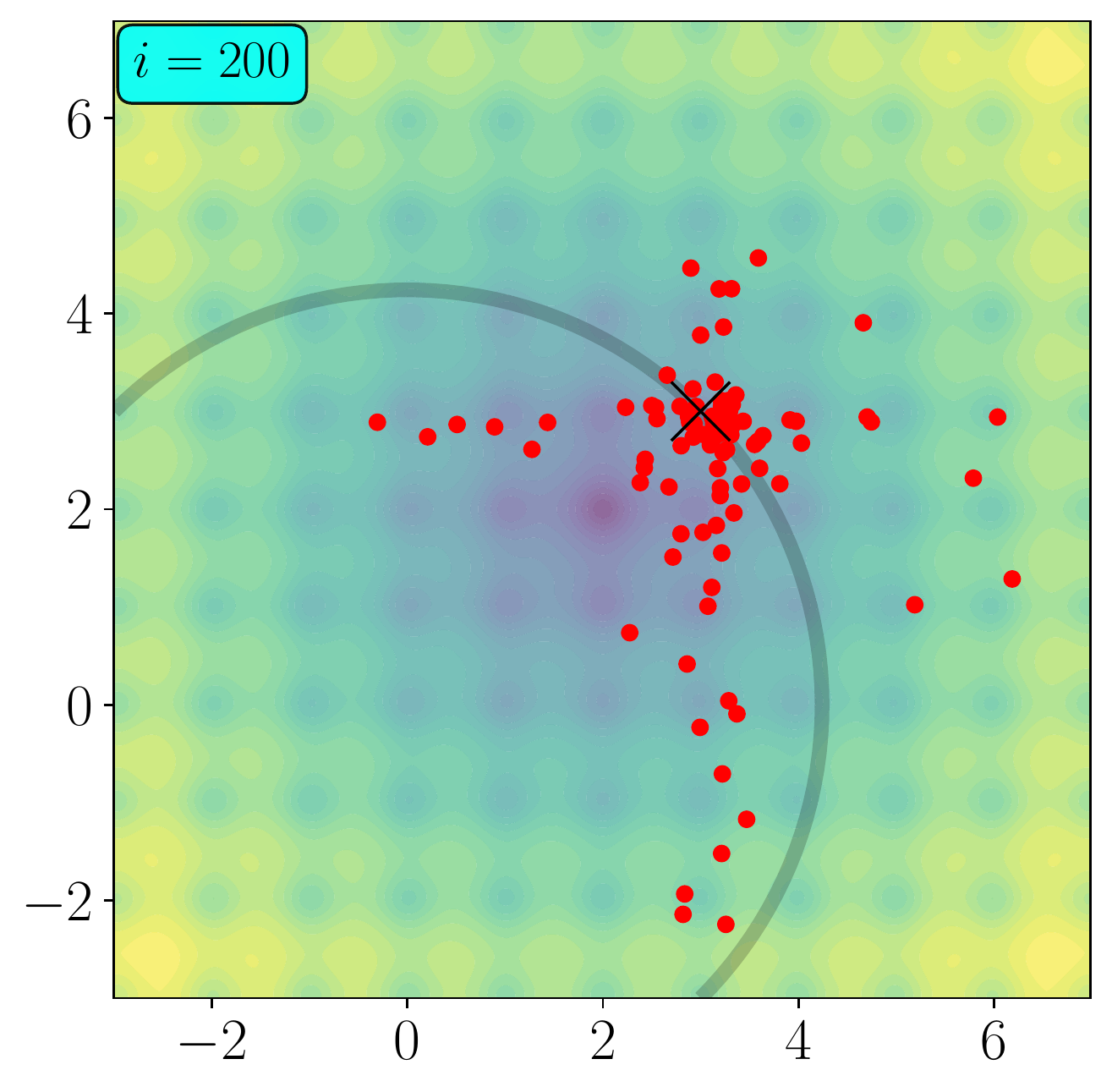}
    \includegraphics[width=.32\linewidth]{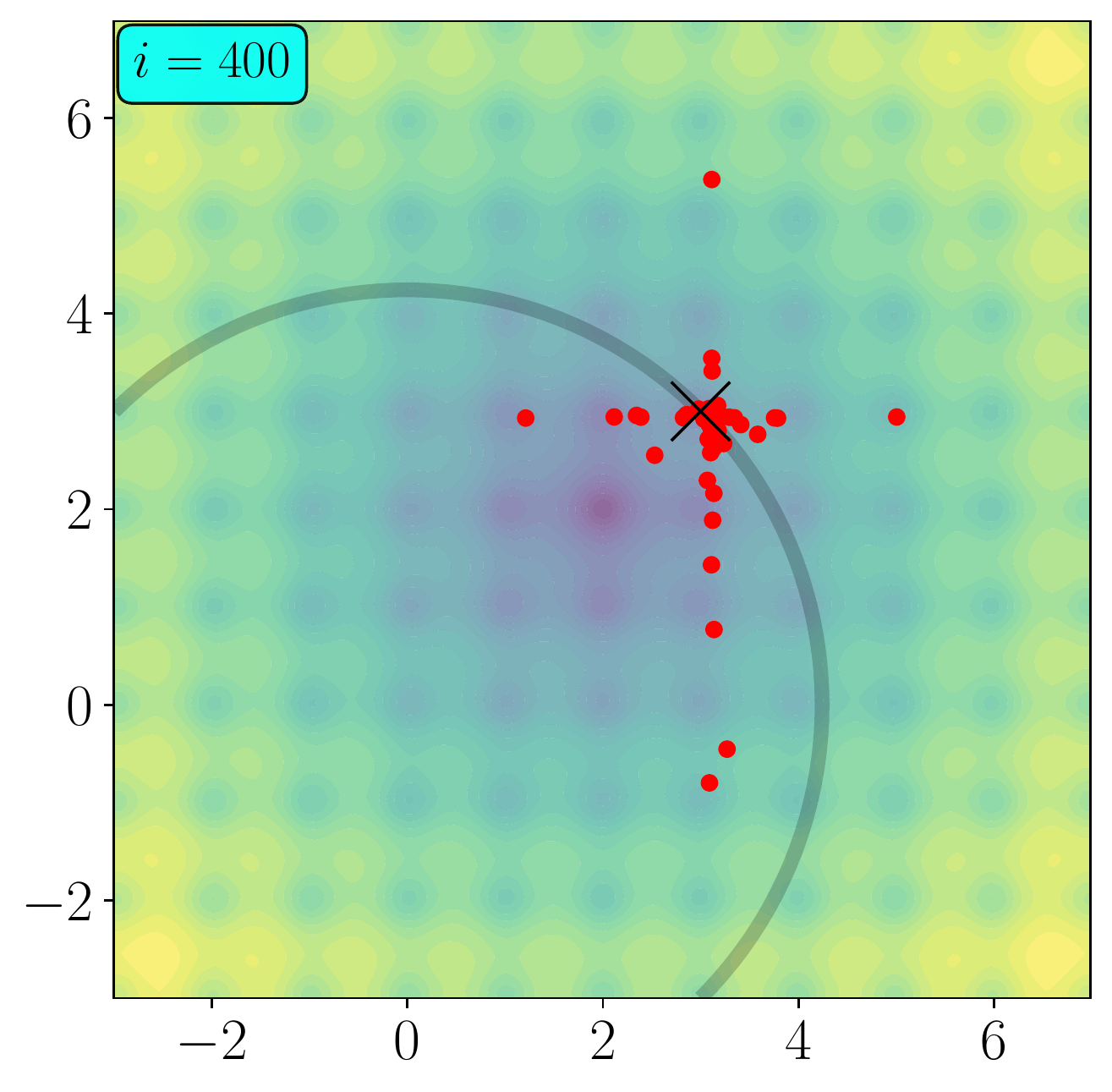}
    \includegraphics[width=.32\linewidth]{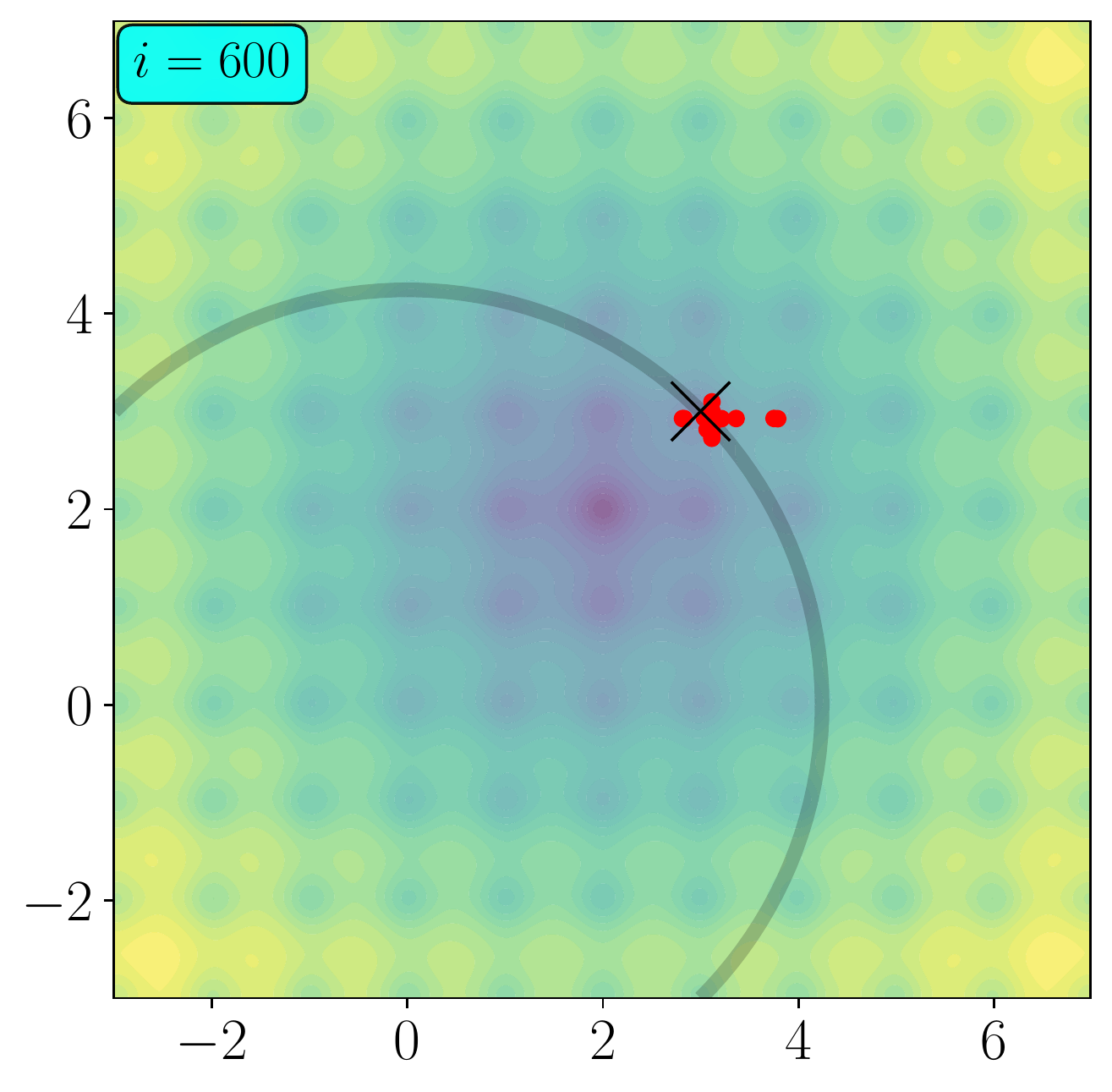}
    \includegraphics[width=.32\linewidth]{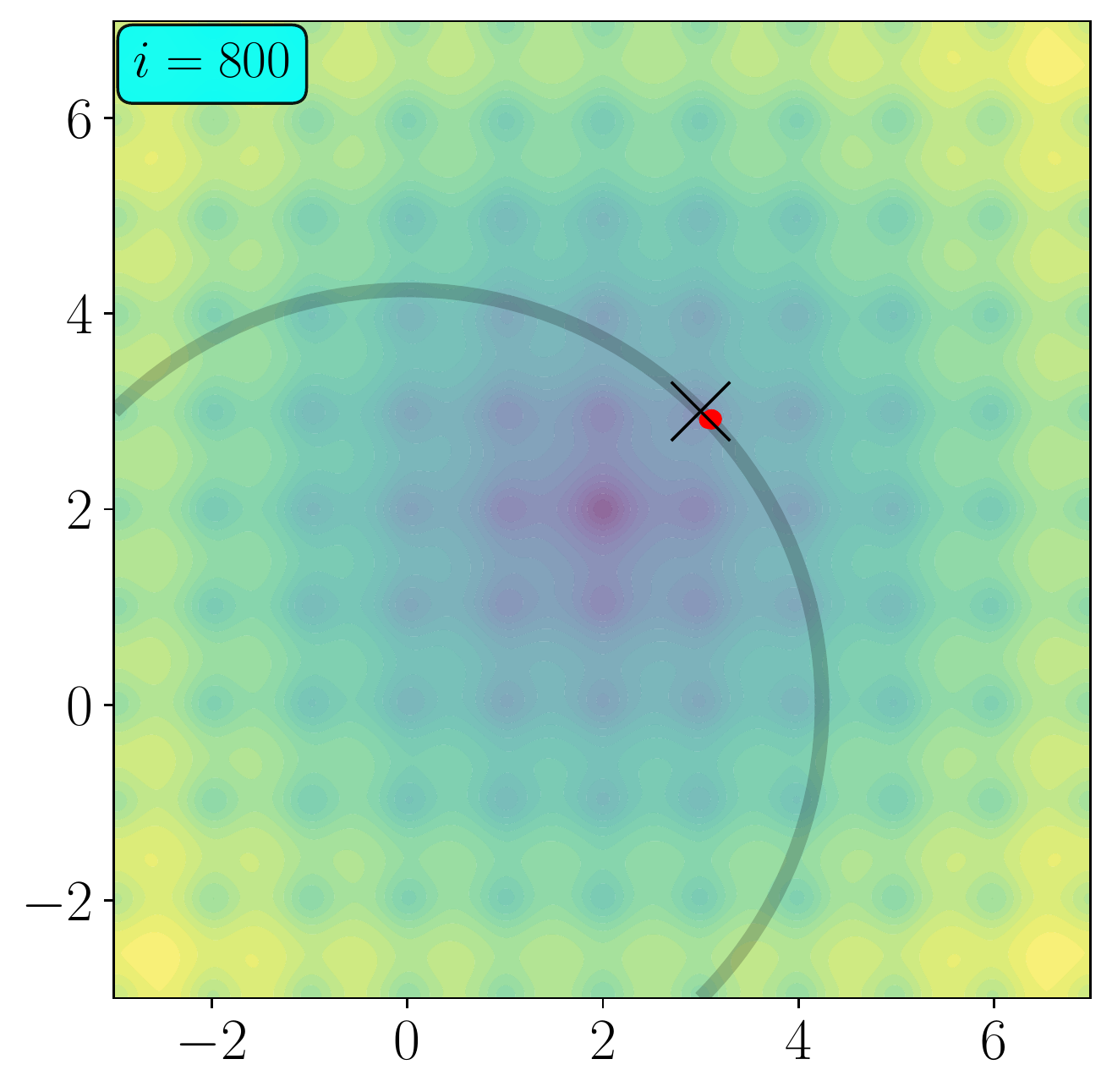}
    \caption{%
        Evolution of the particles for simulations when CBO is employed for minimizing the Ackley function under the constraint that $\bigl\{x^2 + y^2 \geq 18 \bigr\}$.
        The parameters for this example are $\varepsilon = \nu = 1$ and $J = 100$.
    }
    \label{figures:cbo_particles_1000}
\end{figure}

To conclude this section,
we present~ in~\cref{figures:cbo_particles_1000} numerical results for~\eqref{eq:optimization_problem_ackley}
with the inequality constraint $B = \bigl\{ x \in \real^2 \colon x^2 + y^2 \geq 18 \bigr\}$
instead of an equality constraint.
The problem setting is then the same as in~\cref{fig:pdeSdeConstraint}.
Rather than presenting the points of convergence of the method for several parameter choices as in the previous figures,
here we perform only one simulation with a given set of parameters ($\nu = \varepsilon = 1$, $J = 100$),
and we illustrate the evolution in time of the ensemble.
As the figure shows,
the ensemble appears to collapse after roughly $800$ iterations,
that is at time $t \approx 8$.
In addition, the point of convergence is very close to the optimal solution given the constraint.

\section{Numerical results for EKI with constraints}
\label{sec:numerics_for_eki}
In this section,
we present numerical results for EKI with constraints.
We simulate only the particle system and do not present result for the associated mean-field equation.
Our numerical experiments aim at illustrating the performance of the method on a toy problem,
as a proof of concept.

The numerical schemes employed are described in~\cref{sub:numerical_scheme_eki},
and numerical results for a simple inverse problem arising from a PDE application are presented in~\cref{sub:eki}.

\subsection{Numerical schemes}%
\label{sub:numerical_scheme_eki}

The presence of the relaxation term originating from the constraint poses a challenge
also for the discretization of EKI with constraints~\eqref{eq:eki:evolution_of_thetas_constraints},
although to a lesser extent than for CBO
thanks to the affine-invariance of the method~\cite{garbuno2020affine}.
In this section,
we consider specifically the derivative-free dynamics~\eqref{eq:unified_eki},
and we investigate two different approaches for discretizing this dynamics in time:
\begin{itemize}
    \item An explicit discretization using the explicit Euler method with adaptive time step.
        Specifically, the particle positions are evolved according to
        \begin{align}
            \label{eq:unified_eki_explicit}
            \begin{aligned}[b]
                x^{(j)}_{n+1} &= x^{(j)}_{n} - \Delta t_n \sum_{k=1}^{J} M^{kj}_n x^{(k)}_{n}, \qquad j  = 1, \dotsc, J, \\
                M^{kj}_n &= \frac{1}{J} \eip*{\mathcal G(x^{(k)}_n) - \bar {\mathcal G}_n}{\mathcal G(x^{(j)}_n) - \widetilde y}[\widetilde \Gamma],
            \end{aligned}
        \end{align}
        and the time step is calculated dynamically according to the method proposed in~\cite{MR3998631},
        that is
        \begin{equation}
            \label{eq:time_step_adaptation}
            \Delta t_n = \frac{\Delta t_*}{\lVert M_n \rVert_2 + \frac{\Delta t_*}{\Delta t_{\max}}},
        \end{equation}
        where $\Delta t_*$ (base time step) and $\Delta t_{\max}$ (maximum time step) are parameters,
        and where $M_n$ is the matrix with entries $M_n^{kj}$.
        With the notation $X_n = (x^{(1)}_n, \dotsc, x^{(J)}_n)$,
        equation~\eqref{eq:unified_eki_explicit} reads in matrix form as
        \begin{equation}
            \label{eq:update_explicit}
            X_{n+1} = X_n - \Delta t_n X_{n} M_n.
        \end{equation}
    \item A semi-implicit discretization given by
        \begin{equation*}
            \label{eq:semi_implicit}
            X_{n+1} = X_n - \Delta t_n X_{n+1} M_n,
        \end{equation*}
        where the same notation is used as in the previous item.
        The associated update formula is given by
        \begin{equation}
            \label{eq:update_semi_implicit}
            X_{n+1} = X_n (I_d + \Delta t_n M_n)^{-1}.
        \end{equation}
        Although this approach could be employed with a fixed time step,
        we obtain faster convergence when the time step is adapted according to~\eqref{eq:time_step_adaptation},
        and so we consider only the latter setting.
\end{itemize}
\begin{remark}
    Note that $\sum_{k=1}^{K} M_n^{kj} = 0$ for all $n$ and $j$, by definition of~$\bar {\mathcal G}_n$.
    Consequently, the update formulas~\eqref{eq:update_explicit} and~\eqref{eq:update_semi_implicit} can be rewritten equivalently as
    \begin{subequations}
    \begin{align}
        \label{eq:good_explicit}
        X_{n+1} &= X_n - \Delta t_n (X_{n} - \bar x_n) M_n, \\
        \label{eq:good_semi_implicit}
        X_{n+1} &= \bar x_n + (X_n - \bar x_n) (I_d + \Delta t_n M_n)^{-1},
    \end{align}
    \end{subequations}
    where
    \(
        \bar x_n = \frac{1}{J} \sum_{j=1}^{J} x^{(j)}_n
    \)
    and \( X_n - \bar x_n \) denotes,
    by a slight abuse of notation,
    the matrix obtained by subtracting $\bar x_n$ from each column of $X_n$.
    Although mathematically equivalent,
    these discretizations are observed in our numerical experiments to be much less sensitive to roundoff errors,
    and therefore preferable in practice.
\end{remark}

\subsection{Numerical experiments for EKI with constraints}%
\label{sub:eki}
In this section,
we illustrate the performance of the dynamics~\eqref{eq:unified_eki},
discretized according to~\eqref{eq:good_explicit} or~\eqref{eq:good_semi_implicit},
for solving a toy inverse problem with constraints.
The problem considered concerns the recovery of the initial condition of a Fokker--Planck equation based on incomplete observation of the solution.
More precisely, our aim is to find the parameters of a Gaussian mixture
\begin{equation}
    \label{eq:gaussian_mixture}
    \rho_0(x) = \sum_{n=1}^{N} w_n \frac{\displaystyle \exp \left( - \frac{\bigl\lvert x-m_n \bigr\rvert^2}{2 \sigma_n^2} \right)}{\sqrt{2 \pi \sigma_n^2}},
\end{equation}
given noisy observations at time $T$ of the solution $\rho(x, t)$ to the following initial value problem with initial condition $\rho_0$:
\begin{equation}
    \label{eq:initial_value_problem_fokker_planck}
    \left\{
    \begin{aligned}
        \partial_t \rho &= \partial_x (x \rho + \partial_x \rho),
        \qquad &&\text{in } \real \times (0, T), \\
        \rho &= \rho_0,
        \qquad &&\text{on } \real \times \{t = 0\}.
    \end{aligned}
    \right.
\end{equation}
We assume that the data is composed of noisy observations
\[
    y_k = \rho(x_k, T) + \eta_k,
    \qquad 1 \leq i \leq K,
\]
where $x_k \in \real$ are discrete positions
and $\eta_k$ are noise terms drawn independently from~$\mathcal N(0, \gamma^2)$ for some covariance $\gamma^2$.
The observation positions are given by $x_k = -L + (k-1) \frac{L}{K-1}$,
i.e.~they are uniformly spread between $-L$ and $L$,
with both ends included.

\paragraph{Recovering the weights.}%
Assuming first that the means and variances of the Gaussian mixture~\eqref{eq:gaussian_mixture} are known,
we seek to find,
as an approximation of the true weights $w^\dagger$,
a vector of weights~$w = (w_1, \dotsc, w_N)^\t$ that minimizes the least-square misfit
\[
    f(w)
    = \sum_{k=1}^{K} \frac{1}{\gamma^2}
    \bigl\lvert y_k - \rho_e(x_k, T; w)\bigr \rvert^2,
\]
under the constraints that $\sum_{n=1}^{N} w_n = 1$ and $w_n \geq 0$ for all $n \in \{1, \dotsc N\}$,
which guarantee that $\rho_0$ is a probability density.
Here $\rho_e(T, x; w)$ is the exact solution to the initial value problem~\eqref{eq:initial_value_problem_fokker_planck}
when the weights in the initial condition are given by $w$.
This admits the following explicit expression:
\begin{equation}
    \label{eq:exact_fokker_planck}
    \rho_e(x, t; w) = \sum_{n=1}^{N} w_n \frac{\displaystyle \exp \left( - \frac{\bigl\lvert x-m_n(t) \bigr\rvert^2}{2 \sigma_n(t)^2} \right)}{\sqrt{2 \pi \sigma_n(t)^2}},
\end{equation}
with $m_n(t) = \e^{-t} m_n$ and $\sigma_n(t) = \sqrt{1 + (\sigma_n^2 - 1) \e^{-2t}}$.
Note that, in this problem,
the forward model $w \mapsto \{\rho_e(x_k, T; w)\}_{k=1}^{K}$ is linear
and so,
in the absence of constraints,
the ensemble Kalman inversion method performs an exact preconditioned gradient descent.

\begin{figure}[ht]
    \centering
    \includegraphics[width=.7\linewidth]{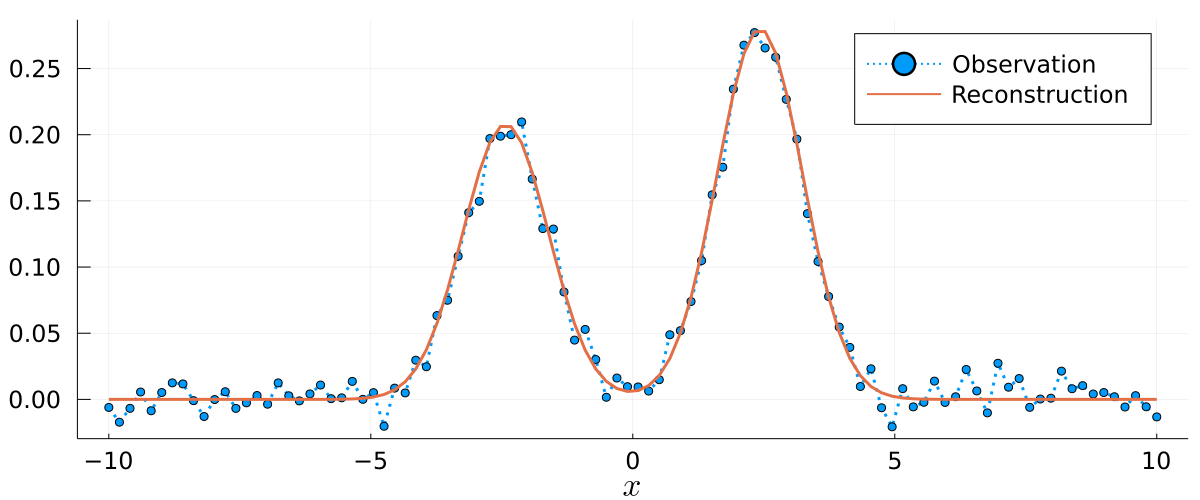}
    \caption{%
        Noisy observations of the solution to the Fokker--Planck equation,
        in the case of only 2 components in the binary mixture,
        and solution at time $T$ obtained by using the reconstructed initial condition.
    }%
    \label{fig:noisy_observation}
\end{figure}
We assume, for simplicity,
that the Gaussian mixture is composed of only $N=2$ components,
and that the known means and variances of these components are given by~$m_1 = - m_2 = 4$ and $\sigma_1 = \sigma_2 = .1$.
For the observations, we take the parameters $L = 10$, $\gamma = .01$ and $K = 100$,
and for the adaptation of the time step~\eqref{eq:time_step_adaptation} we use the parameters $\Delta t_* = 1$ and $\Delta t_{\max} = \infty$.
The ensemble size is taken to be $J=100$,
the final time is $T = .5$,
and the noisy observations are illustrated in~\cref{fig:noisy_observation}.

\begin{figure}[ht]
    \centering
    \includegraphics[width=0.24\linewidth]{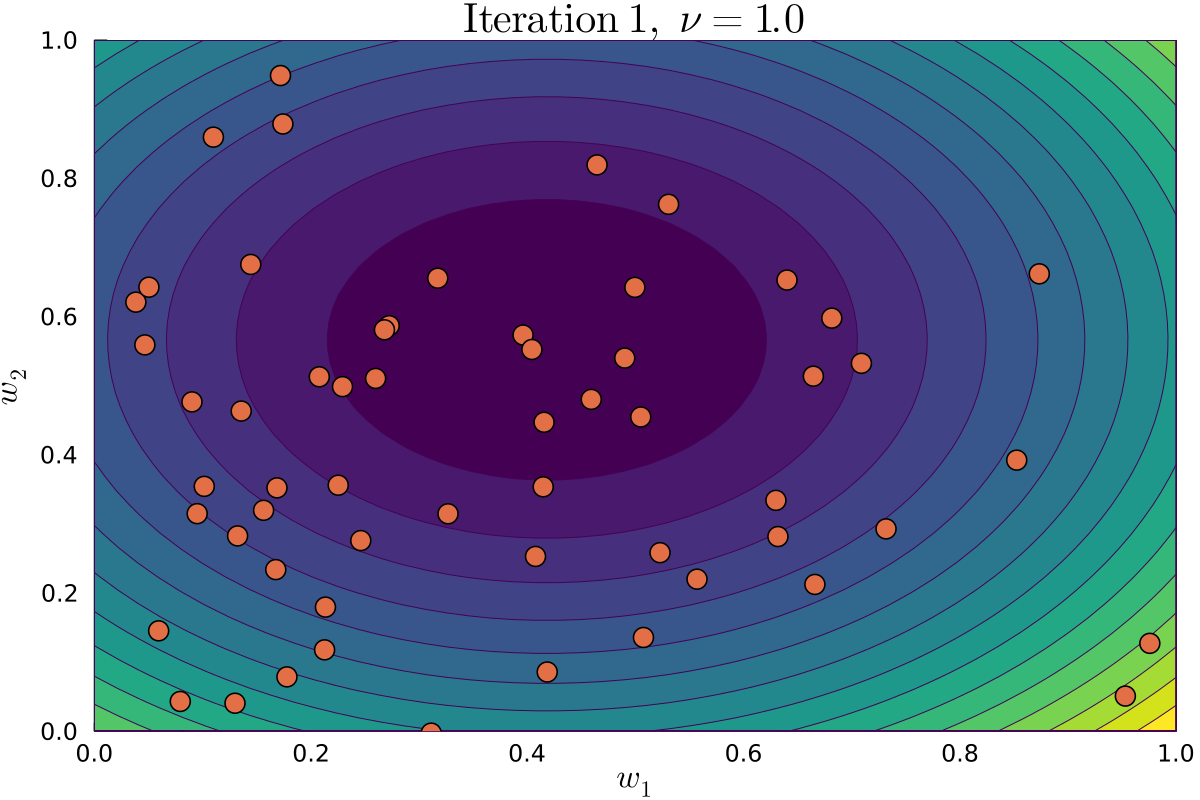}
    \includegraphics[width=0.24\linewidth]{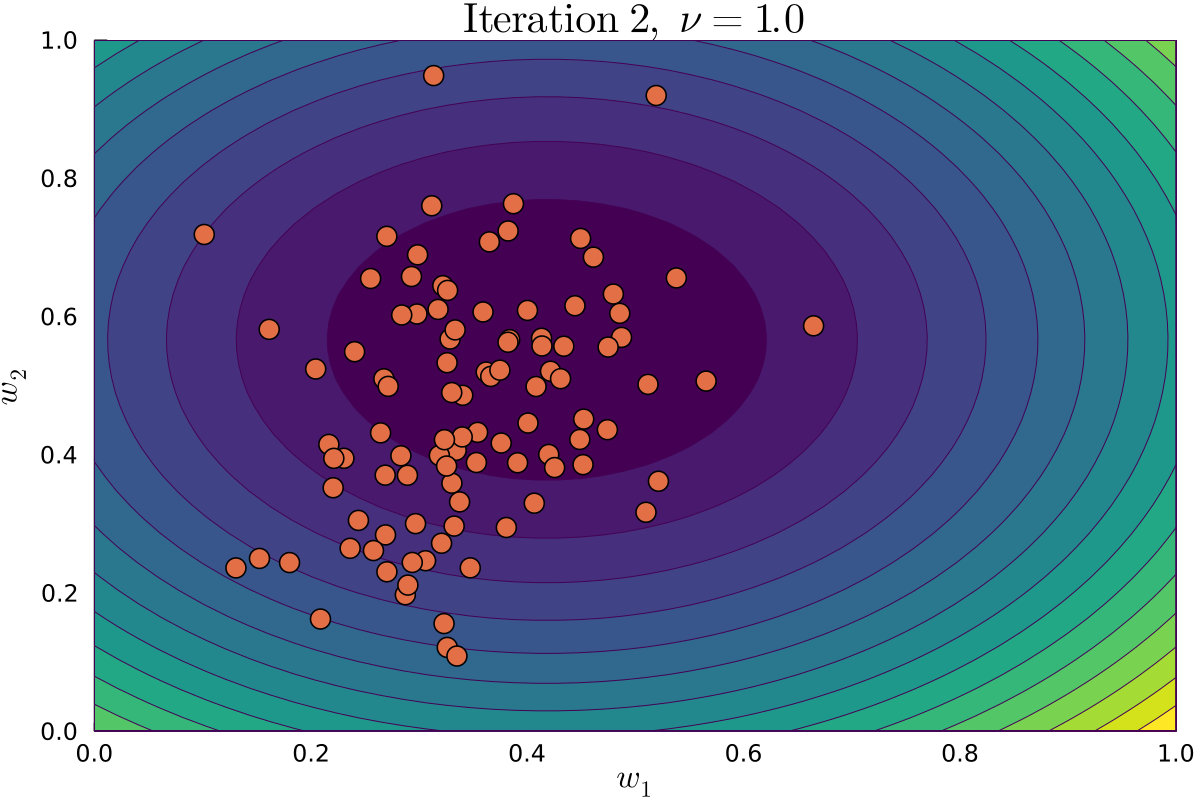}
    \includegraphics[width=0.24\linewidth]{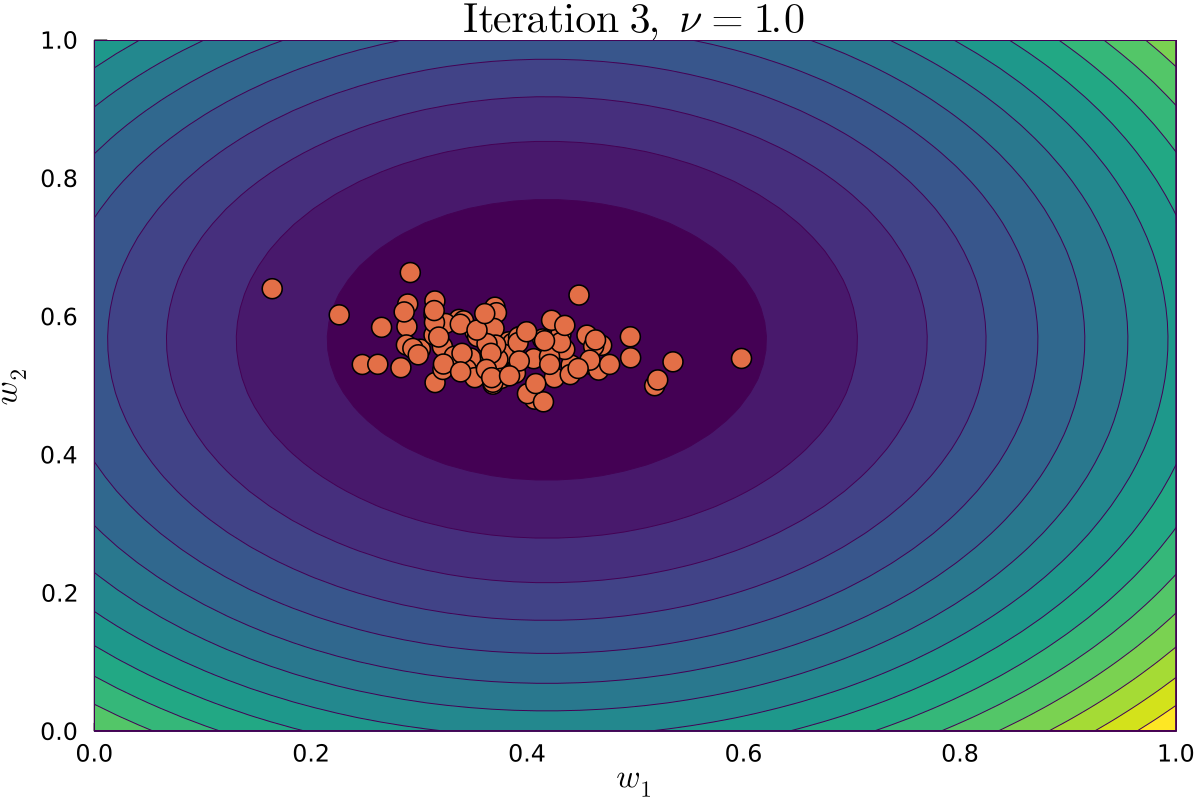}
    \includegraphics[width=0.24\linewidth]{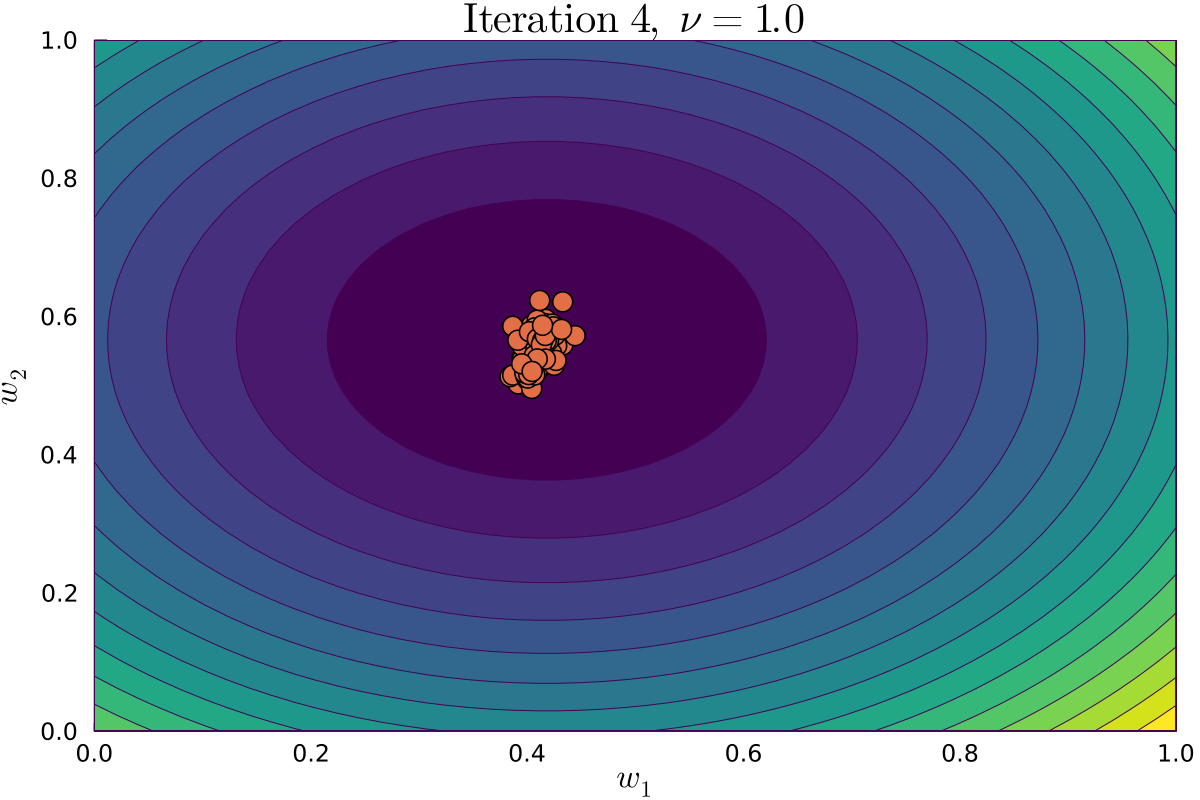}
    \includegraphics[width=0.24\linewidth]{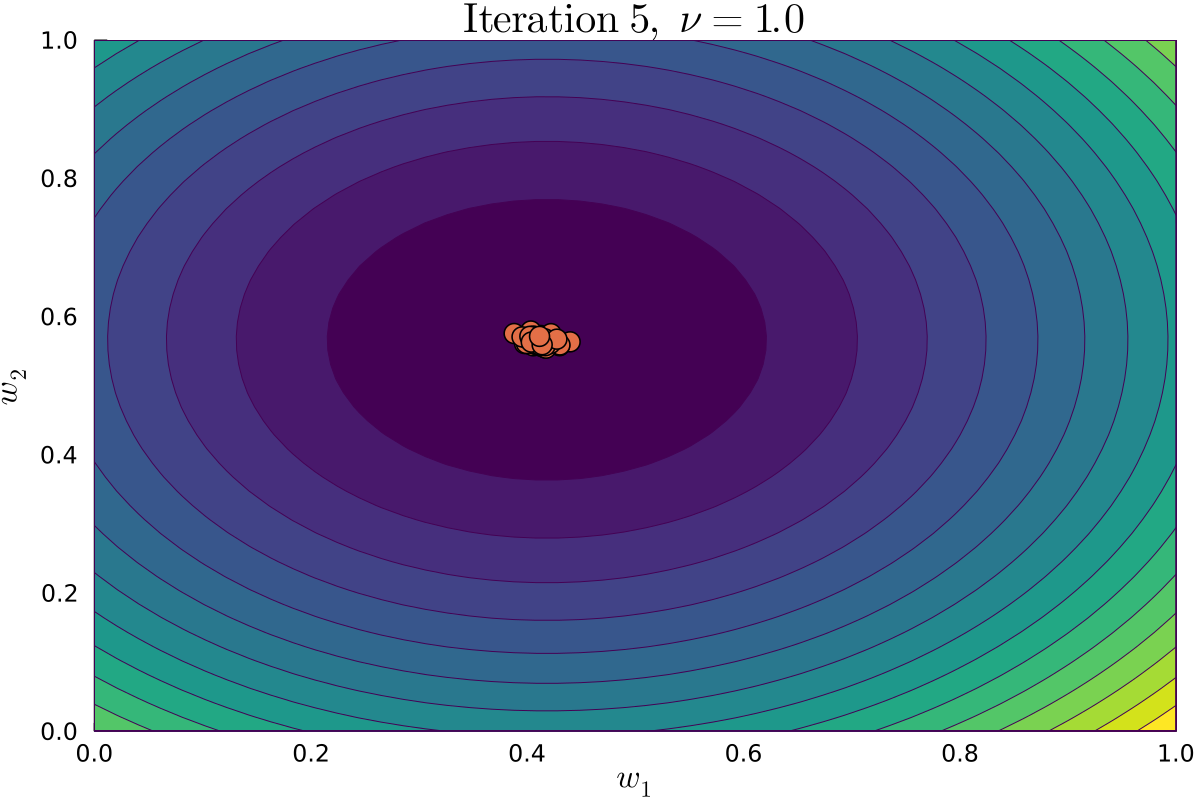}
    \includegraphics[width=0.24\linewidth]{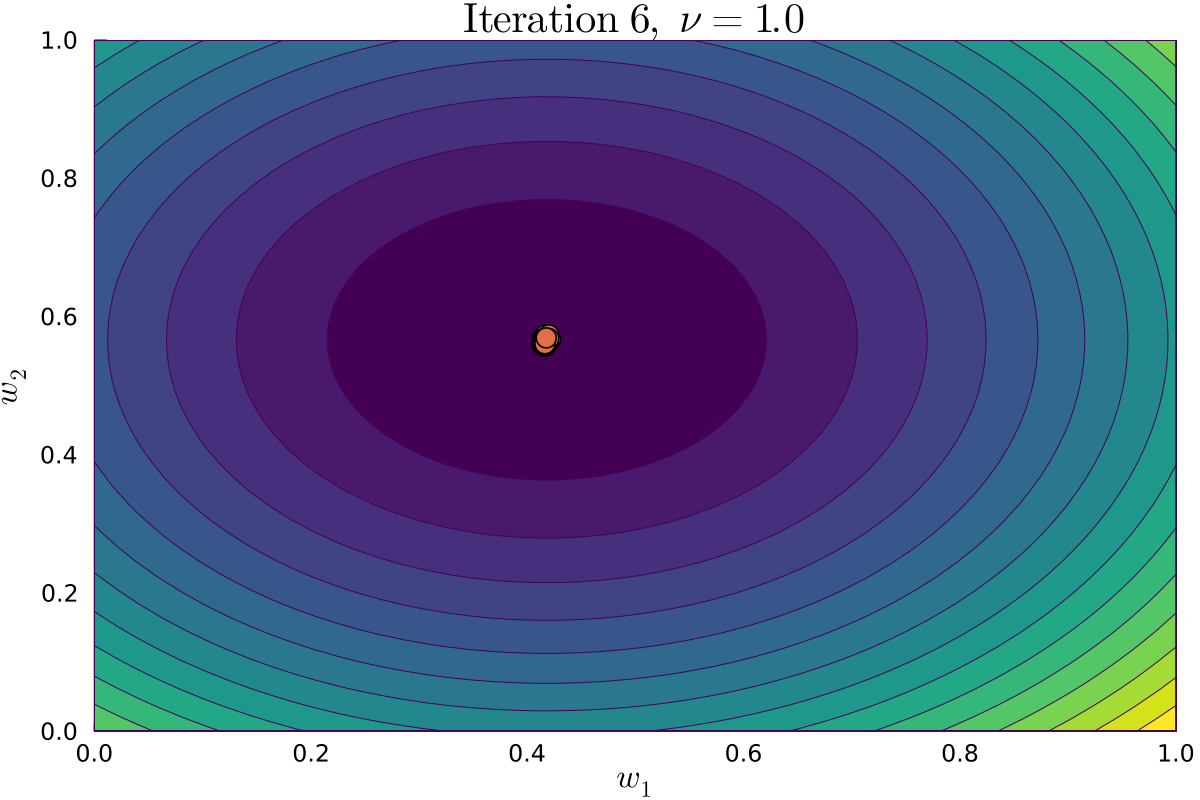}
    \includegraphics[width=0.24\linewidth]{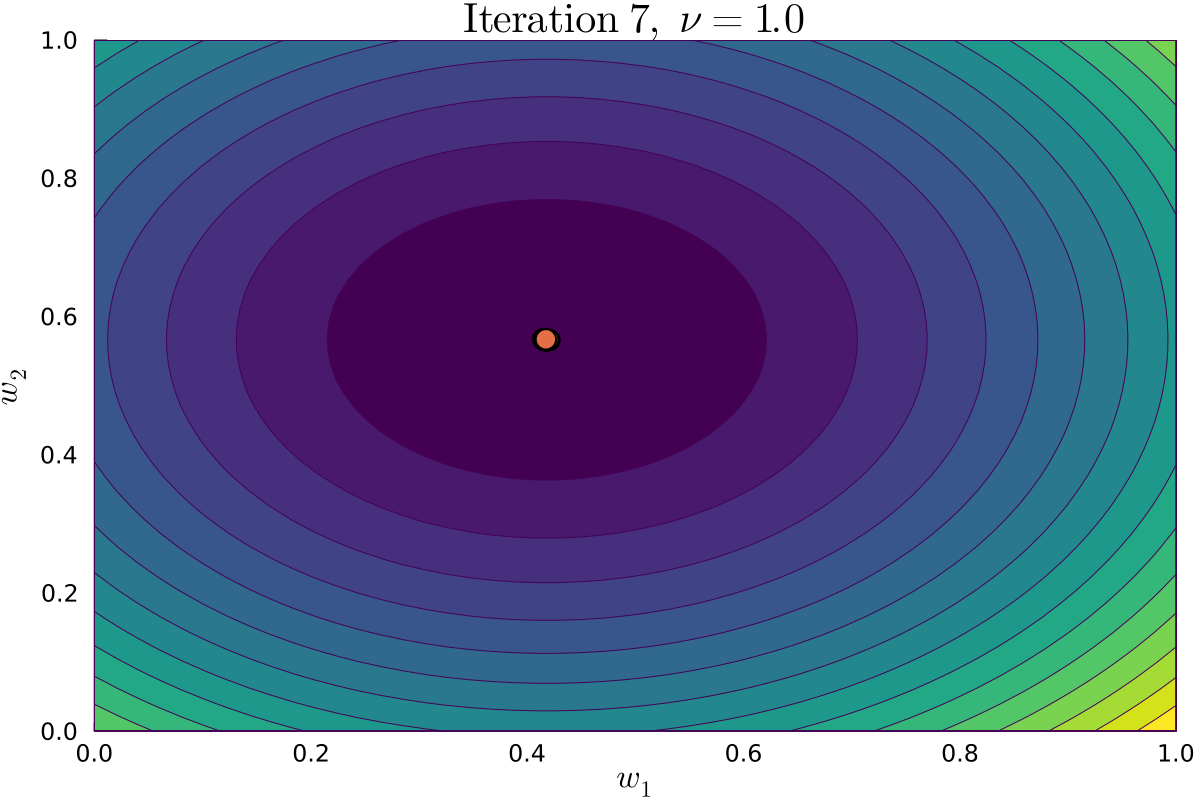}
    \includegraphics[width=0.24\linewidth]{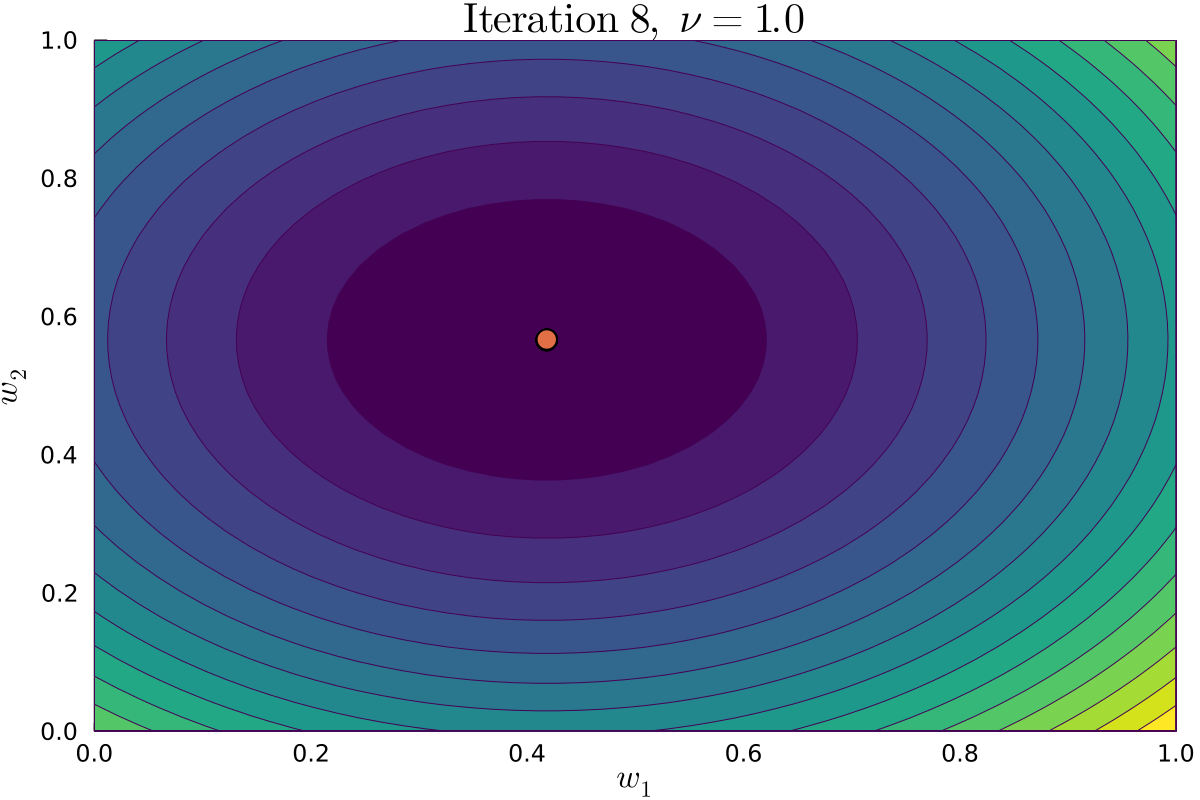}
    \caption{%
        Evolution of the ensemble obtained from the explicit discretization~\eqref{eq:good_explicit} of EKI with constraints,
        in the case where $\nu = 1$.
        The penalized objective function (denoted by $g$ in~\eqref{eq:weight_w_constraint})
        is depicted as a filled contour in the background.
    }%
    \label{fig:eki_with_large_epsilon}
\end{figure}

\begin{figure}[ht]
    \centering
    \includegraphics[width=0.24\linewidth]{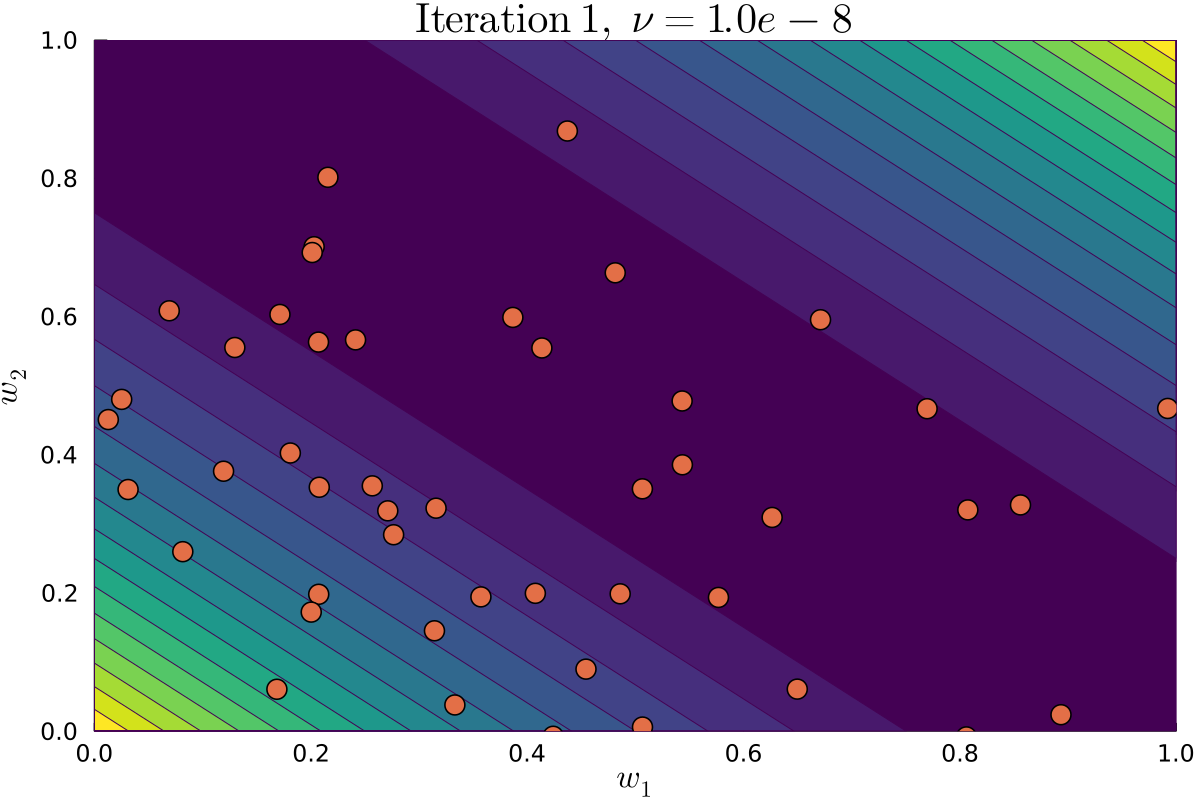}
    \includegraphics[width=0.24\linewidth]{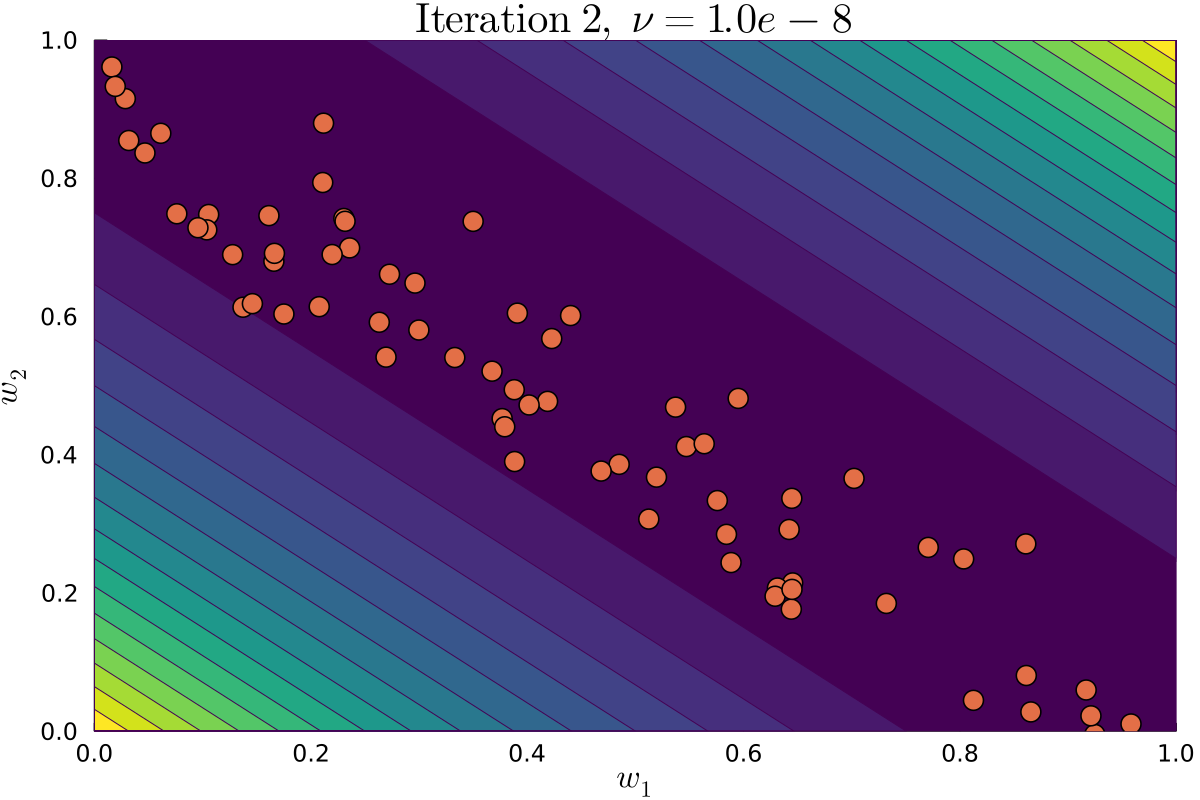}
    \includegraphics[width=0.24\linewidth]{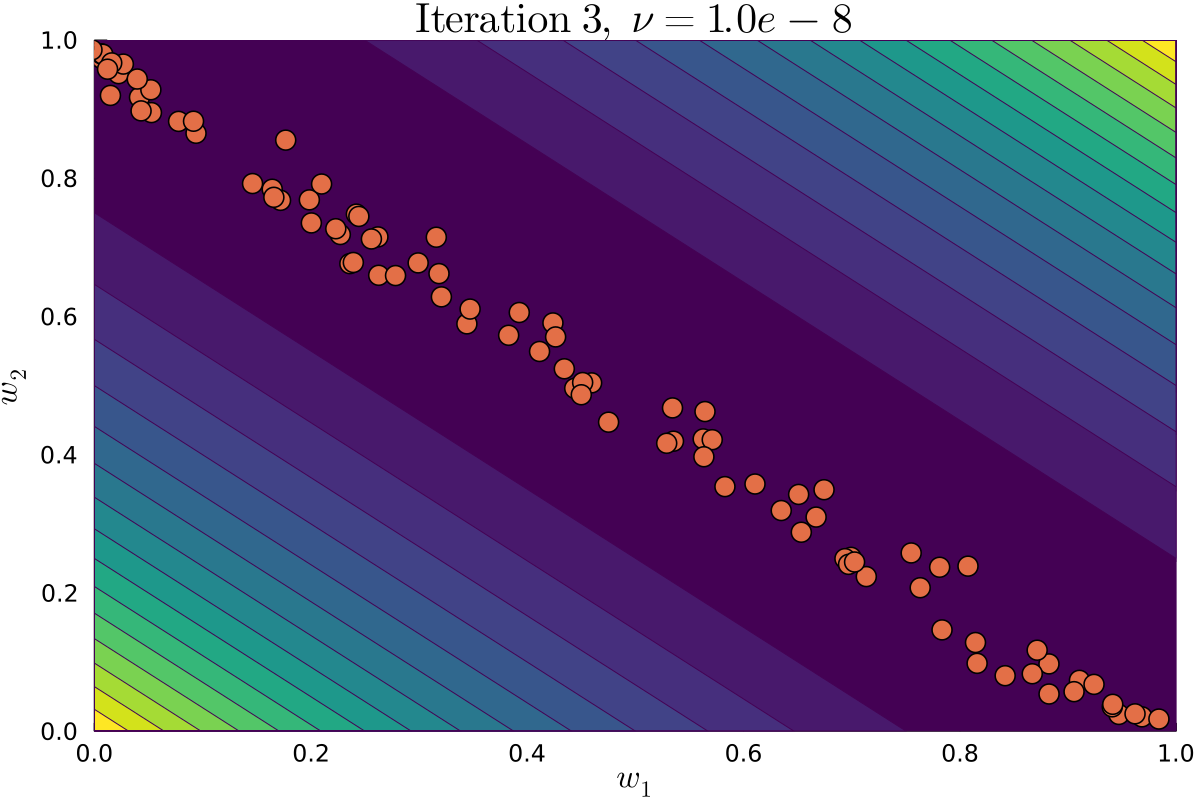}
    \includegraphics[width=0.24\linewidth]{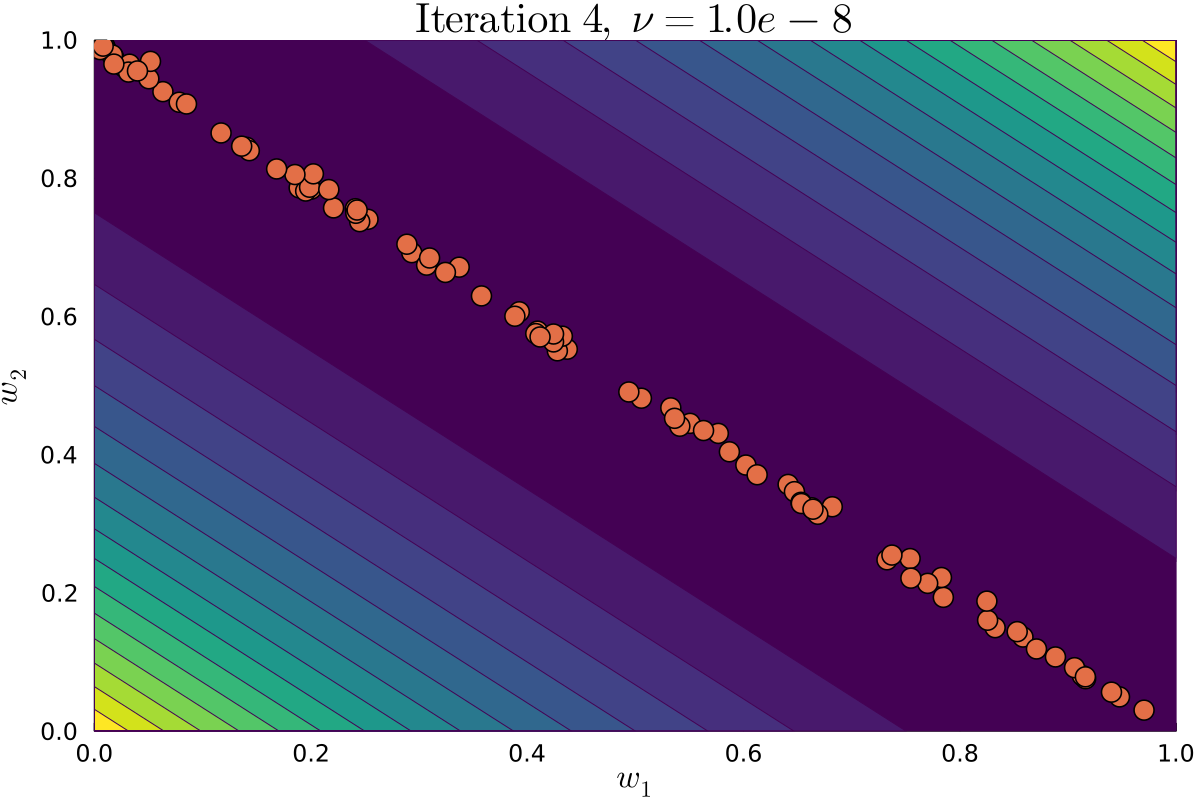}
    \includegraphics[width=0.24\linewidth]{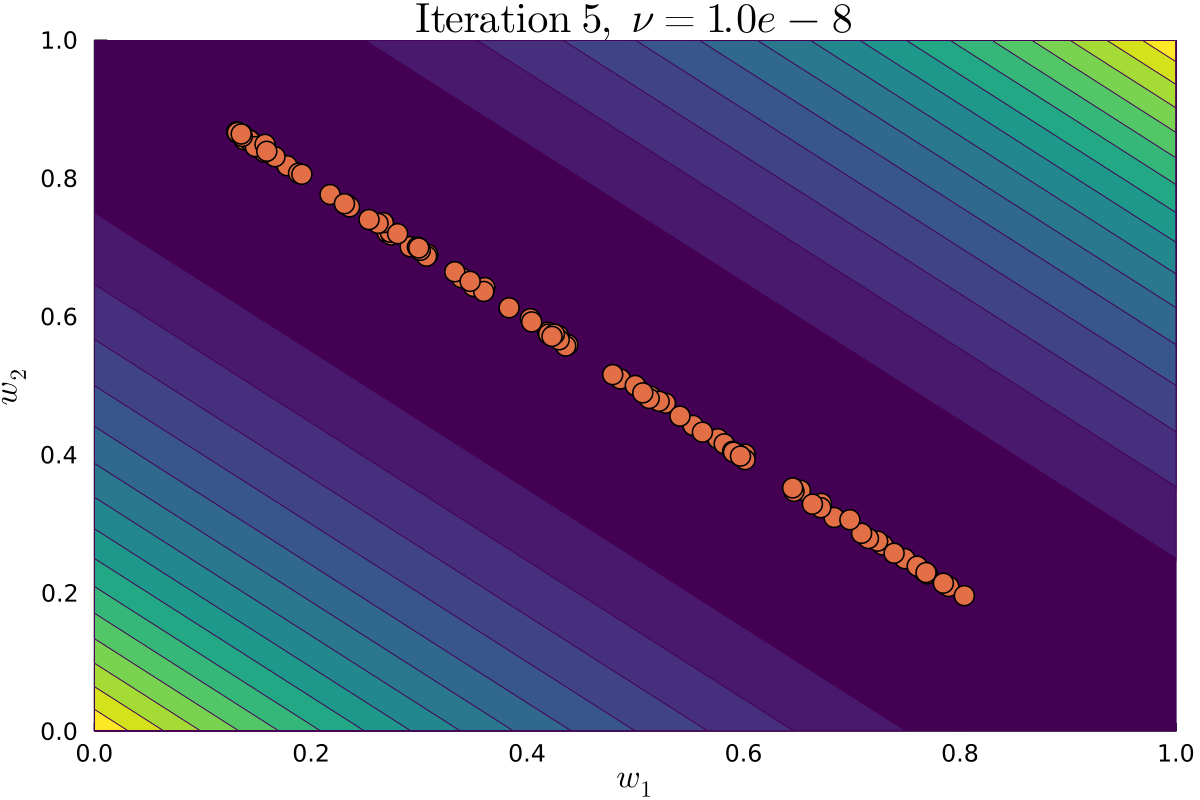}
    \includegraphics[width=0.24\linewidth]{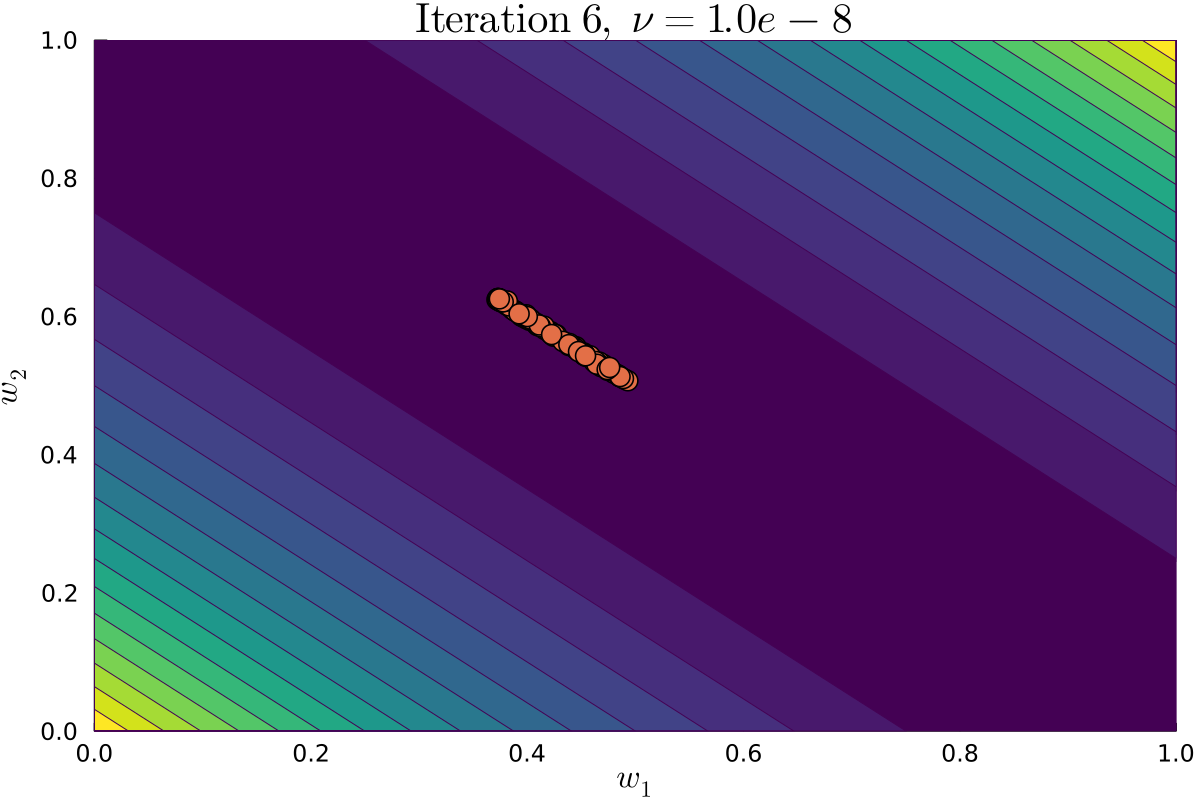}
    \includegraphics[width=0.24\linewidth]{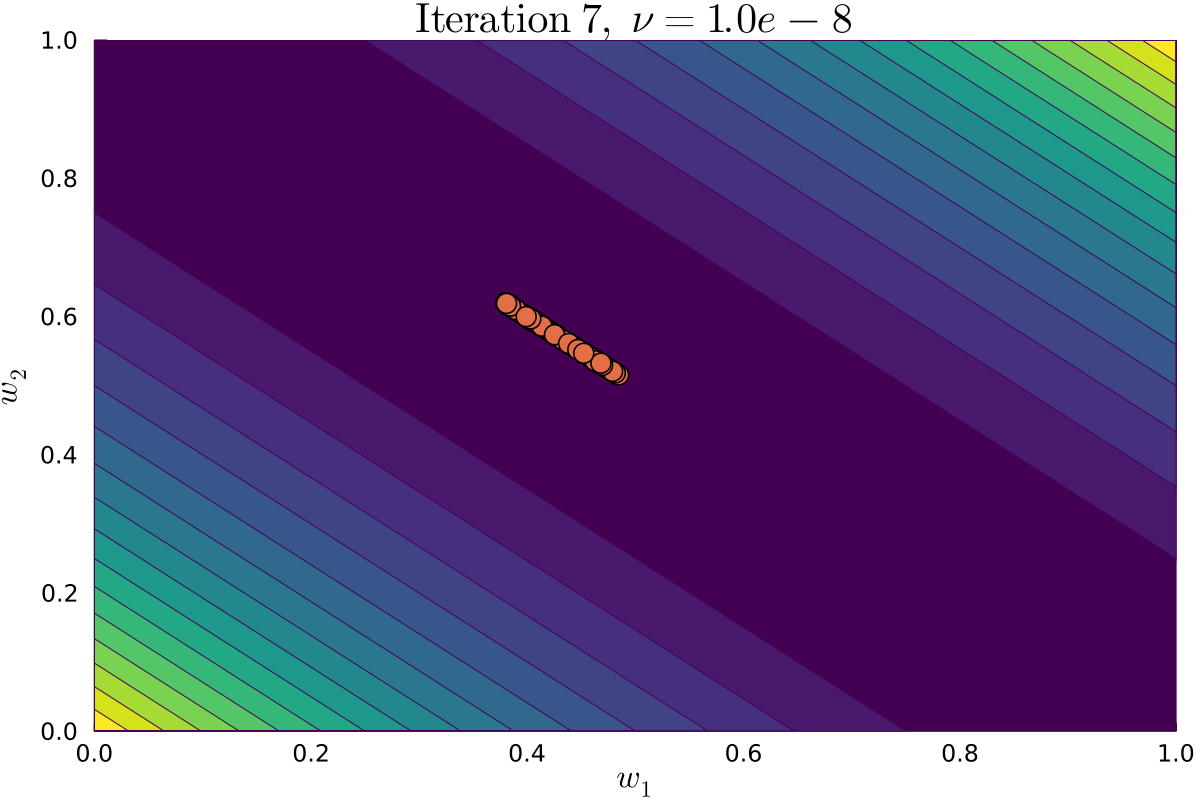}
    \includegraphics[width=0.24\linewidth]{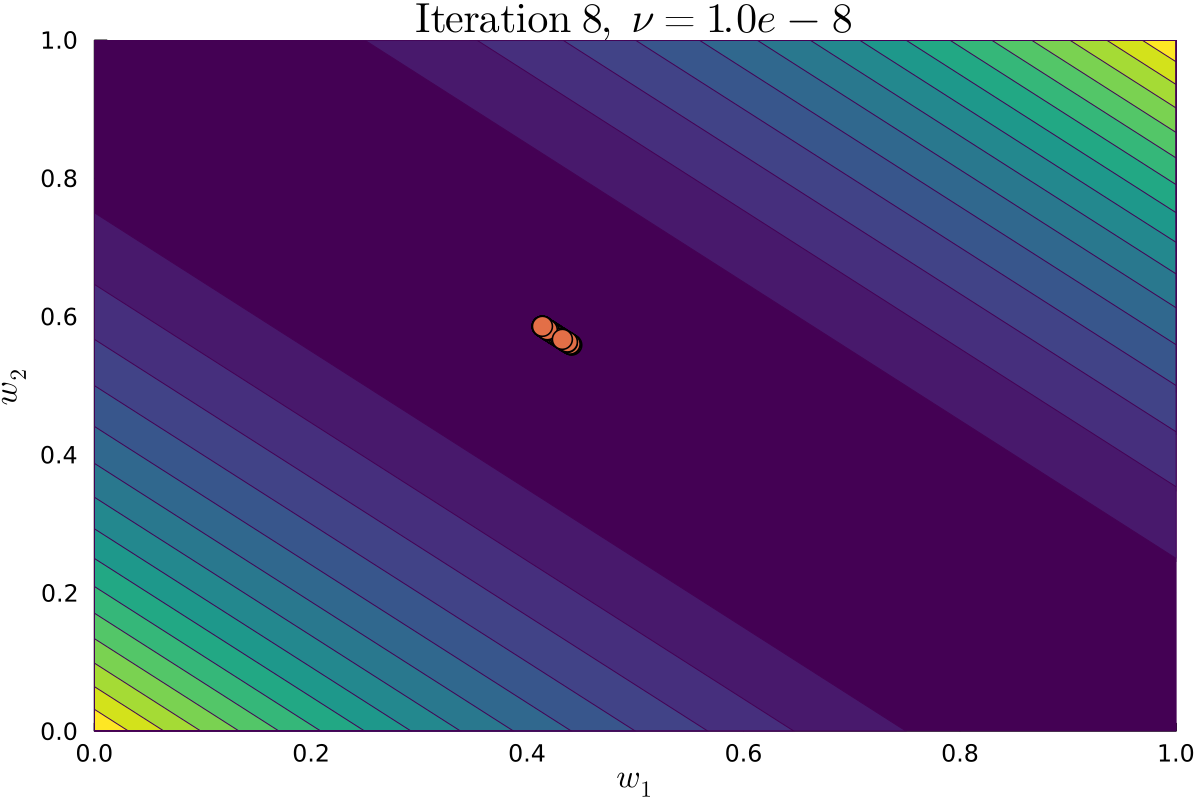}
    \caption{%
        Like~\cref{fig:eki_with_large_epsilon},
        this plot depicts the evolution the discrete-time dynamics~\eqref{eq:good_explicit},
        with now $\nu = 10^{-8}$.
    }%
    \label{fig:eki_with_small_epsilon}
\end{figure}

We begin by analyzing the influence of the small parameter $\nu$ on the solution returned by the method
when the particles forming the initial ensemble,
each corresponding to the couple of weights,
are drawn from $\mathcal N(0, I_2)$.
For this numerical experiment,
we use only the explicit discretization~\eqref{eq:good_explicit}.
Since the evolution~\eqref{eq:unified_eki} is deterministic,
and since we observed empirically that the initial ensemble does not have much influence on the point of convergence of the method,
we run only one simulation per value of $\nu$.
In the table below, we indicate for several values of $\nu$ the value of the sum $w_1 + w_2$ at the point of convergence of the method,
as well as the distance to the true value of the weights from which the observed data was generated.
This was taken to be $w^\dagger = (0.411..., 0.588...)$.
\begin{center}
    \begin{tabular}{|c|c|c|}
         \hline
         $\nu$ & $w_1 + w_2$ & $\lvert w_1 - w_1^\dagger \rvert + \lvert w_2 - w_2^\dagger \rvert$
         \\ \hline
         1 & 0.9845... & 0.0282... \\
         $10^{-2}$ & 0.9847... & 0.0282... \\
         $10^{-4}$ & 0.99289... & 0.0282... \\
         $10^{-6}$ & 0.999869... & 0.0282... \\
         $10^{-8}$ & 0.99999868... & 0.0282...
         \\ \hline
    \end{tabular}
\end{center}
As expected, we observe that the smaller $\nu$,
the closer the point of convergence is to the feasible manifold.
In addition, the proximity to the true weights is hardly affected by changes in the value of $\nu$.
\Cref{fig:eki_with_large_epsilon,fig:eki_with_small_epsilon} depict the evolution of the ensemble for $\nu = 1$ and $\nu = 10^{-8}$.
Although the method converges quickly in both cases,
the dynamics look qualitatively very different:
for $\nu = 1$ the effect of the penalization is not clearly apparent,
while for $\nu = 10^{-8}$ the constraint term appears dominant in the penalized objective function~\eqref{eq:weight_w_constraint}
depicted in the background.
In the latter case,
the particles first converge to a vicinity of the feasible line $w_1 + w_2 = 1$,
and then then move along this line to the optimizer.

\begin{figure}[ht]
    \centering
    \includegraphics[width=.45\linewidth]{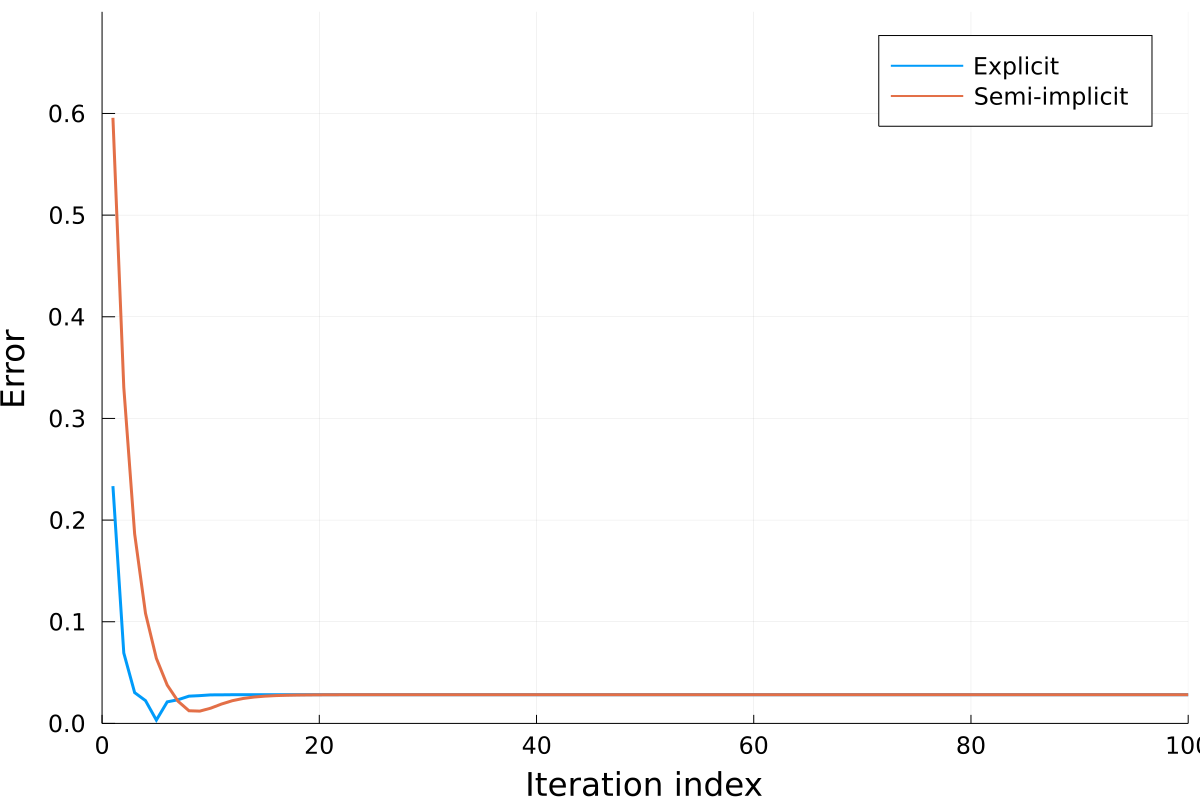}
    \includegraphics[width=.45\linewidth]{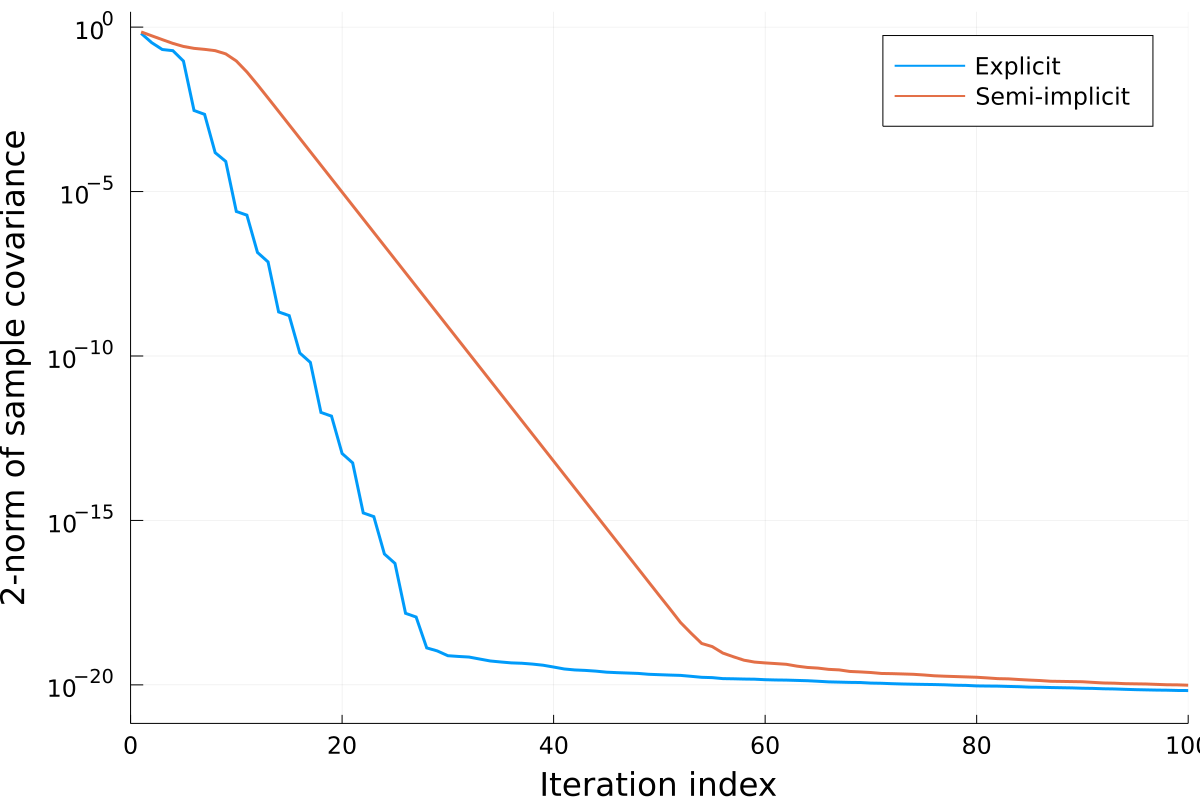}
    \caption{%
        Evolution of the error (\textbf{left}) and sample covariance (\textbf{right}),
        for the explicit and semi-implicit discretizations described in~\cref{sub:numerical_scheme_eki}.
        In both cases, the time step is adapted dynamically according to~\eqref{eq:time_step_adaptation}.
        We emphasize that the error is computed with respect not to the optimal solution to~\eqref{eq:optimization_problem_fp}
        but to the true weights,
        which were employed to generate the data.
        It is not surprising, therefore,
        that the error converges to a strictly positive value in the limit as the iteration index tends to infinity.
    }%
    \label{fig:comparison_time_integration}
\end{figure}

To conclude this paragraph,
we compare the two discretizations proposed in~\eqref{sub:numerical_scheme_eki},
when both methods use the same initial ensemble and with $\nu = 10^{-8}$.
In the left panel of~\cref{fig:comparison_time_integration},
we illustrate the evolution of the error,
measured as $\lvert \bar w_1 - w_1^\dagger \rvert + \lvert  \bar w_2 - w_2^\dagger  \rvert$,
where $\bar w_1$ are $\bar w_2$ are sample averages of the weights over the particles.
In the right panel,
we present the evolution of the matrix 2-norm of the sample covariance for each method.
Both methods converge quickly,
but the figures indicate a slightly faster convergence for the fully explicit discretization.

\paragraph{Recovering the weights and variances.}%
In this paragraph,
we consider the more challenging case where $N = 3$ and both the weights and variances of the initial Gaussian mixture need to be recovered,
in which case the forward model is no longer linear.
For simplicity, we assume that the means $\{m_n\}_{n=1}^N$ are still known and this time given by~$(m_1, m_2, m_3) = (-5, 0, 5)$,
and we seek to minimize
\begin{equation}
    \label{eq:optimization_problem_fp}
    f(w, v)
    = \sum_{k=1}^{K} \frac{1}{\gamma^2}
    \bigl\lvert y_k - \rho_e\bigl(x_k, T; w,  {\rm abs.} (v)\bigr) \bigr \rvert^2,
\end{equation}
where $w$ is the vector of weights and $v = (\sigma_1^2, \sigma_2^2, \sigma_3^2)$ is the vector of variances,
under the constraints
\[
    \sum_{n=1}^{3} w_n = 1,
    \qquad w_n \geq 0 \quad \forall n \in \{1, 2, 3\},
    \qquad \sigma_n^2 \geq 0 \quad \forall n \in \{1, 2, 3\}.
\]
Here the function ``${\rm abs.}$'' denotes the element-wise absolute value,
and the function $(x, t) \mapsto \rho_e(x, t; w, v)$ denotes the exact solution~\eqref{eq:exact_fokker_planck} to the Fokker--Planck equation~\eqref{eq:initial_value_problem_fokker_planck}
when the parameters $w$ and $v$ are employed in the initial condition~\eqref{eq:gaussian_mixture}.
We define the objective function in this manner,
with the presence of the ``abs.'' function,
in order to guarantee that
this function can be evaluated for any choice of parameters $w \in \real^3$ and $v \in \real^3$,
which is a requirement for applying EKI.

\begin{figure}[ht]
    \centering
    \includegraphics[width=0.8\linewidth]{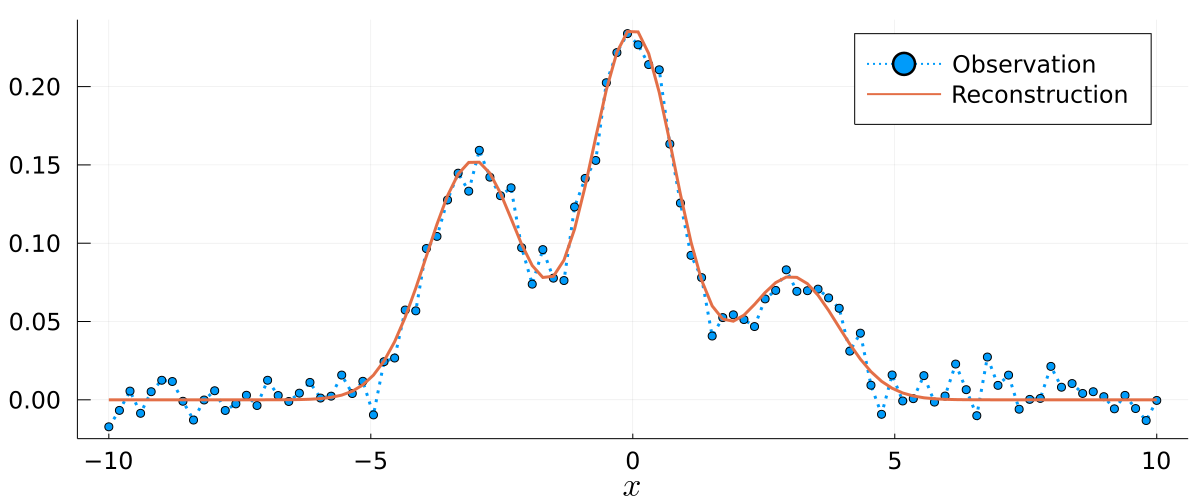}
    \caption{
        Noisy observations (blue dots) and solution to the Fokker--Planck equation~\eqref{eq:initial_value_problem_fokker_planck} at time $T$
        when using the reconstructed initial condition.
    }
    \label{fig:fokker-planck_small_eps}
\end{figure}

Apart from the parameters of the Gaussian mixtures,
and unless otherwise specified,
all the parameters employed to generate the numerical results presented in the rest of this section
are the same as in the previous paragraph.
We begin by examining the influence of the small parameter $\nu$ on the error and convergence point of the method.
We employ to this end the explicit method described in~\cref{sub:numerical_scheme_eki}.
The table below, in which all figures are truncated after three significant digits,
gives the points of convergence of EKI with constraint for different values of $\nu$.
We observe again that small values of $\nu$ lead to a point of convergence closer to the feasible manifold,
as is expected.
\begin{center}
    \begin{tabular}{|c|c|c|c|c|c|c|c|}
         \hline
         & $w_1 + w_2 + w_3$ & $w_1$ & $w_2$ & $w_3$ & $\sigma_1^2$ & $\sigma_2^2$ & $\sigma_3^2$ \\ \hline
        Truth & 1 & 0.333 & 0.476 & 0.191 & 0.400 & 0.100 & 0.500 \\ \hline
        $\nu = 1$ & 0.9844 & 0.337 & 0.480 & 0.167 & 0.433 & 0.095 & 0.303 \\
        $\nu = 10^{-2}$ & 0.9847 & 0.337 & 0.480 & 0.167 & 0.434 & 0.095 & 0.305 \\
        $\nu = 10^{-4}$ & 0.99502 & 0.342 & 0.482 & 0.171 & 0.467 & 0.093 & 0.363 \\
        $\nu = 10^{-6}$ & 0.9999273 & 0.344 & 0.483 & 0.173 & 0.483 & 0.092 & 0.392 \\
        $\nu = 10^{-8}$ & 0.999999270 & 0.344 & 0.483 & 0.173 & 0.483 & 0.092 & 0.393 \\
        \hline
    \end{tabular}
\end{center}
The solution to the Fokker--Planck equation at time $T$,
when using as parameters in the initial condition the point of convergence of EKI with small parameter~$\nu = 10^{-8}$,
is illustrated in \cref{fig:fokker-planck_small_eps}.
A good qualitative agreement between the observed and reconstructed solutions is observed.

\begin{figure}[ht]
    \centering
    \includegraphics[width=.45\linewidth]{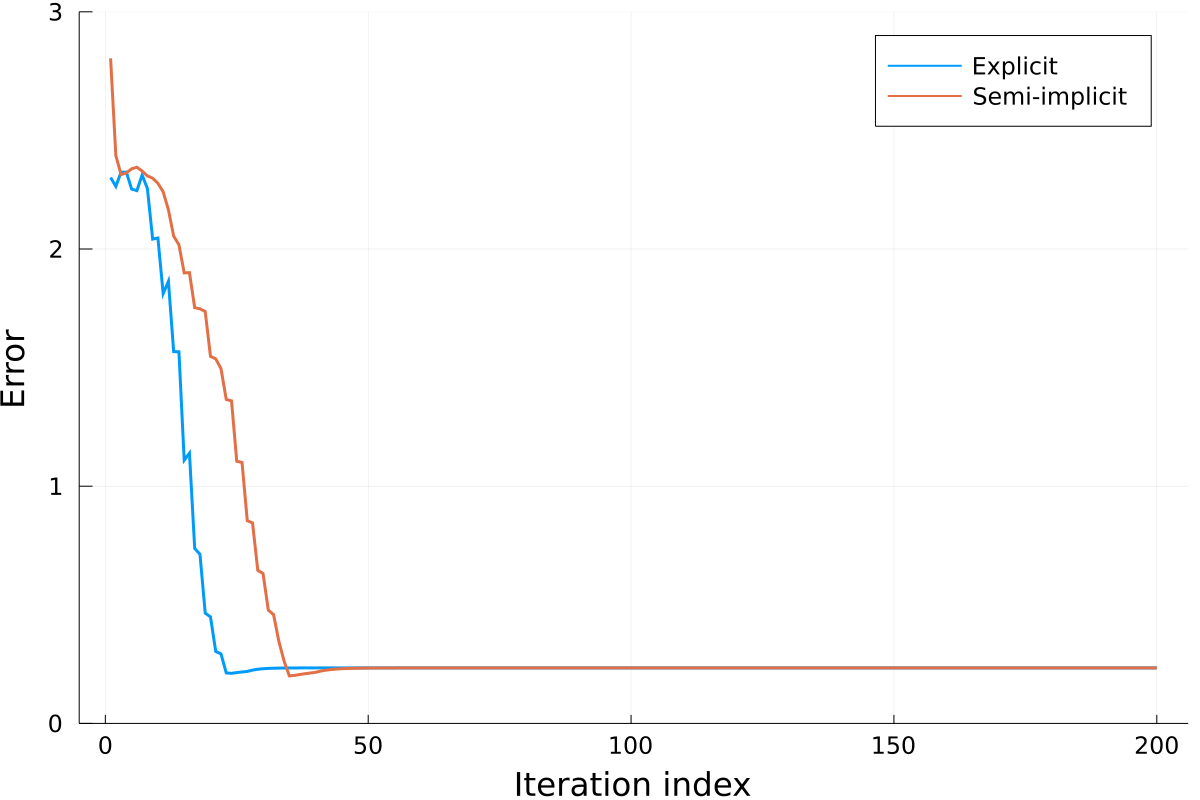}
    \includegraphics[width=.45\linewidth]{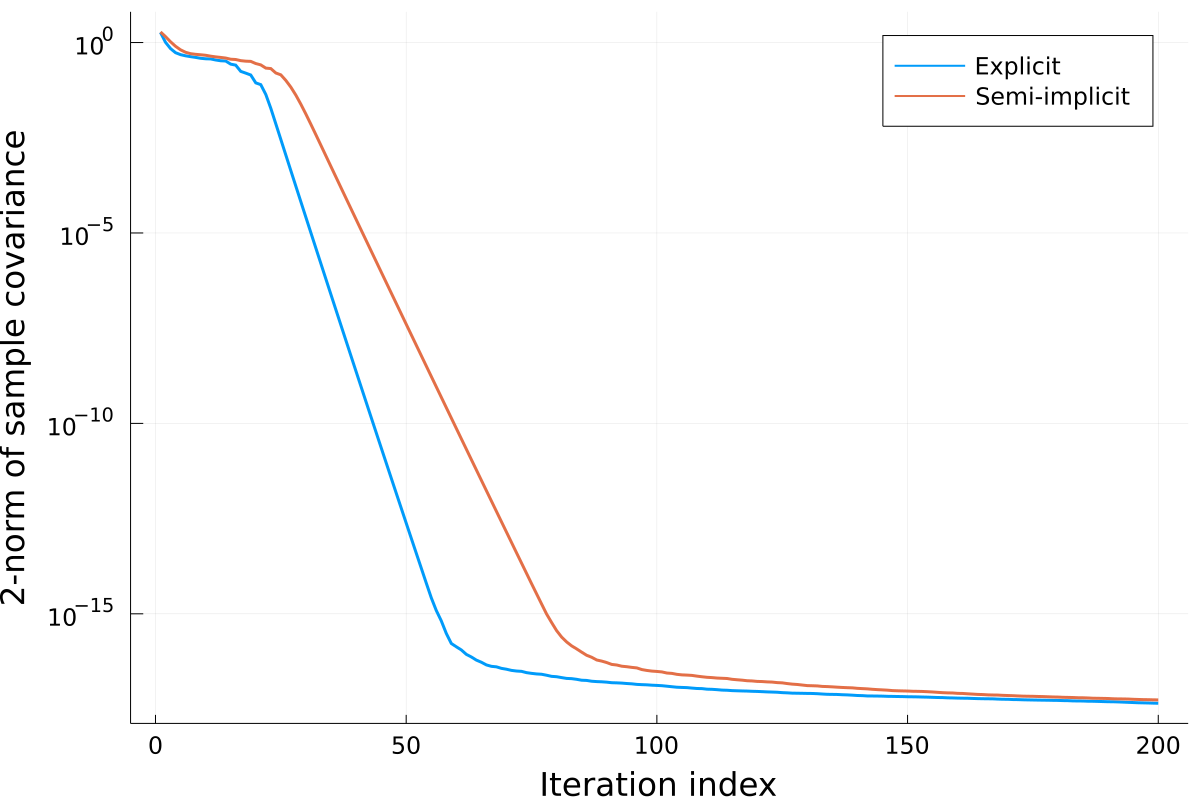}
    \caption{%
        Evolution of the error (\textbf{left}) and sample covariance (\textbf{right}),
        for the explicit and semi-implicit discretizations described in~\cref{sub:numerical_scheme_eki}.
        In both cases, the time step is adapted dynamically according to~\eqref{eq:time_step_adaptation}.
    }%
    \label{fig:comparison_time_integration_highd}
\end{figure}

To conclude this section,
we compare in~\cref{fig:comparison_time_integration_highd} the two discretization methods~\eqref{eq:good_explicit} and~\eqref{eq:good_semi_implicit} in the case where $\nu = 10^{-8}$.
All the ensemble members are initialized according to $\mathcal N(0, I_6)$.
The error depicted in the left panel is computed as
\[
    \sum_{n=1}^{3} \bigl\lvert \bar w_n - w_n^{\dagger} \bigr\rvert + \bigl\lvert \bar \sigma_n^2 - (\sigma_n^2)^\dagger \bigr\rvert,
\]
where $w_n^{\dagger}$ and $(\sigma_n^2)^\dagger$ are respectively the true weights and variances of the components of the initial Gaussian mixture,
whereas $\bar w_n$ and $\bar \sigma_n^2$ are the average weights and variances over the ensemble.
As already remarked in the caption of~\cref{fig:comparison_time_integration},
the error does not decrease to 0 as the number of iterations is increased,
because the true value of the mixture parameters do not coincide with the minimizer of $f(w, v)$ in~\eqref{eq:optimization_problem_fp}.
As shown in the right panel,
the sample variances decrease to a small value below $10^{-15}$ in less than a hundred iterations.

\section{Conclusions}%
\label{sec:conclusion}
In this work,
we study a new approach for incorporating constraints in consensus-based optimization and ensemble Kalman methods.
We demonstrate that, despite their simplicity and ease of implementation,
the proposed methods perform well for a number of toy examples.
Our study is mostly qualitative and aims at providing a general idea of the performance of the methods,
their behavior in the many-particle limit (in the case of CBO),
and the influence of the parameters they contain.

Future investigation would be useful in order to determine the practical value of the proposed approach for realistic high-dimensional constrained optimization problems,
and to further assess its efficiency in comparison with other methods in the literature.
It would also be worthwhile to obtain quantitative results on the link between the particle systems and their mean-field limits,
which could inform the choice of the number of particles employed in practice.

\section*{Acknowledgments}
JAC was supported by the Advanced Grant Nonlocal-CPD (Nonlocal PDEs for Complex Particle Dynamics: Phase Transitions, Patterns and Synchronization) of the European Research Council Executive Agency (ERC) under the European Union's Horizon 2020 research and innovation programme (grant agreement No. 883363) and by EPSRC grant number EP/T022132/1.
JAC and UV were also supported by EPSRC grant number EP/P031587/1.
CT was partly supported by the the European social Fund and by the Ministry Of Science, Research and the Arts Baden-W\"urttemberg (Germany).
UV was also supported by the Fondation Sciences Math\'ematiques de Paris (FSMP),
through a postdoctoral fellowship in the ``mathematical interactions'' program.

\printbibliography
\end{document}